\documentclass[11pt]{article}
\usepackage{amsmath}
\usepackage{amsthm}
\usepackage{amssymb}
\usepackage{graphicx}

\newtheorem{theorem}{Theorem}[section]
\newtheorem{lemma}[theorem]{Lemma}
\newtheorem{proposition}[theorem]{Proposition}
\newtheorem{corollary}[theorem]{Corollary}
\newtheorem{conjecture}[theorem]{Conjecture}

\theoremstyle{definition}
\newtheorem{definition}[theorem]{Definition}
\newtheorem{example}[theorem]{Example}

\newtheorem{remark}[theorem]{Remark}

\numberwithin{equation}{section}

\newcommand{\bq}{\ensuremath{\Bbb Q}}
\newcommand{\bz}{\ensuremath{\Bbb Z}}

\newcommand{\cala}{\ensuremath{\mathcal{A}}}
\newcommand{\calc}{\ensuremath{\mathcal{C}}}
\newcommand{\cald}{\ensuremath{\mathcal{D}}}
\newcommand{\cale}{\ensuremath{\mathcal{E}}}

\newcommand{\calg}{\ensuremath{\mathcal{G}}}
\newcommand{\calh}{\ensuremath{\mathcal{H}}}
\newcommand{\calj}{\ensuremath{\mathcal{J}}}
\newcommand{\calk}{\ensuremath{\mathcal{K}}}
\newcommand{\call}{\ensuremath{\mathcal{L}}}
\newcommand{\calw}{\ensuremath{\mathcal{L}}}
\newcommand{\calm}{\ensuremath{\mathcal{M}}}
\newcommand{\caln}{\ensuremath{\mathcal{N}}}
\newcommand{\calo}{\ensuremath{\mathcal{O}}}
\newcommand{\calq}{\ensuremath{\mathcal{Q}}}
\newcommand{\calr}{\ensuremath{\mathcal{R}}}
\newcommand{\cals}{\ensuremath{\mathcal{S}}}

\newcommand{\calv}{\ensuremath{\mathcal{V}}}

\newcommand{\fts}{finite type invariants}
\newcommand{\fti}{finite type invariant}
\newcommand{\fty}{finite type}

\newcommand{\betti}{\ensuremath{\mathfrak b}_1}
\newcommand{\bettip}{\ensuremath{\mathfrak b}_p}
\newcommand{\de}{\ensuremath{\delta}}
\newcommand{\dep}{\ensuremath{\mathrm{d}}}
\newcommand{\dpp}{\ensuremath{\mathrm{d}_p}}
\newcommand{\e}{\ensuremath{\varepsilon}}
\newcommand{\eye}{\ensuremath{\mathrm{I}}}
\newcommand{\g}{\ensuremath{\gamma}}
\newcommand{\ho}{\ensuremath{\mathrm{Hom}}}
\newcommand{\id}{\ensuremath{\mathrm{id}}}
\newcommand{\intt}{\ensuremath{\mathrm{int}}}
\newcommand{\la}{\ensuremath{\lambda}}

\newcommand{\lra}{\ensuremath{\longrightarrow}}
\newcommand{\m}{\ensuremath{3}}
\newcommand{\opp}{\ensuremath{\mathfrak o}_p}
\newcommand{\ovmu}{\ensuremath{\overline\mu}}
\newcommand{\ov}[1]{\ensuremath{\overline#1}}
\newcommand{\p}{\ensuremath{\partial}}

\newcommand{\rk}{\ensuremath{\mathrm{rk}}}
\newcommand{\sll}{\ensuremath{\mathrm{sl}}}
\newcommand{\sbq}{\ensuremath{\subseteq}}
\newcommand{\spq}{\ensuremath{\supseteq}}
\newcommand{\so}{\ensuremath{\mathrm{SO}}}
\newcommand{\spin}{\ensuremath{\mathrm{Spin}}}
\newcommand{\su}{\ensuremath{\mathrm{SU}}}
\newcommand{\tl}[1]{\ensuremath{\tilde#1}}
\newcommand{\tor}{\ensuremath{\mathrm{Tor}}}
\newcommand{\val}{\ensuremath{\mathrm{v}}}
\newcommand{\what}[1]{\ensuremath{\widehat#1}}
\newcommand{\wt}[1]{\ensuremath{\widetilde#1}}
\newcommand{\x}{\ensuremath{\times}}

\newcommand{\foot}{\setcounter{footnote}{1}\footnote}

\def\place#1#2#3{\text{\kern-#1pt{\raise#2pt\hbox{\rlap{#3}}}\kern#1pt}}

\def\figure#1#2{\includegraphics[scale=#2]{#1.eps}}

\begin{document}
\baselineskip=16 pt plus 1pt minus .5pt

\title{Finite type invariants of $3$-manifolds}
\author{Tim D.\ Cochran\footnote{partially supported
by the National Science Foundation DMS--9400224}${\ }^{\dagger}$ \and
Paul Melvin\footnote{both authors gratefully acknowledge the support
of Research Professorships of the Mathematical Sciences Research
Institute, Berkeley, California}}
\date{This Version: 8/1/99\quad
First version: 11/15/97}
\maketitle

\tableofcontents

\section{Introduction}

The primary objective of this paper is to propose a theory of
invariants of {\it \fty\/} for arbitrary compact oriented
\m-manifolds. We shall also give many examples of such invariants,
including some ``new'' \m-manifold invariants, and investigate the
algebraic and combinatorial structure of the set of all \fts.

At the most naive level, invariants of \fty\ should be thought
of as the {\sl polynomials} among all invariants. As such, they
should be computable (at least in theory) in polynomial time in the
complexity of the objects being studied. In recent years, a number of
different theories of \fts\ have evolved in a variety of topological
settings, with their origins in fields as diverse as singularity
theory and perturbative Chern-Simons theory. Perhaps the best known
of these is the theory for knots in the \m-sphere, which was
initiated by V.~Vassiliev \cite{V} and M.~Gusarov \cite{Gu}, and
developed by many other authors (in particular see \cite{BL}
\cite{BN} and \cite{Ko}). Importing some of the key notions from this
theory, T.~Ohtsuki \cite{O} developed an analogous theory for homology
\m-spheres which has been further studied by S.~Garoufalidis,
M.~Greenwood, N. Habegger, A.~Kricker, T.~Le, J. Levine, X.S.~Lin,
H.~Murakami, J.~Murakami, L.~Rozansky, B.~Spence, E.~Witten, and
others (see references). An extension to rational homology \m-spheres
was proposed by Garoufalidis and Ohtsuki \cite{GO1} (see \S10 for a
discussion of an apparent flaw in this theory).  Attempts to extend
beyond the set of rational homology spheres, however, have failed.
Indeed several authors have proved {\sl non-existence\/} theorems for
such extensions \cite{GO1} \cite{H1}. Moreover the most celebrated
extensions of specific \fts\ for rational homology spheres, namely C.
Lescop's extension of the Casson-Walker invariant and the
``universal'' \fti\ of Le-Murakami-Ohtsuki, vanish identically for
manifolds $M$ with first betti number $\betti(M)$ greater than three
\cite{Les} \cite{LMO} \cite{H2}. Our work seems to overcome these
difficulties.

The theory proposed here extends Ohtsuki's theory for integral homology
spheres, and is highly non-trivial for \m-manifolds of arbitrarily
large betti number.  Indeed much of the complexity of Ohtsuki's theory
embeds in our theory for manifolds of high betti number. It is shown
here that the coefficients of the Conway polynomial of a manifold with
first betti number one, as well as coefficients of the
Witten-Reshetikhin-Turaev quantum invariants for a general \m-manifold,
are of \fty. This provides evidence that the theory is a rich one.

There were several principles that guided us in formulating our
theory: 

\smallskip

1) ({\sl polynomial nature}) An invariant of \fty\ should be
a polynomial in some natural sense, preferably defined --- as in
Vassiliev's original viewpoint for knots --- as a function with
vanishing derivative of some order on a stratified space $X$. 
The``chambers'' of $X$ (components of the non-singular part) should
correspond to \m-manifolds, and the ``walls'' between chambers
correspond to certain singularities, perhaps singular \m-manifolds,
representing elementary transitions from one \m-manifold to another.
Some interesting work from this viewpoint has been done by N.~Shirokova
\cite{Sh}. 

\smallskip

2) ({\sl finiteness}) The set of all \fts\ should have an
algebraic structure, graded by degree, which when properly
interpreted is finite dimensional in each degree.

\smallskip

3) ({\sl non-triviality}) There should exist many independent invariants
in all degrees, including at least the more robust algebraic
topological invariants coming from (co)homology theory.  

\smallskip

4) ({\sl combinatorics}) There should be a combinatorial model for the
set of all \fts, as there is for knots and links \cite{Ko} and homology
spheres \cite{GO1} \cite{Le}. 

\smallskip

We begin with a heuristic definition of finite type invariants in
which their polynomial nature is evident. This requires the notion
of a ``combinatorial tangent bundle" for the set $\cals$ of
\m-manifolds. This point of view will also make it clear how our
definition differs from some previous attempts.

For motivation, first reconsider Ohtsuki's notion of \fts\ for
homology \m-spheres from this point of view. The basic idea is that
the homology spheres which are to be viewed as ``closest'' to $S^3$,
say, are those which are obtained from $S^3$ by $\pm1$ surgery on a
knot in $S^3$, denoted $S^3_K$. To this end, construct a cubical complex
$X(S^3)$ whose vertices are (oriented homeomorphism classes of) oriented
homology spheres $\Sigma$ and whose edges represent ``elemental
cobordisms'' between $\Sigma$ and $\Sigma_K$ (the result of surgery on
$K$ in $\Sigma$), i.e. $\Sigma\times I$ with a $2$-handle attached along a
$+1$ (or $-1$) framed knot $K$ in $\Sigma$. The edges emanating from
$\Sigma$ are the ``tangent vectors" at $\Sigma$ to the set of all
homology spheres. They are parametrized by $\pm1$-framed knots $K$ in
$\Sigma$. For $n>1$, the $n$-dimensional cubes are parametrized by
$\pm1$-framed $n$-component links $L$ in $\Sigma$ which have zero linking
numbers. Note that $X$ is connected. If $\phi$ is an invariant of homology
spheres then the (combinatorial) derivative of $\phi$ at $\Sigma$, in the
direction of $K$, is $\p_K\phi = \phi(\Sigma_K) - \phi(\Sigma)$. If
two such framed knots $\{K_1,K_2\}$ are disjoint and have linking
number zero in $\Sigma$, then one defines the second derivative at
$\Sigma$, $\p_{K_2}\p_{K_1}\phi = \phi(\Sigma_{K_1\cup K_2}) -
\phi(\Sigma_{K_1}) - \phi(\Sigma_{K_2}) + \phi(\Sigma)$, etc.. Given
this notion of the tangent space and given this combinatorial
derivative, Ohtsuki's \fty\ invariants of degree $n$ (for homology
$3$-spheres) are precisely the $n^{\rm th}$ degree polynomials. For
example, a degree zero invariant must have vanishing first
derivative, that is $\phi(\Sigma)=\phi(\Sigma_K)$ for each $\Sigma$
and $K$, and so is constant.

Now in extending this definition to all closed \m-manifolds the
crucial question is what should be the ``tangent vectors" to $\cals$
i.e.\ what are the allowable ``infinitessimal deformations''? In
brief, previous attempts allowed $0$-surgery on a knot in $M$ as a
deformation, and we do not. Clearly allowing more tangent vectors
imposes more conditions and increases the chances that the theory
becomes vacuous. For our theory, an admissible ``infinitessimal
deformation'' of $M$ is $M_K$ where $K$ is a $\pm1$ framed {\sl
null-homologous\/} knot in $M$. This corresponds to a cubical complex
$X$ which is disconnected, where a single path component has as
vertices all those \m-manifolds which can be obtained (one from
another) by a sequence of such ``deformations''. In particular all
such \m-manifolds have isomorphic homology groups. The component
containing $S^3$ is $X(S^3)$ as above. Once having stipulated this set
of deformations, we define a {\it polynomial invariant\/} of degree
at most $n$ to be one whose $(n+1)$-st order mixed partial
derivatives vanish. The mixed partial is defined only in restricted
cases as above. We shall not make this precise. The reader can
extract it from our precise definition of \fty\ which follows below.
But, in summary, there is a natural sense in which our \fts\ are
polynomials, and there is a space $X$ whose vertices (chambers) are
\m-manifolds and whose edges (walls between chambers) are elementary
cobordisms (``singular \m-manifolds''), as in the approach of
Vassiliev.

We shall now give our definition for \m-manifolds, which can be seen
to be formally identical to that of Ohtsuki for homology \m-spheres,
and then discuss the elements of the definition which distinguish it
from other attempts. In section~9 we give several significant
generalizations of our definition.

Let $\cals$ be a set of equivalence classes of \m-manifolds
$(M,\sigma)$ with some additional ``structure'' $\sigma$, modulo
``structure-preserving'' homeomorphisms. Examples of the structures
which may be considered are: orientation, spin structure, a marking of
$\p M$ (i.e. a homeomorphism from $\p M$ to a fixed abstract surface),
an element of $H^1(M;\bz_n)$, a marking of $H_1(M)$ (i.e.\ an
isomorphism from $H_1(M)$ to a fixed abstract abelian group). In
fact all of these theories are discussed herein, but a unified
definition is given below. The type of structure and the set $\cals$
may {\sl not} be chosen entirely arbitrarily; there is a mild
restriction discussed below. 

Let $\calm$ be the free abelian group on the set $\cals$. We define a
decreasing filtration of subgroups $\calm =
\calm_0\supset\calm_1\supset\calm_2\supset\cdots$ below,  and with
respect to this filtration and some fixed Noetherian ring $A$ we
stipulate:

\begin{definition} A function $\phi:\cals\to A$ is {\it \fty\/} of
{\it degree $\ell$\/} if its linear extension to $\calm$ vanishes on
$\calm_{\ell+1}$, but not identically on $\calm_\ell$. Let
$\calo^A_\ell$, or often merely $\calo_\ell$, denote the $A$-module
of all \fts\ of degree {\sl at most $\ell$}, i.e.\
$\ho(\calm/\calm_{\ell+1},A)$, and let $\calo$ denote the union of all
$\calo_\ell$.
\end{definition}

The filtration we use is defined as follows.

\begin{definition} The framed link $L=\{L_1,\dots,L_\ell\}$ in
$M$ is {\it admissible\/} if
\begin{enumerate}
\item[a)] each $L_i$ is null-homologous in $M$
\item[b)] the pairwise linking numbers of $L$ (measured in $M$) are
zero
\item[c)] the framings are $\pm1$ with respect to the longitude
guaranteed by (1).
\end{enumerate}
Such a link in $S^3$ has been called {\sl unit-framed,
algebraically split\/} by some other authors. Clearly any sublink of
an admissible link is itself admissible.
\end{definition}

If $L$ is a framed link in $M$ then $M_L$ will denote the result of
Dehn surgery on $M$ along $L$ \cite{Ro}. If $L$ is an admissible link
in $M$ then $[M,L]$ will denote the element of $\calm$ represented by
the (formal) alternating sum of manifolds $M_S$ over all sublinks $S$
of $L$ (including $S=\phi$ and $S=L$),
$$
[M,L] = \sum_{S<L}(-1)^sM_S.
$$
Here the number of components of a link ($S$ or $L$, for example) is
denoted by the corresponding lower case letter ($s$ or $\ell$). If $L$
is empty then $[M,L]$ is the class of $M$ itself. 

It is also sometimes convenient to use the notation $M_{\de L}$ for
$[M,L]$ where $\de$ is the operator which sends a framed link to the
alternating sum of its sublinks, 
$$
\delta L = \sum_{S<L}(-1)^sS.
$$  
Note that $\delta$ is an involution on the free abelian group $\call$
generated by framed links \cite{CM}.

\begin{definition} Let $\calm_\ell$ be the span of the set $\cals_\ell$
of all $[M,L]$, where $M$ is an element of $\cals$ and $L$ is an
admissible link of $\ell$ components in $M$.  As will be seen below,
this defines a filtration 
$$
\calm = \calm_0\supset\calm_1\supset\calm_2 \supset\cdots
$$  
with intersection $\calm_\infty = \bigcap^\infty_{\ell=0}\calm_\ell$. The
quotients $\calm_\ell/\calm_{\ell+1}$ will be denoted by $\calg_\ell$,
and so $\calg = \calg_0\oplus\calg_1\oplus\calg_2 \oplus\dots$ is the
associated graded group.
\end{definition}

One can think of $\cals_1$ as the set of unit tangent vectors to
$\cals$, of $\calm_1$ as the tangent bundle of $\cals$, and
inductively, of $\cals_{\ell+1}$ as the set of unit tangent vectors to
$\cals_\ell$ and $\calm_{\ell+1}$ its tangent bundle.

The reader should note that the definitions above are incomplete. If
$M$ is a manifold with structure $\sigma$ and $S$ is an admissible
link in $M$ then we must specify how the structure $\sigma$ is
``propagated'' to a structure $\sigma_S$ on $M_S$ in order that the
symbol $[M,L]$ be defined. This functor must be invariant under
structure-preserving homeomorphisms of the pair $(M,S)$. When the
structure is an orientation or a marking of $\p M$ then this
propagation is obvious, but when the structure is a spin structure or a
marking of $H_1$ then more must be said (later). This problem restricts
the type of structures which may be considered under this definition.
It is now evident that the set $\cals$ must have the following
closure property: if $(M,\sigma)\in\cals$ then, for any admissible link
$S$ in $M$, $(M_S,\sigma_S)\in\cals$. With these mild restrictions,
Definitions~1.1--1.3 suffice to define a theory of \fts\ for many
categories of \m-manifolds. For simplicity of exposition we shall
henceforth restrict attention to {\sl compact orientable \m-manifolds}
and to {\sl structures which include an orientation}.

The following combinatorial identity holds and shows immediately that
\linebreak$\calm_{\ell+1}\subset\calm_\ell$.

\begin{lemma} If $L\cup K$ is an admissible link in $M$ and $K$ is a
knot, then $L$ is admissible in $M_K$ and $[M,L\cup
K]=[M,L]-[M_K,L]$. More generally, if $K$ is a link then
$[M_K,L]=[M,L\cup\delta K]$ \ $($where the latter is defined linearly for
arguments in $\call)$.
\end{lemma}

\begin{proof} $[M_K,L]=M_{\delta L\cup K}=M_{\delta(L\cup\delta
K)}=[M,L\cup\delta K]$, since $\delta^2=\id$.
\end{proof}

Definition 1.1, when restricted to the subgroup of $\calm$ spanned
by the set of oriented homology \m-spheres is precisely that of
Ohtsuki. It differs from the definition of Garoufalidis-Ohtsuki on
the span of the set of rational homology \m-spheres
(\cite[Definition~1.2]{GO1}; see \S10).

In general the key difference in our proposed extension lies in the
definition of an admissible link. Note that if $L$ is admissible in
$M$ then $H_1(M_L)\cong H_1(M)$. Moreover if one considers the
cobordism $W$ from $M$ to $M_L$, given by attaching $2$-handles to
$M\x[0,1]$ along the components of $L$, then $H_1(M)\cong
H_1(W)\cong H_1(M_L)$. We say that $M_0$ and $M_1$ are
{\it $H_1$-bordant\/} if there exists an oriented cobordism between them
which is a product on $H_1$. Thus one sees that each term $M_S$ of
$[M,L]$ is $H_1$-bordant to $M$ and consequently the partition of
$\cals$ into $H_1$-bordism classes is respected by the filtration.
It follows that the study of invariants of \fty, in our sense,
largely reduces to the study of such on each fixed $H_1$-bordism
class. 

More precisely, for any fixed $3$-manifold $M$ let $\cals(M)$ denote the
set of all \m-manifolds $H_1$-bordant to $M$, and $\calm(M)$ denote its
span in $\calm$. For example $\calm(S^3)$ is precisely the group studied
by Ohtsuki. One sees that $\cals(M)$ satisfies the required closure
property.  

Now for each non-negative integer $\ell$, let
$\calm_\ell(M)$ be the subgroup of $\calm_\ell$ spanned by
all $[M',L]$ with $M'\in\cals(M)$.  Then by the above remark and Lemma
1.4, there is a decreasing filtration 
$$
\calm(M) = \calm_0(M)\supset\calm_1(M)\supset\calm_2(M) \supset\cdots
$$ 
and we can define a function $\phi:\cals(M)\to A$ to be \fty\ of degree
$\ell$ if its extension to $\calm_{\ell+1}(M)$ is zero and its extension
to $\calm_\ell(M)$ is not identically zero.  As above, set $\calg_\ell(M)
= \calm_\ell(M)/\calm_{\ell+1}(M)$, also denoted $(\calm_\ell /
\calm_{\ell+1})(M)$, and $\calo_\ell(M) =
\ho((\calm/\calm_{\ell+1})(M),A)$.  Then the following are trivial
consequences of the definitions.

\begin{proposition} Suppose $\calh$ is the set of $H_1$-bordism
classes of elements of $\cals$. Choose a representative $M_i$ for each
class $i\in\calh$.  Then for each $\ell\ge0$,
\begin{enumerate}
\item[\rm{a)}] $\calm=\bigoplus\limits_\calh\calm(M_i)$
\qquad\qquad {\rm b)} $\calm_\ell=\bigoplus\limits_\calh\calm_\ell(M_i)$
\item[\rm{c)}] $\calg_\ell=\bigoplus\limits_\calh\calg_\ell(M_i)$
\qquad\qquad\ \,{\rm d)}
$\calo_\ell\cong\prod\limits_\calh\calo_\ell(M_i)$
\end{enumerate}
\end{proposition}

\begin{proof} The partition of $\cals$ into
$H_1$-cobordism classes clearly induces a direct sum decomposition on
free abelian groups on the sets, establishing 1.5a. Since every
element in the sum $[M,L]$ is $H_1$-cobordant to $M$, 1.5b follows
easily. Then 1.5c is an easy algebraic consequence of 1.5b. Finally
$\calo_\ell=\ho(\calm/\calm_{\ell+1},A)\cong\Pi_\calh\ho
((\calm/\calm_{\ell+1})(M_i),A)\equiv\Pi_\calh\calo_\ell(M_i)$.
\end{proof}

The last isomorphism in Proposition 1.5 makes it clear that invariants of
\fty, in our sense, are constructed from invariants of \fty\ on each
$H_1$-bordism class. In fact the degree~$0$ \fts\ are precisely those
which are constant on $H_1$-bordism classes, i.e. the ``locally
constant'' functions on $\cals$. For example it is easy to see
that the function $\phi:\cals\to\bz$ given by the first betti number
is \fty\ of degree~$0$, being constant on each $\cals(M_i)$. Similarly
the function which assigns $|H_1(M)|$ to $M$ if $H_1(M)$ is finite,
and $0$ otherwise, is of degree zero.

Our point of view is that we have ``split'' the classification
problem for \m-manifolds into two parts. First, the problem of
determining if $M_0$ and $M_1$ lie in the same $H_1$-bordism class.
Second, if they lie in the same $H_1$-bordism class, can they be
distinguished by invariants of \fty? Some recent work of A.\ Gerges,
K. Orr and the first author suggests that this may be a good strategy
because $H_1$-bordism is determined by the most understood
\m-manifold invariants, namely the cohomology ring and the torsion
linking form.

\begin{theorem} {\rm (Amir Gerges \cite{Ge}; see
\cite{CGO} for d)}. Suppose $M_0$ and $M_1$ are closed, connected oriented
\m-manifolds. The following are equivalent.
\begin{enumerate}
\item[\rm{a)}] $M_0$ is $H_1$-bordant to $M_1$.
\item[\rm{b)}] $M_1$ is obtained from $M_0$ by surgery on an
admissible framed link $L$ in $M_0$. $($In fact $L$ may be chosen
to be a boundary link {\rm \cite[\S3.17]{CGO})}.
\item[\rm{c)}] There exist \m-manifolds $M_0=X_1$,
$X_2,\dots,X_n=M_1$ such that $X_{i+1}$ is obtained by $\pm1$
surgery on a null-homologous knot in $X_i$.
\item[\rm{d)}] There is an isomorphism $\phi:H_1(M_1)\to H_1(M_0)$
which induces isomorphisms between the $\bq/\bz$ linking forms and
between triple cup product forms $\bigotimes^3 H^1(M_i;\bz_n) \to
H^3(M_i;\bz_n)$ for $n=0$ and each
$n=p^r$\penalty-1000$(p$ prime$)$ where $p^r$ is the exponent of the
$p$-torsion subgroup of $H_1(M_i)$.
\item[\rm{e)}] There are isomorphisms $\phi_i:H_1(M_i)\to G$ $($a
fixed abelian group$)$ such that $(\phi_0)_*([M_0])=(\phi_1)_*([M_1])$
in $H_3(G)$.
\end{enumerate}
\end{theorem}

For example, note that 1.6e shows that for \m-manifolds with $H_1$
isomorphic to $0$, $\bz$ or $\bz^2$, there is only one $H_1$-bordism
class. For $H_1\cong\bz^3$ the non-negative integer
$|H^3(M_0)/(H^1(M_0)\cup H^1(M_0)\cup H^1(M_0))|$ is
a complete invariant. For $H_1\cong\bz_p$ ($p$ prime) there are two
equivalence classes, represented by $L(p,1)$ and $L(p,q)$ for any
mod $p$ quadratic non-residue $q$. For details and more
examples see \cite{CGO}.

Recall that the linking form can be computed directly from the
linking matrix associated to a surgery description of $M$ and that
such linking forms have been completely classified \cite{KK}. The
triple cup product forms can be calculated from the triple Milnor
invariants $\ovmu(123)$ of \m-component sublinks of a surgery
presentation of $M$ (\cite{T1}; Lemma 4.2). Hence, since
$H_1$-bordism is related to classical computable invariants, it
makes sense to separate the classification problem along these
lines. Although one {\sl need\/} not speak about invariants of
\fty\ for specific $H_1$-bordism classes, Proposition 1.5d makes it
clear that it would be more honest to do so.

One now sees that the degree zero \fts\ are precisely those which
are invariants of the isomorphism class of the triple ($H_1$,
linking form, triple cup product forms).

Our first major result, proved in section 2, is the finite generation of
the summands in the graded group $\calg(M)$ for any $M$; the analogous
theorem for spin manifolds is proved in \S6. In case $M$ is a homology
sphere this was proved by Ohtsuki \cite{O}. Henceforth, $\calm$ will
denote the (usual) theory of compact oriented $3$-manifolds (possibly
with boundary), while other theories will carry an adornment (such as
$\calm^{\spin}$ for spin manifolds).

\medskip
\noindent{\bf Theorem 2.1.} {\rm (finiteness theorem) {\it For any
compact oriented \m-manifold $M$ and any non-negative integer $\ell$, the
group $\calg_\ell(M)=(\calm_\ell/\calm_{\ell+1})(M)$ is finitely
generated. Therefore $\calo_\ell^A(M)$ is a finitely generated
$A$-module.}

\medskip

These finiteness results are directly related to the complexity of
calculation of invariants of \fty. Given any degree $n$, there is a
finite set $\{x_1,\dots,x_k\}\subset\calm(M)$, consisting of the union of
generating sets for $\calg_\ell$ for $0\le\ell\le n$, such that any
$\phi\in\calo_n(M)$ is completely determined by its values on
$\{x_i\}$, since any $\alpha\in (\calm/\calm_{n+1})(M)$ is a linear
combination of $\{x_i\}$. The techniques of section~2 suggest a
reasonable ``algorithm'' to calculate the coefficients.

In section 3 we show that the coefficients of the ``Conway
Polynomial'' of a \m-manifold $M$ with $\betti(M)=1$ are
non-trivial invariants of \fty, implying that $\calg_{2\ell}(M)$ has rank
at least~$1$. We also show that these invariants generate a polynomial
subalgebra of $\calo(M)$.

In section 4 we demonstrate that our theory is highly non-trivial, even
for manifolds with large first betti number, by exploiting the
$\bz_{p^k}$-valued invariants $\tau^d_p$ recently introduced by the
authors \cite{CM}.  These invariants were extracted from the quantum
$\so(3)$-invariants $\tau_p$ (for odd primes $p$).  Here it is shown that
they are of finite type and that they determine the quantum
$\so(3)$-invariants. This result appears to be new, even for homology
spheres.  In fact we show the stronger fact that $\tau_p$ is {\it
analytic\/}, which, loosely speaking, means that it is equal to the
``Taylor series'' constructed from its approximating ``polynomials''
$\tau^d_p$. In this regard $\tau_p$ is similar to the Jones and Conway
polynomials for knots.

By considering sequences of these invariants we establish {\sl rational\/}
non-triviality of the filtration on $\calm(M)$ for ``most" $3$-manifolds
$M$. We also provide strong evidence that Ohtsuki's theory for homology
spheres actually embeds in in the theory for manifolds $H_1$-bordant to
$M$.

The strongest results are for $H_1$-bordism classes containing a {\it
robust\/} manifold (see \ref{def:robust}). The list of robust manifolds
includes all rational homology spheres and the \m-torus $T=S^1\x S^1\x
S^1$, and is closed under connected sum. Therefore for any abelian group
$A$ whose rank is a multiple of \m\ there exists a robust \m-manifold $M$
with $H_1(M)\cong A$.

\medskip
\noindent{\bf Corollary \ref{cor:nontriviality}.} (part c) {\it If $M$
is robust, then each $\calg_{3k}(M)$ has positive rank, and so
$\calg(M)$ and $\calo^A(M)$ $($with $A=\bz$ or $\bq)$ are of infinite
rank.}

\medskip
The reader should note that $\calm/\calm_{\ell+1} \otimes \bq \cong
\bigoplus^\ell_{i=0} (\calg_i\otimes\bq)$ and so the non-triviality of
$\calg_i$ for $i\le\ell$ is directly related to the existence of
invariants of degree~$\ell$ (since $\calo_\ell$ with $\bq$ coefficients is
$\ho(\calm/\calm_{n+1},\bq)$).  For example, this result is used to prove
the existence of a finite type lift of the Casson invariant to arbitrary
$3$-manifolds that can detect homology sphere summands in $3$-manifolds
(Theorem \ref{thm:lift}).

For $H_1$-bordism classes $\cals(M)$ which are not robust we can still
show that the filtration $\calm_\ell(M)$ strictly descends as long as
some $\tau_p$ does not vanish identically on $\cals(M)$. If one
assumes that $M$ is {\it normal\/}, defined by the condition that
$\tau_p(M)\ne0$ for infinitely many $p$, then stronger results can be
obtained.  There exist normal manifolds with any prescribed homology;
in fact it is conceivable that all manifolds satisfy this condition.

\medskip\penalty-1000
\noindent{\bf Corollary \ref{cor:nontriviality}.} (parts a,b) {\it
If $\tau_p(M)\neq0$ for some prime $p>3$, then:
\begin{enumerate}
\item[\rm a)] For every positive integer $n$, there exists $m<\infty$
such that each $(\calm_\ell / \calm_{\ell+m})(M)$ has an element of
order at least $n$.
\item[\rm b)] Each $(\calm_\ell/\calm_\infty)(M)$ is of rank at
least $p-1$, and thus of infinite rank if $M$ is normal.  
\end{enumerate}}

\medskip
Finally we state the result which explains in what sense the complexity
of Ohtsuki's theory for homology spheres embeds in the general theory
for manifolds of high betti number.  In particular we paraphrase the part
of this result which relates to Ohtsuki's rational valued finite type
invariants of homology spheres.

\medskip
\noindent{\bf Corollary \ref{cor:istar}.} (parts b,c) 
{\it \begin{enumerate} 
\item[\rm b)] If $\tau_p(M)\neq0$ for some prime $p$, then the mod $p$
reduction of any of Ohtsuki's invariants is a linear combination of
invariants of the form $i^*(\phi)$ for $\phi\in\calo(M)$,  where by
definition $i^*(\phi)(x) = \phi(M\# x)$ \ $($and $M$ is assumed to be of
``minimal $p$-order" in its $H_1$-bordism class$)$.
\item[\rm c)] If $M$ is normal and $\Sigma_1$ and $\Sigma_2$
are homology spheres that can be distinguished by Ohtsuki's
invariants, then $M\#\Sigma_1$ and $M\#\Sigma_2$
can be distinguished by the finite type invariants $\tau^d_p$.
\end{enumerate}}

\medskip
In section 5 we describe an epimorphism from a finitely generated
group of ``Feynman diagrams'' to the graded group $\calg_\ell(M)$.
This is used to evaluate a few examples for small values of $\ell$.
The ``standard'' IHX and AS relations lie in the kernel but we show
that for some $M$ the kernel of this epimorphism is not completely
captured by these relations as is the case for homology spheres
\cite{GO3} \cite{Le}.

In section 6 we show that our theory for spin manifolds $\calo^\spin$
contains all of $\calo$ as well as the Rochlin invariant, which is
shown to be a degree three $\bz_{16}$-valued \fti. 

In section~7 we
briefly discuss several theories for \m-manifolds with non-empty
boundary.

In section 8 we investigate the category of oriented
\m-manifolds with marked $H_1$. We show that the coefficients of the
``Conway polynomial'' of the manifold are of \fty. We claim, but postpone
to a future paper, that Reidemeister torsion for
\m-manifolds with $H_1\cong\bz_{p^k}$ is analytic, in particular
determined by \fts.

In section 9 we sketch generalizations of our theory, in particular,
to a family of theories related to the lower-central-series.

In section 10 we note connections to the theories of \cite{GO1} for
rational homology spheres. We show that the invariant of Lescop
(including that of Casson-Walker) is of \fty\ (see also \S8). We also
indicate a relationship between our approach and a possible approach to
a theory of
\fts\ based on Heegard splittings and the mapping class group, whose
analogue for homology spheres was introduced and investigated in
\cite{GL3}.


\section{Finiteness}



\newsavebox{\band}
\savebox{\band}{\begin{picture}(0,0)
\setlength{\unitlength}{.48pt}
\multiput(0,0)(10,0){2}{\line(0,1){80}}
\end{picture}}

\newsavebox{\one}
\savebox{\one}{\begin{picture}(0,0)
\setlength{\unitlength}{.48pt}
\put(10,40){\oval(20,20)[l]}
\put(10,29.5){\line(1,0){20}}
\put(30,40){\oval(20,20)[r]}
\put(15,0)
  {\begin{picture}(0,0)
   \multiput(0,0)(10,0){2}{\line(0,1){25}}
   \end{picture}}
\put(15,35)
  {\begin{picture}(0,0)
   \multiput(0,0)(10,0){2}{\line(0,1){45}}
   \end{picture}}
\end{picture}}

\newsavebox{\two}
\savebox{\two}{\begin{picture}(0,0)
\setlength{\unitlength}{.48pt}
\put(10,40){\oval(20,20)[l]}
\put(10,29.5){\line(1,0){50}}
\put(60,40){\oval(20,20)[r]}
\put(30,50){\line(1,0){10}}
\multiput(15,0)(30,0){2}
  {\begin{picture}(0,0)
   \multiput(0,0)(10,0){2}{\line(0,1){25}}
   \end{picture}}
\multiput(15,35)(30,0){2}
  {\begin{picture}(0,0)
   \multiput(0,0)(10,0){2}{\line(0,1){45}}
   \end{picture}}
\end{picture}}

\def\longtwo{\begin{picture}(0,0)
\setlength{\unitlength}{.48pt}
\put(10,40){\oval(20,20)[l]}
\put(10,29.5){\line(1,0){80}}
\put(90,40){\oval(20,20)[r]}
\put(30,50){\line(1,0){40}}
\multiput(15,0)(30,0){3}
  {\begin{picture}(0,0)
   \multiput(0,0)(10,0){2}{\line(0,1){25}}
   \end{picture}}
\multiput(15,35)(60,0){2}
  {\begin{picture}(0,0)
   \multiput(0,0)(10,0){2}{\line(0,1){45}}
   \end{picture}}
\multiput(45,35)(10,0){2}{\line(0,1){10}}
\multiput(45,55)(10,0){2}{\line(0,1){25}}
\end{picture}}

\def\three{\begin{picture}(0,0)
\setlength{\unitlength}{.8pt}
\put(10,40){\oval(20,20)[l]}
\put(10,29.5){\line(1,0){80}}
\put(90,40){\oval(20,20)[r]}
\multiput(30,50)(30,0){2}{\line(1,0){10}}
\multiput(15,0)(30,0){3}
  {\begin{picture}(0,0)
   \multiput(0,0)(10,0){2}{\line(0,1){25}}
   \end{picture}}
\multiput(15,35)(30,0){3}
  {\begin{picture}(0,0)
   \multiput(0,0)(10,0){2}{\line(0,1){45}}
   \end{picture}}
\end{picture}}

\def\twist{\begin{picture}(0,0)
\put(10,42.2){\oval(20,22)[l]}
\put(9.8,31){\line(1,0){17}}
\put(44.5,31){\line(1,0){15}}
\put(60,42.2){\oval(20,22)[r]}
\put(27,31){\line(3,4){17}}
\put(27,54){\line(3,-4){7.5}}
\put(37,41){\line(3,-4){7.5}}
\multiput(15,0)(31,0){2}
  {\begin{picture}(0,0)
   \multiput(0,12)(10,0){2}{\line(0,1){15}}
   \end{picture}}
\multiput(15,35)(31,0){2}
  {\begin{picture}(0,0)
   \multiput(0,0)(10,0){2}{\line(0,1){35}}
   \end{picture}}
\end{picture}}


In this section we prove the main finiteness result in
the oriented category. We also show that the group of \fts\ forms a
filtered commutative algebra.

\begin{theorem}{\rm (finiteness theorem)} For any compact
oriented \m-manifold $M$ and any integer $\ell$, the group 
$\calg_\ell(M) = (\calm_\ell/\calm_{\ell+1})(M)$ is finitely
generated. Therefore $\calo^A_\ell(M)$ is a finitely generated
$A$-module.
\end{theorem} 

The proof is very similar to that of the corresponding result of
Ohtsuki \cite{O}, except that one must deal with admissible links in $M$
rather than $S^3$. Philosophically, all of Ohtsuki's local lemmas work
except that the ones whose proofs involve ``blowing up or down'' can
only be applied to $\pm1$ framed circles. Hence the ``braiding
lemma'' and the ``framing lemma'' do not hold in full generality, and in
particular, most of the properties of \cite{GO1} do not hold.

\medskip
\noindent {\it Proof of 2.1}. Fix $M$ and a non-negative integer
$\ell$.  Following \cite{O} we write $\sim$ for the equivalence relation
on $\calm_\ell(M)$ induced by the projection to $\calg_\ell(M)$.  Our
basic tool is Ohtsuki's ``fundamental lemma" (\cite{O}, Lemma~2.2)
which generalizes to the present setting.

\begin{lemma}\label{lem:fundamental} {\rm (fundamental lemma)} If $L\cup
K$ is an admissible link in $M$ then $[M,L]\sim[M_K,L]$ where $M_K$ is
surgery on $K$ and the latter $L$ is the image of $L$ in $M_K$. $($Note
that $K$ may have more than one component\/$)$.
\end{lemma}

\begin{proof} Since $L$ has $\ell$ components,
$[M,L]\sim[M,L\cup\delta K]$, because each of the non-empty terms in
$\delta K=\sum_{S<K}(-1)^sS$ gives rise to an element of
$\calm_{\ell+1}$. But $[M,L\cup\delta K]=[M_K,L]$ by Lemma 1.4.
\end{proof}

Recall that by definition $\calm_\ell(M)$ is spanned by elements of the
form $[M',L']$, where $M'$ is $H_1$-bordant to $M$ and $L'$ is an
admissible $\ell$-component link in $M'$.  If we work modulo
$\calm_{\ell+1}(M)$, however, we need only consider the case $M'=M$.  In
other words $\calg_\ell(M)$ is generated by elements of the form $[M,L]$,
where $M$ is any chosen ``basepoint'' in the $H_1$-bordism class and $L$
has $\ell$~components (cf. \cite{O} Lemma~2.3).

\begin{lemma}\label{lem:basepoint} {\rm (basepoint lemma)} Suppose $M$
and $M'$ are $H_1$-bordant and $L'$ is an admissible link of $\ell$
components in $M'$. Then there exists an admissible link $L$ in $M$ with
$\ell$ components such that $[M',L']\sim[M,L]$.
\end{lemma}

\begin{proof} By Theorem 1.6b we may assume $M\cong M'_K$, where $K$ is
an admissible link in $M'$. $K$ may be varied by an isotopy in $M'$ until
$L'\cup K$ is admissible in $M'$. It then follows from the fundamental
lemma (\ref{lem:fundamental}) that $[M',L']\sim[M'_K,L'] = [M,L]$ where
$L$ is the image of $L'$ in $M$.
\end{proof}

The next result, generalizing Lemma 2.5 of \cite{O}, shows how to
arrange that all framings be $+1$.

\begin{lemma}\label{lem:framing} {\rm (framing lemma)} Suppose $L$ is an
$\ell$-component admissible link in $M$ with framing $-1$ on the component
$K$. Let $L'$ be the link $L$ with the framing on $K$ changed to $+1$.
Then $[M,L]\sim-[M,L']$.
\end{lemma}

\begin{proof} Let $K'$ be a $+1$-framed parallel of $K$ with $\ell
k(K,K')=0$.  Set $J = L-K$, so $L'=J\cup K'$.  Observe that the
pairs $(M,J)$ and $(M_{K\cup K'},J)$ are homeomorphic, since doing $+1$
and $-1$ surgery on parallels of the core of a solid torus $T$ yields a
manifold diffeomorphic to $T$ fixing $\p T$, and so $[M,J] = [M_{K\cup
K'},J]$.  Now by the fundamental lemma, $[M, L] \sim [M_{K'},L] =
[M_{K'},J] - [M_{K\cup K'},J] = [M_{K'},J] - [M,J] = -[M,L']$.
\end{proof}

The ``braiding lemma" of Ohtsuki also generalizes to the present context.
The key proviso is that the unknotted component $K$ (in the statement
below) is $\pm1$-framed. The analogous result of (\cite[Fig.1]{GO1})
without this proviso, is false. In the following, non-integral framings
are allowed on $J$. For convenience we now assume that $M$ is closed. The
modifications necessary in the case of non-empty boundary are discussed
in section~7.

\begin{lemma}\label{lem:braid} {\rm (braiding lemma)} Suppose $J\cup L$
is a framed link in $S^3$ such that $L$ $($with $\ell$ components\/$)$ is
admissible  in $M=S^3_J$, and such that each component of $J$ has zero
linking number with each component of $L$. In addition suppose that
$L$ has an unknotted component $K$, and that the components of
$J\cup L$ which pierce a disk $D$ spanned by $K$ have been divided
into $m$ groups of strands, represented by ``bands'' in
\textup{Figure~2.6a}, in such a way that each component passes
algebraically zero times through each band. Number the bands, and for
each increasing sequence $1\le i_1 < \cdots < i_k \le m$, let
$L_{i_1\cdots i_k}$ be the framed link obtained from $L$ by replacing $K$
with a curve $K_{i_1\cdots i_k}$ in $D$ $($with the
same framing as $K)$ which encircles the bands $i_1,\dots,i_k$ while
passing in front of the other bands.  Then
$$
[M,L]\sim\sum_{i,j=1}^m [M,L_{ij}] - (m-2)\sum_{i=1}^m [M,L_i].
$$
\end{lemma}

The case $m=3$ is illustrated in Figure 2.6.


\begin{center}
\begin{picture}(280,128)
\setlength{\unitlength}{.8pt}

\put(0,10)
{\begin{picture}(0,0) 
\put(0,30)
  {\begin{picture}(0,0)
  \put(-10,20){$K$}
  {\three}
  \end{picture}}
\put(30,-25){a) $J\cup L$}
\end{picture}}

\put(130,10)
{\begin{picture}(0,0)
\put(0,90)
  {\begin{picture}(0,0)
  \put(0,0){\usebox{\two}}\put(45,0){\usebox{\band}}
  \put(-15,10){$K_{12}$}
  \end{picture}}
\put(81,90)
  {\begin{picture}(0,0)
  \put(0,0){\longtwo}
  \put(-15,10){$K_{13}$}
  \end{picture}}
\put(163,90)
  {\begin{picture}(0,0)
  \put(0,0){\usebox{\band}}\put(10,0){\usebox{\two}}
  \put(52,10){$K_{23}$}
  \end{picture}}
\put(0,0)
  {\begin{picture}(0,0)
  \put(0,0){\usebox{\one}}\multiput(28,0)(18,0){2}{\usebox{\band}}
  \put(-15,30){$K_{1}$}
  \end{picture}}
\put(84,0)
  {\begin{picture}(0,0)
  \multiput(0,0)(47,0){2}{\usebox{\band}}\put(17,0){\usebox{\one}}
  \put(8,33){$K_{2}$}
  \end{picture}}
\put(163,0)
  {\begin{picture}(0,0)
  \multiput(0,0)(18,0){2}{\usebox{\band}}\put(30,0){\usebox{\one}}
  \put(53,30){$K_{3}$}
  \end{picture}}
\put(50,-25){b) $J\cup L_{i}$ and $J\cup L_{ij}$}
\end{picture}}

\end{picture}
\end{center}

\medskip
\centerline{Figure 2.6}
\setcounter{theorem}{6}
\medskip

\begin{proof}
Following \cite{GL} we give an ``algebraic" proof. Assume that the framing
on $K$ is $+1$; the other case then follows from the framing lemma
(\ref{lem:framing}).   Let $q = [M,L]$ and $x = [M,\hat L]$, where $\hat
L$ is obtained by ``blowing down" K, that is removing $K$ and putting a
full left twist in all the bands.  Note that $q \in \calm_\ell$ and $x \in
\calm_{\ell-1}$.  Furthermore, if we set $1 = [M,L-K]$ then $q = 1-x$ by
Lemma 1.4.  In a completely analogous way, we define $q_{i_1\cdots
i_k}$ and $x_{i_1\cdots i_k}$ with $q_{i_1\cdots i_k} = 1-x_{i_1\cdots
i_k}$ (note that $q = q_{1\cdots m}$ and $x = x_{1\cdots m}$), and with
this notation, the lemma states that $q \sim \sum q_{ij} - (m-2)\sum
q_i$.  

Now the key to the proof is the elementary observation that a full
left twist in a collection of bands is a product of left twist in pairs
of bands and in the individual bands.  Explicitly 
$$
x = \prod_{i,j=1}^m x_{ij} \prod_{i}^m x_i^{2-m}
$$
with lexicographic ordering in the first product.  Here the product (left
to right) corresponds to the stacking (bottom to top) of the associated
tangles, and $x_i^{-1} = 1 + q_i + q_i^2 + \cdots$ is a right handed
twist in the $i$th band.  Substituting the $q$'s for the $x$'s and
expanding the right hand side, we obtain $1-q = 1 - \sum q_{ij} +
(m-2)\sum q_i +$ quadratic terms (which vanish in $\calg_\ell$), and the
result follows.
\end{proof}

Another useful local result which generalizes to our setting is
Ohtsuki's ``half-twist lemma'' (stated incorrectly in Figure~4.3 of
\cite{O}, but later corrected in Figure~5 of \cite{GO3}).

\begin{lemma}\label{lem:twist} {\rm (half-twist lemma)} Assume the
hypotheses of the braiding lemma $(\ref{lem:braid})$ with $m=2$, and
suppose that $L'$ is obtained from $L$ by replacing $K$ by a
half-twisted unknot $K'$, as shown in {\rm Figure~2.8}.  Then
$$
[M,L']\sim-[M,L]+2[M,L_1]+2[M,L_2].
$$
\end{lemma} 

\noindent (Recall that $L_1$ and $L_2$ are obtained from $L$ by
replacing $K$ with unknots encircling the first and second bands,
respectively.)


\medskip
\begin{center}
\begin{picture}(100,70)
\put(15,0){\twist}
\put(0,40){$K'$}
\put(40,-5){$J\cup L'$}
\end{picture}
\end{center}

\centerline{\rm Figure 2.8}
\setcounter{theorem}{8}
\medskip

\begin{proof}
Adopting the notation of the preceding proof, and letting $q' = 1~-~x'
 = [M,L']$, we must show $q' \sim -q + 2q_1 + 2q_2$. By Lemma 1.4
we compute
$q' = 1-x^{-1}x_1^2x_2^2 = 1 - (1+q+q^2+\cdots) (1-q_1)^2(1-q_2)^2 \sim
-q + 2q_1 + 2q_2$. 
\end{proof}

Recall, following Levine, that the ordered oriented links $L$ and $L'$
in $S^3$ are said to be {\it surgery equivalent\/} if $L\cong L_0\sim
L_1\sim\dots\sim L_k\cong L'$ where $L_i\sim L_{i+1}$ means that
there is a $2$-disk $D_i$ in $S^3$ such that $\p D_i$ is disjoint from
and has zero linking number with each component of $L_i$ and such that
$\pm1$ surgery on $\p D_i$ transforms $L_i$ to $L_{i+1}$ \cite{L}.

\begin{lemma}\label{lem:surgery} {\rm (surgery lemma)} Assume the
hypotheses of the braiding lemma $(\ref{lem:braid})$.  If $J\cup L$ is
surgery equivalent to $J\cup L'$ then $[M,L]\sim[M,L']$, where $M =
S^3_J$ and  the framings on $L'$ are taken equal to the corresponding
framings on $L$.
\end{lemma}

\begin{proof} It suffices to assume the weaker condition that there is a
$\pm1$-framed knot $K$ in $S^3-(J\cup L)$ having zero linking number with
the components of $J\cup L$ such that the pair $(S^3_K,J\cup L)$ is
homeomorphic to $(S^3,J\cup L')$. Hence $(S^3_{J\cup K},L)=(M_K,L)$ is
homeomorphic to $(S^3_J,L')=(M,L')$, and so by the fundamental lemma
$[M,L] \sim [M_K,L] = [M,L']$. 
\end{proof}

We now continue with the proof of Theorem 2.1, using Levine's surgery
equivalence classification for {\sl arbitrary\/} links in $S^3$
\cite{L}.  Consider, as above, $M=S^3_J$.  (What follows is all fairly
easy if $J$ has zero linking numbers --- and in this case was done by
Ohtsuki without Levine's theorem --- but this is not always possible to
assume.\foot{although it is, for example, if
$H_1(M)$ has no
$2$-torsion})

Fix an orientation and an ordering for the components of $J$, and choose
a family of {\it base paths}, i.e.\ disjoint paths from a chosen
basepoint in $S^3-J$ to each of the components of $J$.  (In general we
shall refer to any oriented, ordered, based link simply as a {\it based
link}.)  

Consider the family of based links $J\cup L$, where $L$ has
$\ell$ components.  For later notational convenience, assume that the
ordering index for $J\cup L$ runs from $1$ to $\ell+m$ (so $m$
is the number of components in $J$) with $L$ corresponding to
$1,\dots,\ell$.  Of particular interest is the case when $L = T$, where
$T$ is a {\sl trivial\/} link lying in in a ball disjoint from $J$ (and
its base paths).  We shall define a ``special" class of based links
related to
$J\cup T$.

\begin{definition} A based link $J\cup L$ in $S^3$ is {\it special\/} if
it is obtained \linebreak from $J\cup T$ by replacing some number of
disjoint \m-string trivial tangles $(B^3,\g_i\cup\g_j\cup\g_k)$, by (one
of $2$ possible) ``Borromean tangle(s)" $(B^3,\g'_i\cup\g'_j\cup\g'_k)$
subject to the condition that $\{\g_i\,\g_j,\g_k\}$ are arcs of $3$ {\sl
distinct\/} components of $J\cup T$ with at least one being a component
of $T$. Such a replacement is called a {\it Borromean replacement} of
type $(i,j,k)$. The geometric number of such is denoted $n_{ijk}$.  
\end{definition}

Let $[M,L]$ be an arbitrary generator of $\calg_\ell(M)$. By the
framing lemma (\ref{lem:framing}) we may assume that all components of
$L$ have framing $+1$. Isotope $L$ in $M$ so that $L\subset S^3_J$ is
disjoint from the surgery tori and each component of $L$ has zero
linking with each component of $J$.  

Now consider the link $J\cup L$ in $S^3$.  Order and orient the
components the components of $L$ arbitrarily, and choose base paths
which extend the basing of $J$.  Thus $J\cup L$ becomes a {\sl based\/}
link in the sense defined above.  By \cite[p.51]{L} there is a set
$\{\mu_{ij},a_{ijk}\}=\mu(J\cup L)$ of integers associated to this based
link. The $\mu_{ij}$ are the linking numbers and the $a_{ijk}$ are
``lifts'' of Milnor's triple $\ovmu$-invariants. Compare these to
$\mu(J\cup T)$. Clearly the linking numbers agree. Moreover $a_{ijk}$
depends only on the \m-component based sublinks \cite[p.54, paragraph
3]{L}. A \m-component sublink $\{J_i,J_j,J_k\}$ is independent of $L$ and
hence the corresponding
$a_{ijk}$ for $J\cup L$ and $J\cup T$ agree. Thus, in the following
discussion we restrict to those $(i,j,k)$ corresponding to a
\m-component sublink containing at least one component of $L$ or $T$
(so $i\le\ell$ by our ordering conventions). These may be altered by
Borromean replacements. By the proof of Theorem~C of
\cite{L}, there exists a {\sl special\/} link $J\cup L_s$ such that
$\mu(J\cup L_s)=\mu(J\cup L)$ where each Borromean replacement involves at
least one component from $T$. By Theorem~D of that paper, $J\cup L_s$ is
surgery equivalent to $J\cup L$. By the surgery lemma
(\ref{lem:surgery}) $[M,L]\sim[M,L_s]$. Therefore we have shown that
$\calg_\ell(M)$ is spanned by elements of the form
$[S^3_J,L]$ where $J\cup L$ is special and all framings are $+1$.

By the proof of Theorem~C of \cite{L} the invariants $a_{ijk}$ of a
special link differ from those of $J\cup T$ by precisely the algebraic
number of Borromean replacements of type $(i,j,k)$. Therefore two special links
are surgery equivalent if and only if the {\sl algebraic\/} number of
tangle replacements of type $(i,j,k)$ is the same for each triple
$i\!<\!j\!<\!k$.  Consequently we need only consider {\sl one\/} special
link for each possible value of the collections $\{a_{ijk} \, | \,
i\!<\!j\!<\!k \}$ (with all indices between $1$ and $\ell+m$, and
$i\le\ell$ as usual).  The corresponding set of $[S^3_J,L]$ (using
$+1$ framings) forms a spanning set for $\calg_\ell(M)$, which is still
{\sl infinite\/} since the $a_{ijk}$ can be arbitrary. 

Choose such a set for which the {\sl actual\/} number $n_{ijk}$ of
replacements of type $(i,j,k)$ is equal to $|a_{ijk}|$, for each $i,j,k$.
Now apply the braiding lemma (\ref{lem:braid}), noting that the links on
the right hand side are all special if the one on the left is special,
to show that one need only consider special links for which there are at
most {\sl two\/} replacements involving each component of $L$.  This then
yields a {\sl finite\/} spanning set for $\calg_\ell(M)$, corresponding
to collections $\{a_{ijk} \, | \, i\!<\!j\!<\!k \}$ for which each of
the indices $1,\dots,\ell$ appears in at most two non-zero $a_{ijk}$'s. 
This completes the proof of Theorem~2.1.
\qed 


\begin{remark}
With a little more work it can be seen that only links with each non-zero
$a_{ijk}$ equal to $+1$ are needed in the generating set:
Consider a special link representing one of the generators.  Fix $i<j<k$
and consider the number of replacements $n_{ijk}$ of type $(i,j,k)$.
This number is either $0$,
$1$ or $2$ (according to the construction above) and we are only
interested in the latter two cases. 

If $n_{ijk}=1$ then $a_{ijk}=\pm1$. In case $a_{ijk}=-1$ and $L_k$ is not
involved in any other replacements then simply change the orientation of
$L_k$ to get $a_{ijk}=+1$. In case $L_k$ is involved in one other
replacement, apply the half-twist lemma (\ref{lem:twist}) to reduce to
situations in which it is involved in only one replacement or the
$a_{ijk}$ is changed to $+1$. 

If $n_{ijk}=2$ then $a_{ijk} = \pm2$, and changing the orientation on
$L_k$ if necessary gives $a_{ijk} = 2$.  Now apply \ref{lem:twist}
again to reduce to cases in which $a_{ijk}=0$ (for which we can
substitute a simpler special link) or $n_{ijk} = 1$.  Thus we obtain a
spanning set with each $a_{ijk}$ equal to $0$ or $1$ and
$n_{ijk}=a_{ijk}$. 

In summary, if we think of $L = \{L_{1}, \dots, L_{\ell}\}$ and
$J=\{J_1,\dots,J_m\}$, then we have found a spanning set in one-to-one
correspondence with the {\sl subsets} of the index
set $U = \{(i,j,k)\,|\,1\!\le\!i\!<\!j\!<\!k\!\le\!\ell+m, \
i\!\le\!\ell \}$ in which each of the indices $1,\dots,\ell$ appears
at most twice.
\end{remark}

\medskip
We now prove that $\calo$, the group of all \fts, and $\calo(M)$,
the group of all \fts\ for manifolds in the $H_1$-bordism class of
$M$, have the structure of algebras. As usual, one must be careful to
define $\la\la'$ as the {\sl linear extension\/} to $\calm$ of the usual
product of functions on $\cals$. So for example if $M$ and $N$ are
manifolds, $\la\la'(M+N)=\la(M)\la'(M)+\la(N)\la'(N)$.

\begin{proposition}\label{prop:algebra} If $\la\in\calo_p$,
$\la'\in\calo_q$ then $\la\la'\in\calo_{p+q}$.
\end{proposition}

\begin{proof} We shall show that
$$
\la\la'([M,L]) = \sum_{S<L}\la([M,S])\la'([M_S,L-S])
$$
which will complete the proof since if $\ell>p+q$ then either $s>p$ or
$\ell-s>q$. Rewrite $\la'([M_S,L-S])$ as $\sum_{T>S}(-1)^{t-s}\la'(M_T)$.
Then the right hand side above can be expressed as
$$
\sum_{S<L}\bigl[\sum_{R<S}(-1)^r\la(M_R) \sum_{T>S}(-1)^{t-s}
\la'(M_T)\bigr].
$$
Rearranging the order of summation gives
$$
\sum_{R<T<L} \bigl[(-1)^{r+t}\la(M_R)\la'(M_T) \sum_{R<S<T} (-1)^s\bigr]
$$
The inner sum vanishes unless $R=T$, since it is an alternating sum of
binomial coefficients.  For $R=T$ we get $(-1)^t \la(M_T) \la'(M_T)$,
and summing over $T<L$ gives $\la\la'([M,L])$ as desired.
\end{proof}

Thus if $A$ is a commutative ring then $\calo$ is a filtered
commutative ring in which $A$ occurs naturally as the subring of
constant functions. The multiplication then makes $\calo$ a filtered
commutative $A$-algebra and $\calo(M)$, for any $M$, a subalgebra.


\section{The Conway polynomial}



\def\hcirc{\begin{picture}(0,0) 
\put(10,38){\oval(16,16)[l]}
\put(40,38){\oval(16,16)[r]}
\put(10,29.5){\line(1,0){30}}
\put(20,46){\line(1,0){10}}
\end{picture}}

\def\fcirc{\begin{picture}(0,0) 
\put(25,38){\oval(46,16)}
\end{picture}}

\def\vcirc{\begin{picture}(0,0) 
\put(80,63){\oval(20,20)[t]}
\put(80,17){\oval(20,20)[b]}
\put(69.6,16){\line(0,1){47}}
\put(89.6,25){\line(0,1){30}}
\end{picture}}

\def\vcirco{\begin{picture}(0,0) 
\put(80,63){\oval(20,20)[t]}
\put(80,17){\oval(20,20)[b]}
\put(69.6,16){\line(0,1){47}}
\put(89.6,25){\line(0,1){15}}
\put(89.6,45){\line(0,1){9}}
\end{picture}}

\def\vcirct{\begin{picture}(0,0) 
\put(80,63){\oval(20,20)[t]}
\put(80,17){\oval(20,20)[b]}
\put(69.6,16){\line(0,1){12}}
\put(89.6,25){\line(0,1){3}}
\put(69.6,36){\line(0,1){28}}
\put(89.6,36){\line(0,1){18}}
\end{picture}}

\def\lcirc{\begin{picture}(0,0) 
\put(14.8,15){\line(0,1){10}}
\put(35,15){\line(0,1){10}}
\put(14.8,35){\line(0,1){30}}
\put(35,35){\line(0,1){30}}
\put(54.1,65){\oval(77.8,55)[t]}
\put(54.1,15){\oval(77.8,55)[b]}
\put(83,65){\oval(20,20)[b]}
\put(83,15){\oval(20,20)[t]}
\put(54,74.5){\oval(37,1)[tr]}
\put(54,65){\oval(37,20)[tl]}
\put(52,7.5){\oval(37,5)[br]}
\put(54,15){\oval(37,20)[bl]}
\end{picture}}

\def\lcircs{\begin{picture}(0,0) 
\put(14.8,15){\line(0,1){9}}
\put(35,15){\line(0,1){9}}
\put(14.8,32){\line(0,1){33}}
\put(35,32){\line(0,1){33}}
\put(54.1,65){\oval(77.8,55)[t]}
\put(54.1,15){\oval(77.8,55)[b]}
\put(83,65){\oval(20,20)[b]}
\put(83,15){\oval(20,20)[t]}
\put(54,74.5){\oval(37,1)[tr]}
\put(54,65){\oval(37,20)[tl]}
\put(54,4.5){\line(1,0){12}}
\put(68,8.5){\oval(4,4)[br]}
\put(54,15){\oval(37,20)[bl]}
\end{picture}}

\def\lcircso{\begin{picture}(0,0) 
\put(14.8,15){\line(0,1){9}}
\put(35,15){\line(0,1){9}}
\put(14.8,32){\line(0,1){8}}
\put(14.8,45){\line(0,1){20}}
\put(35,32){\line(0,1){33}}
\put(54.1,65){\oval(77.8,55)[t]}
\put(54.1,15){\oval(77.8,55)[b]}
\put(83,65){\oval(20,20)[b]}
\put(83,15){\oval(20,20)[t]}
\put(54,74.5){\oval(37,1)[tr]}
\put(54,65){\oval(37,20)[tl]}
\put(54,4.5){\line(1,0){12}}
\put(68,8.5){\oval(4,4)[br]}
\put(54,15){\oval(37,20)[bl]}
\end{picture}}

\def\Bor{\begin{picture}(0,0)
\put(0,0){\thicklines\hcirc}
\put(0,0){\lcirc}
\put(0,0){\vcirc}
\end{picture}}

\newsavebox{\bor}
\savebox{\bor}{\begin{picture}(0,0)
\setlength{\unitlength}{.6pt}
\put(0,0){\lcircs}
\put(0,0){\vcirc}
\end{picture}}

\newsavebox{\boro}
\savebox{\boro}{\begin{picture}(0,0)
\setlength{\unitlength}{.6pt}
\put(0,0){\lcircso}
\put(0,0){\vcirco}
\end{picture}}


In this section we will show that $\calg_{2n}=\calm_{2n}/\calm_{2n+1}$ is
infinite for each $n\ge0$ by exhibiting specific \fts\ $C_{2n}$ of
degree $2n$. The invariant $C_{2n}(M)$ will be defined to be the
coefficient of $z^{2n}$ in the ``Conway polynomial'' of $M$ if
$\betti(M) = 1$, and zero otherwise. Since C.~Lescop's invariant
\cite{Les} is $C_2(M)-\frac1{12}|\tor H_1(M)|$ for manifolds with
$\betti = 1$, this shows that her invariant is \fty\ of degree~$2$ on
this $H_1$-bordism class. Moreover we show that the set
$\{C_2,C_4,\dots\}$ is a basis of a polynomial subalgebra of
$\calo$.  (Note that $C_0$ is excluded since it is identically equal to
$1$ on manifolds of first betti number one, whence $C_0^2 = C_0$ is a
polynomial relation in $\calo$.)

A closed oriented \m-manifold $M$ with $\betti(M)=1$ has a unique
Conway polynomial $\nabla_M(z)=1+a_2z^2+a_4z^4+\dots$ defined
as follows. Let $\wt M$ denote the infinite cyclic cover of $M$.
Evidently $H_1(\wt M)$ has two $\bz[t,t^{-1}]$ module structures,
differing by $t\mapsto t^{-1}$. The {\it Alexander polynomial} of $M$
is defined to be the order of (either of) these torsion modules
divided by $|\textup{Tor}(H_1(M))|$. It can also be identified with the
Alexander polynomial of a suitable knot. Indeed $M$ can be constructed by
$0$-framed surgery $\Sigma_K$ on a null-homologous knot $K$ in a rational
homology sphere $\Sigma$ (\cite[\S5.1.1]{Les}), and it is an easy exercise
to see that the Alexander module $H_1(\widetilde{\Sigma-K})$ of $K$ is
isomorphic to $H_1(\wt M)$ (where the module structure is determined by a
choice of orientation on $K$). Now recall that the Alexander polynomial of
$K$ in $\Sigma$ is defined to be the order of this torsion module divided
by $|H_1(\Sigma)|$, and may be computed as $\det(tV-V^T)$ where $V$
is any (rational) Seifert matrix for $K$ in $\Sigma$
(\cite[\S2.3.12--13]{Les}). Since $|H_1(\Sigma)|=|\textup{Tor}(H_1(M))|$,
this coincides with the Alexander polynomial of $M$. Of course this
polynomial is only defined up to a unit $\pm t^n$ in $\bq[t,t^{-1}]$, but
it can be normalized by setting $\Delta_M(t) = \Delta_{K,\Sigma}(t) =
\det(t^{1/2}V-t^{-1/2}V^T)$ so that $\Delta_M(t^{-1}) = \Delta_M(t)$ and
$\Delta_M(1)=1$. This yields a uniquely defined Alexander polynomial, a
Laurent polynomial in $t^{1/2}$ with rational coefficients, which can be
shown to be an honest polynomial in $(t^{1/2}-t^{-1/2})^2$
(\cite[\S2.3.14--15]{Les}). Substituting $z$ for $t^{1/2}-t^{-1/2}$ then
yields the {\it Conway polynomial} $\nabla_M(z)$ of $M$, or equivalently
$\nabla_{K,\Sigma}(z)$ of $K$ in $\Sigma$, an element of
$\bq[z^2]$.\foot{$\nabla_{K,\Sigma}(s^{-1}\!-\!s)$ coincides with the
polynomial defined by Boyer and Lines \cite{BoL}.}  Extending linearly by
setting $\nabla_M=0$ if $\betti(M)\neq1$ yields a polynomial valued
invariant $\nabla:\calm\to\bq[z^2]$.

We shall also need the fact that the Conway polynomial can be defined
for {\sl links} in rational homology spheres (see e.g.\ \cite{BoL}). In
particular if $K$ is a $k$-component null-homologous {\sl
oriented\/} link in a rational homology sphere $\Sigma$, then 
$\nabla_{K,\Sigma}(z)$ is of the form $z^{k-1}(a_0+a_1z^2+\dots)$.
The crucial fact needed here, due to Boyer and Lines, is
that $\nabla_K=\nabla_{K,\Sigma}$ satisfies the familiar recursion formula
$\nabla_{K^+}-\nabla_{K^-}=-z\nabla_{K^0}$ (see
\cite[\S2.3.16]{Les}).

The main result of this section is the following.

\begin{theorem} Let $n$ be a nonnegative integer and $M$ be a closed,
oriented \m-manifold. Consider the \m-manifold invariant
$C_{2n}:\calm \to \bq$ which assigns to $M$ the coefficient of $z^{2n}$ in
the Conway polynomial $\nabla_M$ if $\betti(M)=1$, and zero otherwise.
Then $C_{2n}$ is \fty\ of degree~$2n$.
\end{theorem}

\medskip
\noindent{\bf Remark.} If the domain of $C_{2n}$ is restricted to
integral homology $S^1\x S^2$'s then $C_{2n}$ is an integral
invariant.

\bigskip
The theorem will follow easily from Theorem \ref{thm:divisibility} below
concerning the divisibility of the alternating sum of Conway polynomials
of links in a rational homology sphere. A realization result, Proposition
\ref{prop:realization}, is then also needed to show that $C_{2n}$ has
degree precisely $2n$.

Suppose $K$ is a null-homologous oriented link in a rational homology
sphere $\Sigma$, and $L=\{L_1,\dots,L_\ell\}$ is an admissible framed link
in $\Sigma$ (see 1.2). We say that $L$ is {\it admissible in
$(\Sigma, K)$} if $K$ bounds a Seifert surface in $\Sigma - L$, or
equivalently $L$ is disjoint from $K$ and $\ell k(K,L_i)=0$ for all $i$.
If $S$ is a sublink of such an $L$ then $\Sigma_S$ is again a rational
homology sphere in which the image of $K$ remains a link. For brevity we
continue to denote this image by $K$ whenever possible. We shall also use
the abbreviation $\nabla_K(S)$ for the Conway polynomial of $K$ in
$\Sigma_S$ for any sublink $S$ of $L$,
$$
\nabla_K(S) = \nabla_{K,\Sigma_S}\ ,
$$
and $\nabla_K(\delta L)$ for $\sum_{S<L}(-1)^s\nabla_K(S)$.

\begin{theorem}\label{thm:divisibility} If $K$ is a null-homologous
oriented link in a rational homology sphere $\Sigma$ and $L$ is an
admissible link of $\ell$~components in $(\Sigma, K)$ then $z^\ell$
divides $\nabla_K(\delta L)$.
\end{theorem}

The proof will be given later in this section.

\begin{example}\label{ex:conway} Suppose $K$ is the trivial knot in
$\Sigma=S^3$ (with either orientation) and $L = K_1\cup K_2$ is the
$+1$-framed $2$-component link shown in Figure 3.4. Then
$(\Sigma_{K_1},K) \cong (\Sigma_{K_2},K) \cong
(\Sigma,K)\cong(\Sigma,\text{unknot})$, whereas $(\Sigma_L,K)$ is the
right-handed trefoil knot (most easily seen by ``blowing-down'' $L$
\cite{K}). Thus $\nabla_K(\delta L) = 1-1-1+(1+z^2)=z^2$, which is
divisible by $z^2$ as predicted by Theorem
\ref{thm:divisibility}.


\begin{center}
\begin{picture}(90,55)
\put(12,0){\setlength{\unitlength}{.7pt}\Bor}
\put(0,30){$K$}
\put(80,0){$K_1$}
\put(80,24){$K_2$}
\end{picture}
\end{center}

\medskip
\centerline{\rm Figure 3.4: $L = K_1\cup K_2$}
\bigskip

This example can be generalized by taking ``parallel" copies to obtain the
$+1$-framed $2n$-component link $L_{2n}$ shown in Figure 3.5.


\begin{center}
\begin{picture}(240,80)
\multiput(10,0)(60,0){2}{\usebox{\bor}}
\put(175,0){\usebox{\bor}}
\thicklines
\put(10,43.5){\oval(40,52)[l]}
\put(235,43.5){\oval(40,52)[r]}
\put(10,69.5){\line(1,0){130}}
\put(175,69.5){\line(1,0){60}}
\multiput(10,17)(60,0){2}{\line(1,0){35}}
\put(130,17){\line(1,0){10}}
\multiput(145,17)(5,0){6}{\line(1,0){2}}
\multiput(145,69.5)(5,0){6}{\line(1,0){2}}
\put(175,17){\line(1,0){35}}
\put(-20,65){$K$}
\put(40,-20){$L^1$}
\put(100,-20){$L^2$}
\put(205,-20){$L^n$}
\end{picture}
\end{center}

\bigskip
\centerline{\rm Figure 3.5: $L_{2n} = L^1 \cup \cdots \cup L^n$}
\setcounter{theorem}{5}
\end{example}

\begin{proposition}\label{prop:realization} Let $K$ be an
unknot in $\Sigma = S^3$ $($with either orientation\/$)$ and
$L_{2n}$ be the $+1$-framed $2n$-component link shown in Figure $3.5$,
where each $L^i$ is a copy of the $2$-component link $L$ in Figure $3.4$. 
Set $\la_{2n} = [\Sigma_K,L_{2n}]$, where $K$ is given the zero framing.
$($Note that $\Sigma_K = S^1\x S^2$ since $K$ is unknotted.$)$  Then
\begin{enumerate}
\item[\rm a)] $\nabla_K(\delta L_{2n}) = z^{2n}$.
\item[\rm b)] $C_{2k}(\la_{2n}) = \delta_{kn}$ $($the Kronecker
delta\/$)$.  In particular
$C_{2n}(\la_{2n}) = 1$ and so $\deg(C_{2n})\ge 2n$.
\end{enumerate}
\end{proposition}

\begin{proof} By definition $\nabla_K(\delta L_{2n}) = \sum_{S<L_{2n}}
(-1)^{s} \nabla_K(S)$.  Each $S$ is a union $\cup S_i$ of sublinks $S^i$ of
$L^i$ with $s_i\le2$ components. Since the $S^i$ lie in disjoint balls,
$\nabla_K(S) = \nabla_K(S^1) \dots \nabla_K(S^n)$, and so
$\nabla_K(\delta L_{2n})$ is a sum of products, which can be rewritten as
the product of sums $\prod^n_{i=1}\sum_{S^i<L^i} (-1)^{s_i} \nabla_K(S^i) =
\prod^n_{i=1}\nabla_K(\delta L^i) = \bigl( \nabla_K(\delta L) \bigr)^n
= z^{2n}$ by Example \ref{ex:conway}.  This completes the proof of a), and
b) follows since $\nabla_{\la_{2n}} = \nabla_K(\delta L_{2n})$. 
\end{proof}

\begin{remark}\label{rem:grading} This proposition can also be proved by
expanding $\la_{2n}$ as a linear combination of manifolds, and then
evaluating $C_{2k}$.  This approach, although longer, facilitates the
computation of {\sl products\/} of Conway coefficients and can be used to
establish lower bounds for the ranks of the groups $\calg_{2n}(S^1\x S^2)$
(see \S5).  

We indicate how this is done.  Write $\tau$ for $0$-surgery on the
right-handed trefoil $T$, and more generally $\tau^n$ for $0$-surgery on a
connected sum of $n$ copies of $T$.  Then it is readily seen that
$\la_{2n} = (\tau-1)^n$, where the right hand side is expanded using
the binomial theorem and ``1" is to be interpreted as $S^1\times S^2$. 
Since $\nabla_{\tau^j} = (1+z^2)^j$, it follows that $C_{2k}(\tau^j)$ is
equal to the binomial coefficient
${j\choose k}$, and so
$$
C_{2k}(\la_{2n}) = \sum_{j=0}^n (-1)^{n-j} {n\choose j} {j\choose k}.
$$
Observe that in this formula, $k$ can be a multi-index
$(k_1,\dots,k_m)$, in which case $C_{2k} = \prod C_{2k_i}$ and
${j\choose k} = \prod{j\choose k_i}$.  If $m=1$ then this reduces to the
formula in \ref{prop:realization}b by a well known combinatorial
identity.
%
%
The case $m=n$ with $k=(1,\dots,1)$ gives the formula
$$
C_2^n(\la_{2n}) = \sum_{j=1}^n (-1)^{n-j} {n\choose j} j^n.
$$ 

In particular for $n=2$ we see that $(C_4,C_2^2)(\la_4) = (1,2)$.  A
similar calculation shows that $(C_4,C_2^2)(\hat\la_4) = (0,4)$ for
$\hat\la_4 = [\Sigma_K,\hat L_4] \in \calm_4(S^1\times S^2)$, where $\hat
L_4$ is the $4$-component ``circular link" obtained from $L_8$ by banding
together pairs of components, as shown in Figure 3.8.


\begin{center}
\begin{picture}(255,80)
\multiput(10,0)(60,0){4}{\usebox{\boro}}
\multiput(63.5,24.2)(0,2.6){2}
{\multiput(0,0)(60,0){3}{\line(1,0){15.5}}}
\multiput(11.9,24.3)(0,2.6){2}{\line(1,0){7.1}}
\put(12,37.1){\oval(20,20)[bl]}
\put(12,37.1){\oval(25,25)[bl]}
\multiput(-.8,36.9)(2.5,0){2}{\line(0,1){13}}
\put(12,50){\oval(20,20)[tl]}
\put(12,50){\oval(25,25)[tl]}
\multiput(12,59.6)(0,2.6){2}{\line(1,0){238}}
\multiput(243.6,24.2)(0,2.6){2}{\line(1,0){7}}
\put(250,37.1){\oval(20,20)[br]}
\put(250,37.1){\oval(25,25)[br]}
\multiput(259.8,36.9)(2.5,0){2}{\line(0,1){13}}
\put(250,50){\oval(20,20)[tr]}
\put(250,50){\oval(25,25)[tr]}
\thicklines
\put(10,43.5){\oval(40,52)[l]}
\put(250,43.5){\oval(40,52)[r]}
\put(10,69.5){\line(1,0){240}}
\multiput(10,17)(60,0){4}{\line(1,0){35}}
\put(-20,65){$K$}
\end{picture}
\end{center}

\bigskip
\centerline{\rm Figure 3.8: $\hat L_4$}
\medskip
\setcounter{theorem}{8}

\noindent It follows that $\calg_4(S^1\x S^2)$ has rank at least two,
detected by the degree $4$ linearly independent \fts\  $C_4$ and $C_2^2$. 
In \S5 it will be shown to have rank exactly two.
\end{remark}

We now return to the proof of the main theorem (3.1).

\begin{proof}[Proof that $\ref{thm:divisibility}$ and
$\ref{prop:realization}$ $\Rightarrow$ $3.1$:] Suppose $\betti(M)=1$ and
$L$ is a $(2n+1)$-component admissible link in $M$. To show that $C_{2n}$
is \fty\ of degree {\sl at most $2n$} it suffices to show that
$C_{2n}([M,L])=0$, that is that
$z^{2n+1}$ divides $\nabla_{[M,L]}$ (the latter is an abbreviation for
$\sum_{S<L}(-1)^s\nabla_{M_S}$). As mentioned above, $M=\Sigma_K$
for some rational homology sphere $\Sigma$ and some $0$-framed
null-homologous knot $K$ in $\Sigma$. By general position we may
assume $L\sbq\Sigma-K$. The epimorphism $H_1(\Sigma-K)\cong
H_1(M)\twoheadrightarrow\bz$ is given by linking number with $K$.
Since each component of $L$ is null-homologous in $M$, it must have
zero linking number with $K$. Thus $L$ is admissible in $(\Sigma,
K)$. Now $M_S=\Sigma_{S\cup K}=(\Sigma_S)_K$ so $\nabla_{M_S} =
\nabla_{K,\Sigma_S} = \nabla_K(S)$, by definition. Therefore
$\nabla_{[M,L]} = \sum_{S<L}(-1)^s\nabla_K(S)=\nabla_K(\delta L)$
which is divisible by $z^{2n+1}$ by 3.2. Hence $C_{2n}$ is \fty\ of
degree at most~$2n$, and so in fact of degree exactly $2n$ by
\ref{prop:realization}.
\end{proof}

It follows immediately from Theorem 3.1 and the previous proposition that
$\calg_{2n}$ is infinite for all $n$.

\begin{corollary}\label{cor:deg} The element $\la_{2n}$ $($in
$\ref{prop:realization}\/)$ is of infinite order in $\calg_{2n}(S^1\x
S^2)$.
\end{corollary}

\begin{proof} If $\la_{2n}$ or some non-zero multiple lay in
$\calm_{2n+1}$ then $C_{2n}(\la_{2n})$ would vanish by Theorem 3.1,
contradicting Proposition \ref{prop:realization}.
\end{proof}

More generally, if the knot $K$ of Figure 3.5 is replaced by an arbitrary
null-homologous knot $K^*$ in a rational homology sphere $\Sigma$,
with the link $L$ living in a small ball, then$\nabla_{K^*}(\delta
L)=\nabla_K(\delta L)\cdot\nabla_{K^*,\Sigma}= z^{2n}(1+\dots)$.  Thus we
have

\begin{corollary} For any $3$-manifold $M$ with $\betti(M)=1$ and any
$n\ge0$, the group $\calg_{2n}(M)$ is of positive rank.  Thus
$\calo_{2n}(M)$, the group of rational valued finite type invariants on
$\calm(M)$ of degree at most $2n$, has rank greater than $n$.
\end{corollary}

\begin{proof} Any such $M$ equals $\Sigma_{K^*}$ for some $0$-framed
null-homologous knot $K^*$ in a rational homology sphere $\Sigma$. The
construction of $L$ above yields a $2n$-component link such that
$\nabla_{[M,L]}=\nabla_{K^*}(\delta L)=z^{2n}\,+$ higher order terms
so $C_{2n}([M,L])=1$.  Thus $C_{2n}$ is of infinite order in
$\calo_{2n}(M)$.  The last statement follows since $\calo_{2n} =
\calg_0\oplus\cdots\oplus\calg_{2n}$.
\end{proof}

In fact much larger bounds for the ranks of these groups can be
deduced from the algebraic independence of the Conway polynomial
coefficients (as functions on the set of knots in $S^3$).

\begin{corollary}\label{cor:rank} Suppose $\betti(M) = 1$.  Then the Conway
invariants freely generate a polynomial algebra $P[C_2,C_4,\dots]$ in
$\calo(M)$.\foot{Coefficients are in $\bq$, but can be taken in $\bz$ if
$H_1(M)$ is torsion free.}  Therefore the rank of $\calo_{2n}(M)$ is
at least $p(0) + \cdots + p(n)$, where $p(k)$ is the
number of unordered partitions of $k$.
\end{corollary}

\begin{proof}  Assume to the contrary that there is a non-zero
rational polynomial $p(x_1,\cdots,x_m)$ such that $p(C_2,\dots,C_{2m})$ is
identically zero on $\calm(M)$.  Since $p\ne0$, there exist integers $n_i$
for which $p(n_1,\dots,n_m) \ne 0$.  Let $K$ be a knot in $S^3$ whose
Conway polynomial is $1 + n_1z^2 + \cdots + n_m z^{2m}$; it is well known
that such knots exist.  

Now recall that $M$ can be described as $0$-framed surgery on a suitable
null-homologous knot $J$ in a rational homology sphere $\Sigma$.  Moreover
all such manifolds, for varying $J$, are $H_1$-bordant since any Seifert
surface for $J$ can be ``unknotted" by $\pm1$-framed surgeries on small
circles that link the bands of the surface.  In particular, the manifold
$M_0$ obtained by
$0$-surgery on $K$ in $\Sigma$ (i.e.\ put $K$ inside a small ball in
$\Sigma$) lies in $\calm(M)$.  But $p(C_2,\dots,C_{2m})(M_0) =
p(n_1,\dots,n_m) \ne 0$, a contradiction.

Finally observe that for every $k$, the degree $2k$ part of
$P[C_2,C_4,\dots]$ lies in $\calo_{2k}(M)$, by Proposition
\ref{prop:algebra}, and is of rank $p(k)$.  The stated bound on
$\rk(\calo_{2n}(M))$ follows.
\end{proof}

\medskip
\noindent{\bf Remark.} It is not being claimed in \ref{cor:rank} that the
grading on $P[C_2,C_4,\dots]$ is preserved under its embedding in
$\calo(M)$.  Showing this would require more work.  However Remark
\ref{rem:grading} establishes this for the elements of degree $4$ or less,
i.e.\ any non-trivial linear combination of $C_4$ and $C_2^2$ is of degree
$4$.
\medskip

We now proceed with the proof of Theorem \ref{thm:divisibility}, which
will be based on the following result.

\begin{theorem}\label{thm:smooth} Suppose $\Sigma$, $K$ and $L$
are as in the hypothesis of $\ref{thm:divisibility}$ with $\ell\ge1$.
Let $J$ be a component of $L$ and let $L'=L-J$. Then there exist oriented
links $K_i$ in $\Sigma-L'$ and signs $\e_i = \pm1$ 
such that $L'$ is admissible in $(\Sigma,K_i)$ for each $i$, and 
$$
\nabla_K(S) - \nabla_K(S\cup J) = z \sum \e_i\nabla_{K_i}(S)
$$ 
for every sublink $S$ of $L'$.
\end{theorem}

To understand this theorem, the reader should think of the simplest
case when $J$ bounds an embedded disk in $\Sigma$ which is punctured
twice by $K$ and not at all by $L'$. Then the difference between
performing $\pm1$ surgery on $J$ or not doing so is a local ``crossing
change'' of $K$. If we let $K_0$ denote the usual ``smoothing" of $K$ then
$\nabla_K(S\cup J)-\nabla_K(S) = \e_0z \nabla_{K_0}(S)$ where $\e_0$ is
the framing on $J$, and clearly $L'$ remains admissible in $(\Sigma,K_0)$.
In general $J$ might be knotted and might have a more complicated
interaction with $K$ and $L'$. Thus the strategy of the proof is to show
that the general case reduces to this simple case, and that the effect on
the Conway polynomial of surgery on $J$ is to add or subtract terms of
the form $z$ times the Conway polynomial of a smoothing.  It is crucial,
however, that these smoothings $K_i$ (as well as the signs $\e_i$) be {\sl
independent of $S$}. By this we mean that $K_i$ is disjoint from $L$ so
that for any sublink $S$ of $L$ we may use the symbol $K_i$ to denote the
image of this single link in $\Sigma_S$.

\begin{proof}[Proof that $\ref{thm:smooth} \Rightarrow
\ref{thm:divisibility}$]  We induct on $\ell$, assuming $\ell\ge1$ since
the case $\ell=0$ is trivial. Choose a component $J$ of $L$ and set
$L'=L-J$.  Then $\nabla_K(\delta L) = \sum_{S<L'} (-1)^s \bigl(\nabla_K
(S) - \nabla_K(S\cup J)\bigr) = z \sum_{S<L'} (-1)^s \sum
\e_i\nabla_{K_i}(S)$ by \ref{thm:smooth}. Reversing the order of
summation, using that $\e_i$ and $K_i$ are independent of
$S$, this gives $z\sum^r_{i=1} \e_i \nabla_{K_i} (\delta L')$, and by
induction each $\nabla_{K_i}(\delta L')$ is divisible by $z^{\ell-1}$.
Hence $\nabla_K(L)$ is divisible by $z^\ell$.
\end{proof}

\begin{proof}[Proof of $\ref{thm:smooth}$] Let $\e_J$ denote the framing
of $J$.  A knot in $\Sigma-(K\cup L')$ will be called {\it simple\/}
if it bounds an embedded disk $D$ in $\Sigma - L'$ which intersects $K$
transversely in algebraically zero points.  Clearly $J'\cup L'$ is
admissible in $(\Sigma,K)$ if $J'$ is simple.
 
First assume that $J$ is simple. Then surgery on $J$ puts a full
$(-\e_J)$-twist in all the strands of $K$ passing through $D$ -- this can
be seen by ``blowing down" $J$ \cite{K}.  What results is an oriented link
$K'$ in $\Sigma-L'$ with $\nabla_{K'}(S) = \nabla_K(S\cup J)$ for all
$S<L'$.  This link can also be obtained from $K$ by a finite sequence of
crossing changes, which we assume have been specified.  Let $K^i$ be the
link obtained by changing the first $i$ crossings of $K$, and $K_i$ be
the link obtained from $K^i$ by smoothing the $i$th crossing.  Then 
$$
\nabla_K(S) - \nabla_K(S\cup J) = \sum \bigl(
\nabla_{K^{i-1}}(S) - \nabla_{K^i}(S) \bigr) = z\sum
\e_i\nabla_{K_i}(S)
$$
where $\e_i$ is the sign of the $i$th crossing ({\sl after\/} it is
changed).  Note that $L'$ is admissible in $(\Sigma,K_i)$ since
changing or smoothing a self-crossing of a link does not change its
linking numbers with other knots.

Now assume that $J$ is not simple.  We claim that there exists a
simple knot $J'$ with $d_J(S) = d_{J'}(S)$ for all $S<L'$, where by
definition $d_*(S) = \nabla_K(S) - \nabla_K(S\cup *)$.  The theorem
would then follow from the simple case.  

To establish the claim, we appeal to a well known fact about the behavior
of linking numbers under surgery (cf.\
\cite{Ho2}).

\begin{lemma}
Let $A$, $B$ be disjoint null-homologous knots in a rational homology
sphere $\Sigma$ and $J$ be a knot in $\Sigma-(A\cup B)$ with framing $\e_J
= \pm1$.  Then 
$$
\ell k_J(A,B) = \ell k(A,B) - \e_J \ell k(A,J)\ell k(J,B)
$$
where $\ell k$ and $\ell k_J$ denote linking numbers in $\Sigma$ and
$\Sigma_J$ respectively.
\end{lemma}

\begin{proof}
Set $\la = \ell k(A,B)$, $\la_J = \ell k_J(A,B)$, $\alpha = \ell k(A,J)$
and $\beta = \ell k(J,B)$.  Let $m_B,\ell_B$ be a meridian and longitude
of $B$ in $\Sigma$, and similarly define $m_J,\ell_J$.  Then $A$ is
homologous in $\Sigma-(B\cup J)$ to $\la m_B + \alpha m_J$.  But $m_J$
is homologous in the surgery torus to $-\e_J\ell_J$, and so $A$ is
homologous in $\Sigma_J-B$ to  $\la m_B - \e_J \alpha \ell_J =
(\la - \e_J \alpha\beta)m_B$.  Thus $\la_J = \la-\e_J\alpha\beta$.
\end{proof}

Using this result, it is easy to compare the Seifert form of $K$ (which
determines its Conway polynomial) in $\Sigma_S$ and $\Sigma_{S\cup J}$
as follows.  Choose a connected Seifert surface $F\sbq\Sigma-L$ for $K$
(it is often helpful to view $F$ as a disk with one-handles attached),
and for each sublink $S$ of $L'$, let $V_S$ denote the corresponding
Seifert form for $K$ in $\Sigma_S$.  In other words $V_S(a,b) = \ell
k_S(a,b^+)$ for $a,b\in H_1(F)$, where $\ell k_S$ denotes linking number
in $\Sigma_S$. Now consider the symmetric bilinear form 
$$
\Lambda_J:H_1(F)\x H_1(F) \to \bz
$$ 
sending $(a,b)$ to $\ell k(a,J) \ell k(J,b)$, where $\ell k$ is the
linking number in $\Sigma$.  We will call this the {\it linking form\/} of
$K$ associated to $J$.\foot{Note that this form is well defined,
independent of a choice of orientation on $J$.}  Then   
$$
V_{S\cup J} = V_S - \e_J \Lambda_J.
$$
Indeed the lemma applied to knots $A$ and $B$ representing $a$ and $b^+$
in $\Sigma_S$, for $a,b\in H_1(F)$, shows that 
$
V_{S\cup J}(a,b) = V_S(a,b) - \e_J \ell k_S(a,J) \ell k_S(J,b),
$
but linking numbers with $J$ in $\Sigma$ and $\Sigma_S$ coincide since $J$
bounds a surface in $\Sigma-S$ (or by repeated application of the lemma).

It follows that if $J'$ is any oriented knot in $\Sigma-(F\cup
L')$ which has the same framing and linking form as $J$ (the latter
holds for example if $J'$ has the same linking number as $J$ has with each
one-handle of $F$) and zero linking numbers with the components of
$L'$, then $d_J(S) = d_{J'}(S)$ for {\sl all\/} $S<L'$.  But it is
obvious that there exists such a knot $J'$ which is simple, chosen for
example to lie in a neighborhood of the zero-handle of $F$.  This
establishes the claim, and thus completes the proof of Theorem
\ref{thm:smooth}.
\end{proof}

We conclude this section with a conjectured generalization of Theorem
\ref{thm:divisibility} to links which can be used to study the
``Conway polynomials" of manifolds of higher first betti number (see \S8).

\begin{conjecture}\label{conj:div} If $K$ is a null-homologous oriented
$k$-component link with zero pairwise linking numbers in a rational
homology sphere $\Sigma$ and $L$ is an admissible link of $\ell$
components in $(\Sigma,K)$ then $z^{2k-2+\ell}$ divides
$\nabla_K([\Sigma,L])$.
\end{conjecture}

\noindent{\bf Remarks.} The case $\ell=0$ was recently proved by Levine
\cite{L2}. The case $k=1$ is covered by Theorem \ref{thm:divisibility}, and
the case $k=2$ follows from the methods of \S5 (the proof is sketched in
Remark \ref{rem:conjecture}).  Added in proof: The full conjecture has now
been established by Amy Lampazzi.


\section{Finite type invariants from quantum invariants}


In this section it is shown that the theory of \fts\ is highly
non-trivial, even for 3-manifolds with large first betti
number\foot{By contrast the \cite{LMO}
invariant, which provides a universal finite type invariant for homology
3-spheres \cite{Le}, gives quite restricted information for manifolds
with first betti number
$\betti>0$, and is in fact identically zero if $\betti>3$ \cite{H2}.}.
To accomplish this, we use the $\bz_{p^k}$-valued invariants $\tau^d_p$
introduced by the authors in \cite{CM}, that are extracted from
the quantum $\so(3)$-invariants. By studying these invariants as $p$ and
$d$ approach infinity, we establish the {\sl rational\/} non-triviality
of the theory and provide strong evidence that much of Ohtsuki's theory
$\calo(S^3)$ of \fts\ of homology 3-spheres embeds in $\calo(M)$ for any
$M$. In addition, it is shown that for arbitrarily high betti number,
the theory exhibits all of the complexity of \fts\ of homology spheres
which ``come from $\sll(2)$-weight systems'' --- namely Ohtsuki's
rational valued invariants of homology spheres.

Recall the {\it quantum invariants\/} $\tau^G_p$ of 3-manifolds
associated with a compact {\it gauge group\/} $G$ and a positive integer
{\it level\/} $p$.  They were first discovered in a physical context
by Witten \cite{W}, and developed mathematically by Reshetikhin and
Turaev for $G=\su(2)$ \cite{RT}, and by Kirby and Melvin for $G=\so(3)$
\cite{KM}. Following the notation of \cite{CM} (rather than \cite{KM}) we
will use the abbreviation $\tau_p$  for the $\so(3)$-invariant
$\tau^{\so(3)}_p$ (denoted $\tau_p'$ in \cite{KM}), which can be viewed
either as a function on
$\cals$ or as a {\sl linear\/} function on $\calm$.  This invariant is
defined for all odd levels $p$ and, when normalized as in our
discussion of the proof of Lemma \ref{lem:op} at the end of this section,
takes values in the cyclotomic field $Q_p = \bq(q)$ where $q$ is a fixed
primitive $p^{\rm{th}}$ root of unity.  In fact, Hitoshi Murakami
\cite{M2} has shown that for {\sl prime\/} $p$, it takes values in the
ring of integers $\Lambda_p = \bz[q]$ in $Q_p$ (see also \cite{MR}), and
so in this case we have a $\bz$-linear map
$$
\tau_p : \calm \to \Lambda_p.
$$
Furthermore, $\tau_p$ is an $\bz$-algebra homomorphism with respect to
the connected sum operation $\# : \calm \times \calm \to
\calm$ (the bilinear extension of the corresponding operation on
$\cals$), i.e.\ $\tau_p(x\#y) = \tau_p(x)\tau_p(y)$. 

Henceforth we assume that $p$ is an odd prime.  Then $\Lambda_p$ (as an
abelian group) is free on $h^j$ for $0\le j\le p-2$, where $h=q-1$, and
so any element $a\in\Lambda_p$ can be written uniquely as $a =
a_0 + a_1 h + \cdots a_{p-2} h^{p-2}$.  Consider the projection
$
\pi^{j+(k-1)(p-1)} : \Lambda_p\to\bz_{p^k},
$ 
for $0\le j\le p-2$ and $k\ge1$, which maps $a$ to $a_j \pmod
{p^k}$.  Clearly any $a\in\Lambda_p$ is determined by the sequence
$\pi^d(a)$ for $d\ge0$.  Now define
$$
\tau^d_p:\calm\to\bz_{p^k}
$$ 
to be the composition $\tau_p^d = \pi^d\circ\tau_p$.  Then the following
is obvious but stated for emphasis.

\begin{proposition} For any odd prime $p$, the sequence of invariants
$\tau^d_p$ for $d\ge0$ determines and is determined by the quantum
$\so(3)$-invariant $\tau_p$.
\end{proposition}

The main result of this section is:

\begin{theorem} \label{thm:ft} For any odd prime $p=2n+3$ and any
integer $d\ge0$, the closed oriented $3$-manifold invariant $\tau^d_p$ is
a \fti\ of degree at most $3d$, in fact of degree at most $3d-n\bettip(M)$
when restricted to $\calm(M)$, where $\bettip(M) = \rk(H_1(M;\bz_p))$. 
\end{theorem}

Before giving the proof, we discuss a number of applications.

It is known that the full quantum invariant $\tau_p$ is not of
finite type for $p>3$ \cite[\S4]{CM} (note that $\tau_3 \equiv 1$), 
but Theorem \ref{thm:ft} shows that it is nevertheless a limit of finite
type invariants in the same sense that an analytic function is the limit
of its Taylor polynomials.  The Conway and Jones polynomials for
knots are also of this nature.  If one pursues the analogy
that \fts\ are the ``polynomials'', then such limits of \fts\ should
be called ``analytic'' invariants.  

We make this more precise.  An
invariant $\phi:\calm \to A$ is {\it weakly analytic\/} if
$\phi(\calm_\infty)=0$.\foot{Thus the set
$\calo_\infty^A$ of
$A$-valued weakly analytic invariants is the dual space
$\ho(\calm/\calm_\infty,A)$, in analogy with the corresponding sets
$\calo_\ell^A = \ho(\calm/\calm_{\ell+1},A)$ of \fts.}  The reader
can check that this is equivalent to the statement that
$\phi$ is {\it dominated} by \fts, in the sense that any classes in
$\calm$ which can be distinguished by $\phi$ can be distinguished by a
\fti\ (namely one of the projections $\calm \to \calm/\calm_\ell$).

We say that $\phi$ is {\it analytic\/} if there is an
inverse system $\{A_k\}$ of abelian groups and finite type invariants
$\phi_k:\calm\to A_k$ such that $A\subset\varprojlim A_k$ and
$\pi_k\circ\phi = \phi_k$ for all $k$.  Here $\pi_k:A\to A_k$ are the
restrictions of the natural projections.  

Observe that finite type
$\Rightarrow$ analytic (take $A_k = A$ and $\phi_k = \phi$ for all $k$)
while the reverse implication fails; for example the projection $\calm
\to \calm/\calm_\infty$ is analytic but not of finite type (also see
below).  Similarly analytic $\Rightarrow$ weakly analytic (since
$x\in\calm_\infty \Rightarrow \pi_k\phi(x) = \phi_k(x) = 0$ for all $k$,
and so $\phi(x) = 0$) while the converse presumably fails (although we
do not know an example).

In this language, we have the following consequence of Theorem
\ref{thm:ft}, which seems to be new even for homology spheres.

\begin{corollary}\label{cor:analytic} If $p$ is an odd prime, then
$\tau_p$ is analytic, and therefore dominated by finite type invariants.
\end{corollary}

\begin{proof} Let $A=\Lambda_p$, $A_k=\oplus^{p-1}\bz_{p^k}$, $\phi =
\tau_p$ and $\phi_k = \oplus_{j=0}^{p-2}\tau_p^{j+k(p-1)}$. Then the
$\phi_k$ are of finite type (by Theorem \ref{thm:ft}), $\Lambda_p \cong
\oplus^{p-1}\bz\subset\varprojlim A_k = \oplus^{p-1} \bz_{(p)}$ (where
$\bz_{(p)}$ is the $p$-adic integers) and $\pi_k\circ\phi = \phi_k$ for
all $k$.  Thus $\tau_p$ is analytic.  
\end{proof}

As another consequence of Theorem \ref{thm:ft}, we have:

\begin{corollary}\label{cor:bordism} If $\rk\,H_1(M;\bz_p) \equiv 0$ mod
$3$ for some odd prime $p=2n+3$, then the invariant
$\tau^{n\bettip/3}_p$ is constant on the entire $H_1$-bordism class of
$M$.
\end{corollary}

\begin{proof} Degree zero invariants are constant on
the $H_1$-bordism classes.
\end{proof}

This is interesting since $H_1$-bordism is fairly well understood in
terms of triple cup products and linking forms \cite{CGO}. Therefore
it should be possible to calculate the precise topological meaning of
these invariants. For example among
manifolds with $H_1\cong\bz^3$, the invariant $\tau^n_p$ is
completely determined by its values on the family of manifolds
$M_k$ given by $0$-surgery on the links obtained from the
Borromean rings by cabling one component $(1,k)$ times, for $k\ge0$. 
(These manifolds represent all the $H_1$-bordism classes \cite{CGO}.) One
has the strong feeling that there should be a single integral invariant
which determines the $\tau^n_p$ for a fixed surgery equivalence class and
varying $p$. Lescop's invariant for $M_k$ is $k^2$ since it is given by
the coefficient of
$z^3$ in the Conway polynomial (\S5
\cite{Les}) (\S5
\cite{C1}).

Note that $\tau_p^n$ is not degree zero on $\calm(\#^2S^1\x S^2)$, since
it is zero for $\#^2S^1\x S^2$ but non-zero for zero surgery on a
Whitehead link \cite{CM}, and any two manifolds with
$H_1\cong\bz^2$ are $H_1$-bordant.

\medskip
We now head towards a proof of the main theorem (\ref{thm:ft}),
discussing along the way its applications to the study of the structure
of the filtered group $\calm$.  The proof we give follows from a
divisibility result for $\tau_p$ which extends the work of \cite{CM}. 
Our measure of divisibility is the {\it
$p$-order\/}
$$
\opp : \calm \to \bz \cup \{\infty\}
$$
defined by $\opp(x) = \val_h(\tau_p(x))$, where $\val_h$ is the
$h$-adic valuation on $\Lambda_p$.  Thus $\opp(x) = m$ if $\tau_p(x)$
is written as $c_mh^m + O(h^{m+1})$ with $(c_m,p) = 1$ (see \cite{CM}). 
Equivalently, $\opp(x)$ can be defined to be the minimum $d$ for
which $\tau^d_p(x)\ne0$, or the maximum $d$ for which $h^d$
divides $\tau_p(x)$ in $\Lambda_p$.   

Observe that $\opp(x)$ is infinite if and only if $\tau_p(x) = 0$,
and so it is only by means of elements of {\sl finite} $p$-order that
$\tau_p$ can be brought to bear on the study of the filtration of
$\calm$.

\begin{definition}\label{def:normal} An element $x$ in $\calm$ is {\it
normal\/} if $\opp(x)$ is finite (i.e.\ $\tau_p(x)$ is non-zero) for
arbitrarily large $p$.  Let $\caln$ denote the set of all normal elements,
and $\cala$ denote its complement, the set of all {\it abnormal}
elements.
\end{definition}

Evidently $\calm_\infty \subset \cala$.  (In fact the inclusion is
proper: the difference of any two manifolds with equal quantum invariants
clearly lies in $\cala$, but if carefully chosen can be shown not to lie
in $\calm_\infty$ \cite{CM2}.)  It is not known, however, whether there
exist any abnormal {\sl
manifolds\/}.\foot{i.e.\ manifolds with
$\tau_p = 0$ for all but finitely many $p$\,; manifolds with $\tau_p = 0$
for infinitely many $p$ are known to exist, for example $0$-surgery on
the trefoil \cite[\S5]{CM}.}  

The collection of normal manifolds includes examples with any prescribed
$H_1$ (e.g.\ connected sums of rational homology spheres with copies of
$S^1\times S^2$); it is conceivable that every $3$-manifolds is normal,
or at least $H_1$-bordant to a normal manifold.  For normal manifolds $M$
it will be seen that the filtration of $\calm(M)$ is very rich.

\begin{remark}\label{rem:op} The reader is warned that $\opp$ is
highly non-linear.  Indeed it follows from properties of valuations
and the multiplicativity of $\tau_p$ that 
\begin{enumerate}
\item[\rm{a)}] $\opp(x+y) \ge \min\{\opp(x),\opp(y)\}$
\item[\rm{b)}] $\opp(mx) = \opp(x) + \val_h(m) =
\opp(x) + (p-1)\val_p(m) $ \\ 
(where $\val_p$ is the $p$-adic valuation on $\bz$)
\item[\rm{c)}] $\opp(x\#y) = \opp(x) + \opp(y)$.
\end{enumerate}
\end{remark}

The mod\,$p$ first betti number $\bettip = \rk H_1(-,\bz_p)$ similarly
extends from $\cals$ to $\calm$ in a non-linear fashion by setting
$\bettip(\sum m_iM_i) = \min(\bettip(M_i))$. The main result of \cite{CM} gives a
lower bound for $\opp$ in terms of $\bettip$, namely
$$
3\opp(x) \ge n\bettip(x)
$$ 
for all $x\in\calm$, where $n=(p-3)/2$.  (See Theorem 4.3 in \cite{CM}
where this is proved for manifolds; the result extends to linear
combinations of manifolds by Remark \ref{rem:op} and the definition of
$\bettip$.) \ Here we refine this result, taking into account where $x$ lies
in the filtration of $\calm$.

\begin{lemma}\label{lem:op} {\rm ($p$-order bound)} If $x\in\calm_\ell$,
then
$
3\opp(x)\ge n\bettip(x)+\ell
$
for any odd prime $p = 2n+3$.
\end{lemma}

The proof of this lemma, which is quite technical, is
postponed until the end of the section.  Meanwhile we explore its many
consequences.  First observe that Theorem \ref{thm:ft} follows easily.  

\begin{proof}[Proof of $\ref{thm:ft}$]  If $x=M_{\delta L}$ where $L$ is a
link with $\ell>3d-n\bettip(x)$ components, then $\opp(x)>d$ by the
lemma, and so $\tau^d_p(x)=0$ by definition of $o_p$. Therefore
$\tau^d_p$ is \fty\ of degree at most $3d-n\bettip(M)$ on $\calm(M)$.
\end{proof}

We now wish to use these results to investigate the structure of the
filtered group $\calm$.  For conceptual reasons, it is convenient first
to reformulate Lemma \ref{lem:op}.  This lemma relates the $p$-order of
$x\in\calm$ to where $x$ lies in the filtration.  In particular, if we
define the {\it depth\/} of $x$ to be
$$
\dep(x) = \max\{\ell \, | \, x\in\calm_\ell\}
$$ 
(a non-negative integer or $\infty$), then the lemma can be viewed as
giving an upper bound for $\dep(x)$ based on information garnered from
$\tau_p(x)$.  This upper bound, called the {\it $p$-depth\/} of $x$, is
given by
$$
\dpp(x) = 3\opp(x)-n\bettip(x).
$$   
It should be thought of as a (quantum) measure of the depth of $x$, and
so $1/\dpp(x-y)$ is a measure of the difference between $x$ and $y$.

The basic properties of the $p$-depth function $\dpp:\calm \to \bz \cup
\{\infty\}$ are collected in the following lemma.  The first property
is just a restatement of Lemma \ref{lem:op}, and the last three follow
from Remark \ref{rem:op} and the definition and elementary properties of
$\bettip$.

\begin{lemma}\label{lem:pdepth} {\rm ($p$-depth properties)} For any odd
prime $p$ and $x,y \in \calm$,
\begin{enumerate}
\item[\rm{a)}] $\dpp(x) \ge \dep(x)$
\item[\rm{b)}] $\dpp(x+y) \ge \min\{\dpp(x),\dpp(y)\}$
\item[\rm{c)}] $\dpp(mx) = \dpp(x) + 3(p-1)\val_p(m)$ \ \ $($\text{for
any integer $m$}$)$
\item[\rm{d)}] $\dpp(x\#y) = \dpp(x) + \dpp(y)$. \qed
\end{enumerate}
\end{lemma}

Of particular interest are the elements in $\calm$ for which the
bound in Lemma \ref{lem:pdepth}.1 is sharp.

\begin{definition}\label{def:robust}
An element $x$ of finite depth in $\calm$ is {\it robust\/} if $\dpp(x)
= \dep(x)$ for all sufficiently large primes $p$ (and {\it strongly
robust\/} if this equality holds for all $p>3$).  In particular, a
manifold $M$ is robust if and only if $\dpp(M)=0$ for all large $p$. 
\end{definition}

Robust elements are clearly normal (\ref{def:normal}) but not
conversely (see below).  They enjoy a number of other special properties,
including the following.

\begin{proposition}\label{prop:robust} {\rm (properties of robust
elements)} 
\begin{enumerate}
\item[\rm{a)}] If $x$ and $y$ are robust, then $x\#y$ is
robust with $\dep(x\#y) = \dep(x) + \dep(y)$.
\item[\rm{b)}] If $M$ and $N$ are $H_1$-bordant $3$-manifolds, then $M$
is robust if and only if $N$ is robust.  Thus one may speak of robust or
nonrobust bordism classes.
\end{enumerate}
\end{proposition}

\begin{proof}
For a) we have $\dep(x\#y) \ge \dep(x) + \dpp(y) = \dpp(x)+\dpp(y)  = 
= \dpp(x\#y)$ (by \ref{lem:pdepth}.4).  Since $\dpp(x\#y) \ge
\dep(x\#y)$ for large $p$ (by \ref{lem:pdepth}a) this implies
$\dep(x\#y) =
\dpp(x\#y) = \dep(x) + \dep(y)$.  For 2) assume $M$ is robust, so
$\dpp(M)=0$.  But $\dpp(M) \ge \min(\dpp(M-N),\dpp(N))$ (by
\ref{lem:pdepth}b) and $\dpp(M-N) \ge 1$ (since $M$ and $N$ are
$H_1$-bordant) so $\dpp(N) = 0$.\foot{For a slightly different point of
view, one can prove b) using the invariant $\tau = \tau^{\opp(M)}_p$,
where $p$ is  chosen large enough so that $\dpp(M) = 0$.  Indeed
$\tau$ is constant by Corollary \ref{cor:bordism}.  Hence
$\tau(N) = \tau(M)\neq0$, and so $\opp(N)\le\opp(M)$.  Since $\bettip(M) =
\bettip(N)$, it follows that $\dpp(N) \le \dpp(M) = 0$ and so $\dpp(N)
= 0$.}
\end{proof}

\begin{example}\label{ex:mfd} A manifold $M$ is robust if and only if
$3\opp(M)=n\bettip(M)$ for all large $p$, and this forces the first
betti number $\betti(M)$ to be a multiple of $3$ (since $n=(p-3)/2$ is
not).  In fact all rational homology spheres (the case $\betti=0$) are
robust by a result of Murakami \cite{M2}, and it is well known that
the $3$-torus $T$ (with $\betti=3$) is robust (see e.g.\ 
\cite[\S5]{CM}).  It follows from \ref{prop:robust}a that for any
$b\equiv0\pmod3$ and any finite abelian group $A$, there is a robust
$3$-manifold with $H_1\cong\bz^b\x A$, obtained by connected summing
$b/3$ copies of $T$ with a suitable rational homology sphere.

On the other hand, the connected sum of manifolds one of which is
non-robust is itself non-robust, as the reader may easily check. Thus
for example $M_0 = \#^3 (S^1\x S^2)$ is not robust even though
$\betti(M_0)=3$. In fact, for manifolds with betti number $3$ and torsion
free homology, it is expected that the set of non-robust manifolds is
precisely the $H_1$-bordism class of this manifold.  The other bordism
classes are represented by the $3$-manifolds $M_k$ (for $k>0$) given by
$0$-surgery on the link obtained from the Borromean rings by performing a
$(1,k)$-cable on one component, and it has been confirmed that these are
robust classes at least for $k=1$ (since $M_1=T$) and $k=2$
\cite[\S5.4]{CM}.
\end{example}  

\begin{example}\label{ex:delta} An example of a (strongly) robust element
of positive depth is the difference 
$$
\Delta = S^3-P
$$ 
where $P$ is the Poincar\'e homology sphere.  To see this, recall that
$\Delta = S^3_{\delta L}$ where $L$ is $+1$ surgery on the Borromean
rings, and so $\dep(\Delta)\ge3$.  But Murakami has shown that
$\tau_p(\Delta) = -6\la(P)h+O(h^2)$, where $\la$  is Casson's invariant,
and so $\opp(\Delta) = 1$ for $p>3$.  Thus $\dpp(\Delta) = \dep(\Delta) =
3$ for {\sl all\/} $p>3$.  More generally, for each $k>0$ the connected
sum
$$
\Delta_k = \Delta \# \cdots \# \Delta \quad \text{($k$ copies)}
$$
is (strongly) robust of depth $3k$ by
Proposition
\ref{prop:robust}a.
\end{example}

We now return to the investigation of the filtration on $\calm$.
As an immediate consequence of Lemma \ref{lem:pdepth} we have the
following estimates for the orders of an element of finite $p$-depth
in the filtered quotients of $\calm$.

\begin{theorem}\label{thm:order} {\rm (order)}
Any $x\in\calm$ of finite $p$-depth $($i.e.\ $\tau_p(x)\ne0)$ has order
at least $p^r$ in $\calm/\calm_s$ for all $s > \dpp(x) + 3(p-1)(r-1)$. 
In particular
$x$ has infinite order in $\calm/\calm_\infty.$ 
Furthermore, if $x$ is robust of depth $d$, then it has infinite order in
the graded summand $\calg_d = \calm_d/\calm_{d+1}$.
\end{theorem}

\begin{proof}
Suppose that $mx = 0$ in $\calm/\calm_s$.  This means that $mx \in
\calm_s$ and so $s \le \dep(mx) \le \dpp(mx) = \dpp(x) +
3(p-1)\val_p(m)$ by properties a) and c) in Lemma \ref{lem:pdepth}.  This
leads to a contradiction unless $m$ is divisible by $p^r$.  The last
statement follows from the first by taking $r=1$ and $p\to\infty$.  
\end{proof}

From this theorem, it is apparent that non-triviality results for the
filtration on $\calm(M)$ will follow from the existence of suitable
elements of finite $p$-depth.  This existence is guaranteed, at least for
$M$ of finite $p$-depth, by the
following

\begin{theorem}\label{thm:existence} {\rm (existence)}
For any $3$-manifold $M$, there exist elements $x_k$ in $\calm_{3k}(M)$
for each positive $k$ such that  $\dpp(x_k) = \dpp(M) + 3k$ for every
prime $p>3$.  In particular the $x_k$ are $($strongly\/$)$ robust if $M$
is.
\end{theorem}

\begin{proof}
For $M=S^3$ the elements $\Delta_k$ constructed in Example
\ref{ex:delta} will do, and for general $M$, set $x_k = M\#\Delta_k$ and
apply Lemma \ref{lem:pdepth}d.
\end{proof}

One can now deduce a variety of non-triviality results for the filtered
group $\calm(M)$ under the mild (and perhaps vacuous) condition that $M$
--- or some manifold $H_1$-bordant to $M$ --- have finite $p$-depth for
some $p>3$.  At the least, one would hope that the filtration does not
stabilize, or equivalently that $(\calm_\ell/\calm_\infty)(M) \ne 0$ for
all $\ell\ge0$.  In fact it turns out that these groups are all of
positive rank (for $M$ as above), and in fact of infinite rank if $M$ is
normal (i.e.\ of finite $p$-depth for arbitrarily large $p$); this
establishes a kind of {\sl rational non-triviality} of the theory for
normal manifolds.  

One can also investigate how {\sl fast} the filtration descends, measured
by the sizes of the associated graded summands $\calg_\ell(M) =
(\calm_\ell / \calm_{\ell+1})(M)$, and more generally $(\calm_\ell /
\calm_{\ell+m})(M)$ for a fixed $m>0$.  The best results are
obtained for robust $M$, in which case the associated graded group
$\calg(M)$ is of infinite rank; this is a stronger form of rational
non-triviality establishing the strict descent of the filtration
over the rationals.  

These results are summarized in the following
 
\begin{corollary}\label{cor:nontriviality} {\rm (non-triviality)}
Let $M$ be a $3$-manifold of finite $p$-depth $($i.e.\ $\tau_p(M)\ne0)$
for some prime $p>3$.  Then: 
\begin{enumerate}
\item[\rm a)] For every positive integer $n$, there exists $m<\infty$
such that each $(\calm_\ell / \calm_{\ell+m})(M)$ has an element of
order at least $n$.
\item[\rm b)] Each $(\calm_\ell/\calm_\infty)(M)$ is of rank at
least $p-1$, and thus of infinite rank if $M$ is normal.  
\item[\rm c)] If $M$ is robust, then each $\calg_{3k}(M)$ has positive
rank, and so $\calg(M)$ and $\calo^A(M)$ $($with $A=\bz$ or
$\bq)$ are of infinite rank.\foot{To prove that $\rk(\calg(M))$ is
infinite, it is only necessary to assume $\dpp(M)$ is uniformly
bounded for infinitely many $p$, but we do not know any examples of
this which do not also satisfy the stronger condition of robustness.}
\end{enumerate}
\end{corollary}

\begin{proof}
For a), choose $r$ and $k$ with $p^r\ge n$ and $3k \ge \ell$.  Then the
element $x_k$ from Theorem \ref{thm:existence} lies in $\calm_{3k}(M)
\sbq \calm_{\ell}(M)$ and is of $p$-depth $\dpp(M) + 3k \ge \dpp(M)+
\ell$.  By Theorem \ref{thm:order}, $x_k$ has order at least $n$ in
$(\calm_\ell / \calm_s)(M)$ for any $s > \dpp(M) + \ell + 3(p-1)(r-1)$,
so any $m > \dpp(M) + 3(p-1)(r-1)$ will satisfy the required condition.

For b), it suffices to show that $x_\ell,\dots,x_{\ell + p-2}$
(provided by \ref{thm:existence}) are linearly independent in
$(\calm_\ell /
\calm_\infty)(M)$, or equivalently that any nontrivial integer linear
combination $c = \sum a_kx_k$ (summed over $\ell \le k \le \ell+p-2$)
does {\sl not\/} lie in $\calm_\infty(M)$.  Since $\tau_p$ is analytic
(\ref{cor:analytic}), it is enough to show that $\tau_p(c) = \sum
a_k\tau_p(x_k)$ is a non-zero element in the cyclotomic ring $\Lambda_p$. 

It can be assumed that the coefficients $a_k$ have no common factor. 
Choose the first one $a_m$ which is prime to $p$. Now observe that each
$x_k$ has $p$-order $k+n$, where $n = \opp(M)$, and so can be written in
the form $b_k h^{k+n} + O(h^{k+n+1})$ with
$b_k$ prime to $p$.  Since $p$ is divisible by $h^{p-1}$ in $\Lambda_p$, 
$\tau_p(c)$ can be written in the form $a_m b_m h^{m+n} +
O(h^{m+n+1})$.  Thus $\tau(c)$ has $p$-order $m+n$, since $a_mb_m$ is
prime to $p$, and so in particular is non-zero.

For c), note that $x_k$ is robust (by \ref{thm:existence}) and so of
infinite order in $\calg_{3k}(M)$ (by \ref{thm:order}).  Thus
$\rk(\calg_{3k}(M))>0$, and so $\calg(M) = \oplus\calg_{\ell}(M)$ and
$\calo^A(M) \cong \ \ensuremath{\mathrm{Hom}} (\calg(M),A)$ (since
$A=\bz$ or $\bq$) both have infinite rank.  
\end{proof}


In the preceding proof, a key role is played by the connected sum of
$M$ with elements in $\calm(S^3)$. There is a convenient way to
formalize this which sheds light on the relationship between the
theory of finite type invariants for homology spheres and the theory for
manifolds which are $H_1$-bordant to $M$.  Indeed, it will be shown below
that for ``most" $M$, this theory exhibits all of the complexity of
\fts\ of homology spheres which come from ``$\sll(2)$-weight systems'',
namely Ohtsuki's rational valued invariants $\la_0,\la_1,\la_2,\dots$
\cite{O1}.

For a fixed $3$-manifold $M$, consider the embedding 
$$
i:\calm(S^3) \hookrightarrow \calm(M)
$$ 
given by $i(\Sigma) = M\#\Sigma$.  Clearly $i$ respects the filtration
on $\calm$,\foot{This means that $i$ does not {\sl decrease\/} depth;
however in some instances $i$ may {\sl increase\/} depth.  For example
for $M=S^1\x S^2$, the depth of $i(2\Delta) = 2((S^1\x S^2) -
(S^1\x S^2)\#P)$ is at least $4$ (but no greater than $5$ by Lemma
\ref{lem:pdepth}), while $2\Delta$ has depth $3$.  Indeed it is shown in
\S5 that $\calm(S^1\x S^2)$ has {\sl no\/} (even) elements of depth
$3$.} and therefore induces a map
$$
i_*:(\calm/\calm_\infty)(S^3) \to (\calm/\calm_\infty)(M)
$$ 
and $A$-module maps 
$$
i^*:\calo^A(M)\to\calo^A(S^3)
$$ 
for each ring $A$.  Explicitly $i_*[x] =
[M\# x]$ (where $[x]$ denotes the coset $x+\calm_\infty$) and
$i^*(\phi)(x) = \phi(M\# x)$.  

It is an interesting (and presumably
difficult) problem to determine when $i_*$ is injective, and when
$i^*$ is surjective.  Injectivity of $i_*$ would mean that elements
of finite depth in $\calm(S^3)$ are never mapped to elements of infinite
depth in $\calm(M)$.  In particular if two homology spheres were
distinguished by some \fti\ (say with values in $A$) then some other
\fti\ (possibly with different values) would distinguish their connected
sums with $M$.  The surjectivity of $i^*$ would show that the latter
could be chosen with values in $A$.  Also, if surjectivity were known for
$A = \bz$ and all prime power cyclic groups, then the injectivity of
$i_*$ would follow.

Now observe that if $\tau_p(M)\ne0$, then $i$ maps elements of
finite $p$-depth in $\calm(S^3)$ to elements of finite $p$-depth (and
therefore finite depth) in $\calm(M)$ (by Lemma \ref{lem:pdepth}d), or
put differently, if a pair of (linear combinations of) homology spheres
can be distinguished by $\tau_p^d$ for some $d$ then so can their
connected sums with $M$, using a possibly larger choice for $d$.  It
follows that $\ker(i_*)$ lies in the set $\calq_p$ of all
classes in $(\calm/\calm_\infty)(S^3)$ of {\sl infinite $p$-depth\/},
that is
$$
\calq_p \ \equiv \ \{ [x] \, | \, \dpp(x) = \infty \},
$$
and this can be used to show that if $M$ is {\sl normal\/} then
$\ker(i_*)$ lies in the set $\calq$ of all classes of {\sl
infinite Ohtsuki depth\/},
$$
\calq  \ \equiv \ \{ [x] \, | \, \la_j(x) = 0 \text{ for all } j\ge0 \}.
$$
With a little more work, one can show (for suitable $M$) that 
im$(i^*)$ contains the subspace $\calo^p$ of $\bz_p$-valued homology
sphere invariants generated by the mod $p$ reductions of the first
$(p-1)/2$ Ohtsuki invariants,
$$
\calo^p \ \equiv \ \text{span}\{ \la_j \text{ mod } p \, | \, j =
0,\dots,n\}
$$ 
where $n = (p-3)/2$.  These results,
summarized below, provide evidence for the injectivity of $i_*$ and the
surjectivity of $i^*$.

\begin{corollary}\label{cor:istar}  Let $M$ be a $3$-manifold of finite
$p$-depth, and consider the maps $i_*$ and $i^*$ $($as above$)$ induced
by taking connected sums with $M$. Then:
\begin{enumerate}
\item[\rm a)] $\ker(i_*) \sbq \calq_p$, the set of classes of
infinite $p$-depth $($defined above$)$.    
\item[\rm b)] {\rm im}$(i^*) \spq \calo^p$ provided $M$ is of
minimal $p$-depth in its $H_1$-bordism class.
\item[\rm c)] If $M$ is normal then $\ker(i_*) \sbq \calq$, the set of
classes of infinite Ohtsuki depth $($defined above$)$.  In particular,
if $\Sigma_1$ and $\Sigma_2$ are homology spheres that can be
distinguished by the $($rational valued$)$ Ohtsuki invariants, then
$M\#\Sigma_1$ and $M\#\Sigma_2$ can be distinguished by the invariants
$\tau^d_p$ for some $p$.\foot{By contrast, the \cite{LMO} invariant,
which includes the Lescop invariant as its degree $1$ term, cannot
distinguish any
$M\#\Sigma_1$ from $M\#\Sigma_2$ if $\betti(M)$ is positive.}
\end{enumerate}
\end{corollary}

\begin{proof}
As remarked above a) is immediate from the additivity of $p$-depth (Lemma
\ref{lem:pdepth}d), and c) follows since $\calq \spq \cap \calq_p$
\ (where the intersection is over all $p$ for which $\tau_p(M)
\ne0$) \  when $M$ is normal.  To see this, recall that $\tau_p^d(x)
\equiv \la_d(x) \pmod p$ for large $p$ \cite{O1}.  Now if $[x] \in
\cap\calq_p$, then $\tau_p(x)=0$ for arbitrarily large $p$ (since
$M$ is normal) and so all the Ohtsuki invariants of $x$ vanish.  For the
last statement in c), consider the difference $\Sigma_1-\Sigma_2$.  

It remains to prove b).  Let $m=\opp(M)$, the $p$-order
of $M$.  Then $\opp(N) \ge m$ for every manifold $N\in\cals(M)$, the
bordism class of $M$, since $\bettip$ is constant on $\cals(M)$).  It
follows that $\tau_p(N)$ can be expressed {\sl uniquely} as a polynomial
$\sum_{j=0}^{p-2}c_{j}(N)h^{m+j}$ with integer coefficients.  Reducing
mod $p$ gives a family of invariants $$t^j:\cals(M)\to\bz_p$$ defined by
$t^j(N) = c_j(N) \pmod p$.  Observe that $t^j$ can be
identified with the invariant $\tau_p^{m+j}$ under the natural inclusion
$\bz_p \hookrightarrow \bz_{p^k}$ (where $k = \lfloor(m+j) /
(p-1)\rfloor+1$) and so is of \fty\ by Theorem \ref{thm:ft}. One specific
case is for $M=S^3$ and $m=0$, and then the $t^j$ are the just the mod
$p$ reductions of Ohtsuki's invariants $\la_j$ for $0\le j\le n$
\cite{O1}. Let us continue to use $\la_j$ to denote these so as to avoid
confusion. Then it suffices to show that $\{\la_j\}$ lie in the span of
$\{i^*t^j\}$ for $0\le j\le n$.

We compute $i^*(t^k)(x)=t^k(M\#x) = \sum^{p-2}_{j=0}t^j(M)\la_{k-j}(x)$.
Since $p$ and $M$ are fixed, the constants $c^j=t^j(M)$ satisfy
$i^*(t^k)=\sum^{p-2}_{j=0}c^j\la_{k-j}$ for $0\le k\le n$.
Since $\opp(M)=m$, the lowest order coefficient $c^0$ is
invertible in $\bz_p$. It follows that this system of equations can
be inverted, and so $\{\la_j\}$ lie in the span of $\{i^*t^j\}$.
\end{proof}

The theory $\calo(M)$  of finite type invariants on certain $H_1$-bordism classes $\cals(M)$ also has
connections with theory of Vassiliev invariants of knots.  We illustrate
this for $M=S^1\x S^2$. Consider the set $\calk$ of isotopy classes of
knots in $S^3$ and the map $\calk\overset{\psi} {\longrightarrow}
\cals(S^1\x S^2)$ which sends a knot $K$ to the homology $S^1\x S^2$
obtained by performing $0$-surgery on  $K$. Composition with any
invariant of homology $S^1\x S^2$'s yields an (unoriented) knot
invariant. In fact we have:

\begin{proposition}  The map $\psi:\calk\to\cals(S^1\x S^2)$ given by
$0$-surgery induces an algebra homomorphism $$\psi^*:\calo_\ell(S^1\x
S^2)\to\calv_\ell$$ from \fts\ for homology $S^1\x S^2$'s to Vassiliev
invariants of degree at most $\ell$ (both with values in a fixed ring
$A$).
\end{proposition}

\begin{proof} Crossing changes on a knot $K$ may be
achieved by performing $\pm1$ surgery on circles (trivial in $S^3$)
which link $K$ zero times. The collection of $\ell+1$ ``crossing change
circles'' forms an admissible link in the $0$-surgered manifold.
\end{proof}

It is an interesting question to characterize the image of $\psi^*$.

\begin{proposition} The image of $\psi^*$ contains
all of the Vassiliev invariants arising from the coefficients of the
Conway polynomial.  Moreover, the $\bz_5$ invariants $\psi^*(\tau^d_5)$
distinguish the right and left-handed trefoil knots, and so the image of
$\psi^*$ is not just the algebra generated by the Conway coefficients.
\end{proposition}

\begin{proof} The first statement is obvious given
the definition of the Conway polynomial of a manifold as in
section~3. The second statement is a calculation done in \cite{KM}.
\end{proof}

We conclude with an application of the basic properties of robust
elements to show how to construct ``interesting" degree $3$ lifts of the
Casson-Walker invariant $\la$.

\begin{theorem}\label{thm:lift} Fix a ``base manifold" in each
robust $H_1$-bordism class of $3$-manifolds of positive first betti
number.  Then there exists a finite type invariant $\tl\la:\calm\to\bq$
of degree
$3$ which satisfies
\begin{enumerate}
\item[\rm a)] $\tl\la$ is a ``lift'' of the
Casson-Walker invariant, that is 
$
\tl\la(\Sigma) = \la(\Sigma)
$ 
for any rational homology sphere, and
\item[\rm b)] $\tl\la$ detects homology sphere summands in all
other robust $H_1$-bordism classes, that is
$
\tl\la(M\#\Sigma) = \la(\Sigma)
$
for each chosen base manifold $M$ and $($integral$)$ homology sphere
$\Sigma$.
\end{enumerate}
\end{theorem}

\begin{proof} Set $\tl\la = \la$ on all $H_1$-bordism classes of
rational homology spheres, and $\tl\la = 0$ on all non robust classes. 
Now consider a robust class of positive first betti number, with chosen
base manifold $M$.  It suffices to construct a map $\tl\la :
(\calm/\calm_4)(M)\otimes\bq \to \bq$ satisfying b).  To do this, we
choose a basis for $(\calm/\calm_4)(M) \otimes \bq \cong \oplus^3_{i=0}
(\calg_i(M)\otimes\bq)$ containing $M$ (which generates
$\calg_0$) and $M\#\Delta$ (which represents a non-zero
element in $\calg_3$ by \ref{thm:order}); here $\Delta$ is the robust
element $S^3-P$ in $\calm(S^3)$ of depth $3$ discussed in Example
\ref{ex:delta}, and so $M\#\Delta$ is also robust of depth $3$ by
\ref{prop:robust}a.  Now define $\tl\la(M\#\Delta) = -1$, and
$\tl\la = 0$ on all other basis elements (including $M$).  Then
$\tl\la(M\#\Sigma) = \tl\la(M\#(\Sigma-S^3))$ for any integral homology
sphere $\Sigma$.  But $\Sigma-S^3$ is known to be of depth at least $3$,
and in fact $\Sigma-S^3 = \la(\Sigma)\cdot(P-S^3) = -\la(\Sigma)\Delta$
in $\calg_3$ \cite{O}.  Hence $\tl\la(M\#\Sigma) =
-\la(\Sigma)\tl\la(M\#\Delta) = \la(\Sigma)$ as desired.
\end{proof}

We now return to the key result:

\medskip
\noindent{\bf Lemma \ref{lem:op}.} {\rm ($p$-order bound)} {\it If
$x\in\calm_\ell$, then
$
3\opp(x) \ge n\bettip(x)+\ell
$
for any odd prime $p=2n+3$.}
\medskip

Before giving the proof, it is useful to review the definition of
the quantum $\so(3)$ invariant $\tau_p$.  Recall from \cite{KM} the {\it
$p$-bracket} $\langle L \rangle = \sum [k]J_{L,k}$ of a framed link $L$ in
$S^3$, a certain linear combination of colored Jones polynomials which
is invariant under ``handle-slides" \cite{K}.  It is a priori an integral
Laurent polynomials in an indeterminant $t$, but is to be viewed
as an element of the cyclotomic ring $\bz(q)$ (where $q$ is a primitive
$p^{\rm{th}}$ root of unity) by identifying $t$ with $q^{4^*}$ where
$4^*$ is any mod $p$ inverse of $4$.  The $p$-bracket can also be written
in terms of Ohtsuki's version $\phi$ of the Jones polynomial as
$$
\langle L \rangle = \sum_{c=0}^n(a|c)\phi_{L^c}
$$
(see Proposition 1.5 in \cite{CM}).  Here $a = (a_1,...,a_\ell)$ is a
multi-index of integers recording the framings of the components of $L$,
$c = (c_1,...,c_\ell)$ is a multi-index cabling for $L$ with associated
cable $L^c$, obtained by replacing each component $L_i$ of $L$ with $c_i$
zero-framed push-offs, and the sum is over all cables with $0\le c_i\le
n$.  The reader is referred to \cite{CM} for the precise definition of
$\phi$ and the coefficients $(a|c) = \prod_{i=1}^{\ell}(a_i|c_i)$, which
are all to be viewed as elements of $\Lambda_p$.     

Now to obtain a $3$-manifold invariant, one must normalize the
$p$-bracket to make it invariant under ``blow-ups" \cite{K}.  This is
achieved by dividing by a factor which depends only on the linking
matrix of $L$.  In fact there is some flexibility in the choice of this
factor according to what properties one wishes the quantum invariant to
have.  The most common choice is $b_{+1}^{\ell_+} b_{-1}^{\ell_-}
b_{0}^{\ell_0/2}$, where $b_a$ is the $p$-bracket of the $a$-framed
unknot, $\ell_+$ and $\ell_-$ are the number of positive and negative
eigenvalues of the linking matrix of $L$, and $\ell_0$ is its nullity
(or equivalently the first betti number of $S^3_L$). 
This leads to an invariant $\tau_p'$ which is {\it multiplicative} under
connected sums and {\it involutive} (with respect to $t\mapsto \bar t =
t^{-1}$) under orientation reversal \cite{KM}.  However because of the
square root $b_0^{1/2}$ this invariant does not in general take values in
$\Lambda_p$ but rather in $\Lambda_{4p} = \Lambda_p[i]$ where $i^2=-1$,
and this obscures some of its number theoretic properties.  For the
present purposes it is more convenient to define the {\it $p$-norm} of
$L$ to be
$$
|L| = b_{+1}^{\ell_+} b_{-1}^{\ell_-} b_{0}^{\ell_0}/h^{n\ell_0}
$$ 
where $h=q-1=t^4-1$ (in contrast with \cite{CM} where $h=t-1$).  We will
need the fact that 
$$
|L|=(a|0) \quad \text{if $M$ is admissible.} 
\eqno(1)
$$ 
This is an easy consequence of the definitions in \cite{CM}.

Now set
$$
\tau_p(S^3_L) = \langle L \rangle / |L|.
$$
It is easily seen, using the well known fact that $b_0$ is a unit times
$h^{2n}$, that $|L|$ is an element of $\Lambda_p$.  In fact $|L|$ is a
divisor of $\langle L \rangle$ \cite{M2} \cite{MR} (see also \cite{CM}
where a stronger result is proved) and so $\tau_p$ takes values in
$\Lambda_p$.  Evidently
$\tau_p$ is multiplicative under connected sums, and with this
normalization $\tau_p(S^3) = 1$ and $\tau_p(S^2\times S^1) = h^n$. 
Unfortunately $\tau_p$ is no longer involutive; indeed $S^2\times S^1$
is amphicheiral, while $h^n \ne \bar h^n$ is not real.  (Note that
$\tau_p$ and $\tau_p'$ differ by a unit in $\Lambda_{4p}$.  In
particular they have the same $p$-order, cf.\ the discussion in
\cite{CM}.)

\begin{proof}[Proof of Lemma $\ref{lem:op}$]
First observe that it suffices to prove the result for generators
$M_{\delta L}$ ($= [M,L]$) where $L$ is an $\ell$-component
admissible link in $M$.  Indeed any $x\in\calm_\ell$ can be written as a
sum $\Sigma n_ix_i$ where $x_i=[M_i,L_i]$ and $L_i$ has $\ell$
components. Suppose that we proved the lemma for the generators $x_i$,
that is to say $3\opp(x_i)-n\bettip(x_i)\ge\ell$ for all $i$. Since
$\opp(x)$ is the minimum $d$ for which $\tau^d_p(x)\neq0$, some
$\tau^{\opp(x)}_p(x_i)\neq0$, which implies $\opp(x_i) \le
\opp(x)$ for some $i$. Hence $\dpp(x) \ge3 \opp(x_i)-n\bettip(x)$ for
some $i$. But $\bettip(x)\le \bettip(x_i)$ for all $i$ so $\dpp(x) \ge
3\opp(x_i)-n\bettip(x_i)\ge\ell$. It follows that $\dpp(x)\ge\dep(x)$.
So we may assume that $x = M_{\delta L}$. 

{\sl Case $1$}:  Suppose that $M=S^3_J$ for some {\it diagonal} framed
link $J$ (i.e.\ all pairwise linking numbers vanish).  Then $\bettip(M) =
j_p$, the number of components in the sublink $J_p$ of $J$ consisting of
all $J_i$ with framings $a_i$ divisible by $p$.  We must show that
$$
3\opp(S^3_{J\cup\delta L}) \ge nj_p+\ell.
\eqno(2)
$$

By definition $\opp(S^3_{J\cup\delta L})$ is the $p$-order of
\begin{eqnarray*}
\tau_p(S^3_{J\cup\delta L}) &=& \sum_{S<L} (-1)^s
\tau_p(S^3_{J\cup S}) \\
&=& \sum_{S<L} (-1)^s \sum_{c,c_{L-S}=0} (a_{J\cup S}|c_{J\cup S})
\phi_{(J\cup S)^{c_{J\cup S}}}/|J\cup S|
\end{eqnarray*}
where $a_T$ and $c_T$ denote the restrictions of (multi-index) framings
$a$ and cablings $c$ of $J\cup L$ to a sublink $T$ of $J\cup L$.  (Thus
the inner sum is over all cablings $c$ of $J\cup L$ with $c_{L-S}=0$, or
effectively cablings of $J\cup S$.)   But if $c_{L-S} = 0$, then
$(a_{J\cup S}|c_{J\cup S}) = (a|c)/(a_{L-S}|0) = (a|c)/|L-S|$, by (1). 
Substituting this into the last displayed expression gives
$$
\sum_{S<L} (-1)^s \sum_{c,c_{L-S}=0} (a|c) \phi_{(J\cup
S)^{c_{J\cup S}}}/|J\cup L|
\eqno(3)
$$
since clearly $|J\cup S||L-S| = |J\cup L|$.  Now this sum can be
rewritten as a sum over {\sl all} cablings $c$ of $J\cup L$,
$$
\sum_c (-1)^{\#c_L} \left(\sum_{k=0}^m (-1)^k {m\choose k}\right) (a|c)
\phi_{(J\cup L)^c}/|J\cup L|
$$
where  $\#c_L$ is the number of components of $L$ whose cabling index is
positive (the {\it support} of $c_L$) and $m=\ell-\#c_L$.  Indeed the
number of times $(J\cup L)^c$ occurs in (3) is computed by fixing $c$ and
counting how many $S$'s there are which contain the support of $c_L$,
and the number of such $S$'s with $\#c_L+k$ components is clearly
$m\choose k$.  Finally, noting that the inner sum of signed binomial
coefficients vanishes unless $m=0$ (i.e.\ $\ell=\#c_L$, whence
$c_L\ge1$) we have
$$
\tau_p(S^3_{J\cup\delta L}) =    
\sum_{c,c_L\ge1} (-1)^{\ell} (a|c)\phi_{(J\cup L)^c}/|J\cup L|.
\eqno(4)
$$

A lower bound for the $p$-order of $\tau_p(S^3_{J\cup\delta L})$ can
now be obtained easily from the results of \cite{CM}.  It is shown there
(Propositions 3.6 and 3.7) that $\opp(a|c) \ge n(j+j_p+\ell)-|c|-|c|_p$,
where $|c|=\sum c_i$ is the total number of cables of $c$, and $|c|_p$ is
the total number of cables of the sublink $J_p$ (of components of $J$
with framings divisible by $p$).  Also $\opp(\phi_{(J\cup L)^c}) \ge
4|c|/3$ (Theorem 3.5, which follows from a result of Kricker and Spence
\cite{KS}), and $\opp|J\cup L| = n(j+\ell)$ (Proposition 3.11).  Hence
any term in the sum (4) has order at least $nj_p+|c|/3-|c|_p$. 
This clearly achieves its minimum value when $c_{J_p}=n$, $c_{J-J_p}=0$
and $c_L=1$, and this value is then $nj_p+(nj_p+\ell)/3-nj_p =
(nj_p+\ell)/3$.  This proves (2).

{\sl Case $2$}:  Consider an arbitrary $M_{\delta L}$.  We must show
$3\opp(M_{\delta L}) \ge n\bettip(M)+\ell$.  By Corollary 2.3 of \cite{M2},
there exists a $\bz/p\bz$-homology sphere $\Sigma$ such that $M\#\Sigma$
can be obtained by surgery on a diagonal link, and so
$3\opp(M_{\delta L}\#\Sigma) \ge n\bettip(M)+\ell$ by the previous
case.  But $\opp$ is additive under connected sums, since $\tau_p$
is multiplicative, and the main theorem of \cite{M2} shows that
$\opp(\Sigma) = 0$.  Thus $\opp(M_{\delta L}) = \opp(M_{\delta
L}\#\Sigma)$ and the lemma is proved. 
\end{proof}


\section{Combinatorial structure of finite type invariants}



\def\bolink{\begin{picture}(0,0)
\setlength{\unitlength}{.8pt}
\put(0,0){\thicklines\hcirc}
\put(55,3){\hcirc}
\put(0,0){\lcirc}
\put(0,0){\vcirct}
\end{picture}}

\newsavebox{\st}
\savebox{\st}{\begin{picture}(0,0)
\multiput(0,-.4)(0,.8){2}{\thicklines\line(1,0){20}}
\end{picture}}

\newsavebox{\stp}
\savebox{\stp}{\begin{picture}(0,0)
\multiput(1.69,1.13)(-.56,.56){2}{\thicklines\line(1,1){14.1}}
\end{picture}}

\newsavebox{\stn}
\savebox{\stn}{\begin{picture}(0,0)
\multiput(1.69,-1.13)(-.56,-.56){2}{\thicklines\line(1,-1){14.1}}
\end{picture}}

\newsavebox{\Th}
\savebox{\Th}{\begin{picture}(0,0)
\multiput(2,16)(24,0){2}{\circle{4}}
\put(4,16){\line(1,0){20}}
\put(14,18){\oval(24,24)[t]}
\put(14,14){\oval(24,24)[b]}
\end{picture}}

\newsavebox{\Sq}
\savebox{\Sq}{\begin{picture}(0,0)
\multiput(17,17)(24,0){2}{\multiput(0,0)(0,24){2}{\circle{4}}}
\multiput(0,0)(58,0){2}{\multiput(0,0)(0,58){2}{\circle*{4}}}
\multiput(19,17)(0,24){2}{\line(1,0){20}}
\multiput(17,19)(24,0){2}{\line(0,1){20}}
\multiput(.2,.2)(41,41){2}{\usebox{\stp}}
\multiput(.2,58.1)(41,-41){2}{\usebox{\stn}}
\end{picture}}

\newsavebox{\Wh}
\savebox{\Wh}{\begin{picture}(0,0)
\multiput(0,16)(68,0){2}{\circle*{4}}
\multiput(22,16)(24,0){2}{\circle{4}}
\multiput(0,16)(48,0){2}{\usebox{\st}}
\put(34,18){\oval(24,24)[t]}
\put(34,14){\oval(24,24)[b]}
\end{picture}}

\def\wh{\begin{picture}(0,0)
\multiput(0,16)(68,0){2}{\circle*{4}}
\multiput(22,16)(24,0){2}{\circle{4}}
\put(0,16){\usebox{\st}}
\put(48,16){\line(1,0){20}}
\put(34,18){\oval(24,24)[t]}
\put(34,14){\oval(24,24)[b]}
\end{picture}}

\newsavebox{\Y}
\savebox{\Y}{\begin{picture}(0,0)
\setlength{\unitlength}{1.5pt}
\put(22,14){\circle{4}}
\thicklines
\put(0,14){\line(1,0){20}}
\put(23.5,15.5){\line(1,1){14}}
\put(23.5,12.5){\line(1,-1){14}}
\end{picture}}

\newsavebox{\y}
\savebox{\y}{\begin{picture}(0,0)
\put(22,14){\circle{4}}
\put(-2,14){\circle*{4}}
\multiput(38,-2)(0,32){2}{\circle*{4}}
\put(0,14){\line(1,0){20}}
\put(23.5,15.5){\line(1,1){14}}
\put(23.5,12.5){\line(1,-1){14}}
\end{picture}}

\newsavebox{\ncurl}
\savebox{\ncurl}{\begin{picture}(0,0)
\setlength{\unitlength}{1.5pt}
\put(22,14){\oval(20,20)[t]}
\put(22,14){\oval(20,20)[bl]}
\put(32,12){\vector(0,-1){0}}
\end{picture}}

\newsavebox{\pcurl}
\savebox{\pcurl}{\begin{picture}(0,0)
\setlength{\unitlength}{1.5pt}
\put(22,14){\oval(20,20)[t]}
\put(22,14){\oval(20,20)[bl]}
\put(24.5,3.7){\vector(1,0){0}}
\end{picture}}

\newsavebox{\tri}
\savebox{\tri}{\begin{picture}(0,0)
\setlength{\unitlength}{1.5pt}
\put(-2,14){\circle{4}}
\multiput(-4,16)(-2,2){3}{\circle*{.7}}
\multiput(-4,12)(-2,-2){3}{\circle*{.7}}
\end{picture}}

\newsavebox{\uni}
\savebox{\uni}{\begin{picture}(0,0)
\setlength{\unitlength}{1.5pt}
\put(-2,14){\circle*{4}}
\put(6,15){$^e$}
\end{picture}}

\def\figI{\begin{picture}(0,0)
\setlength{\unitlength}{1.3pt}
\multiput(16,0)(0,36){2}{\circle{4}}
\thicklines
\multiput(0,0)(18,0){2}{\multiput(0,0)(0,36){2}{\line(1,0){14}}}
\put(16,2){\line(0,1){32}}
\end{picture}}

\def\figH{\begin{picture}(0,0)
\setlength{\unitlength}{1.3pt}
\multiput(4,18)(24,0){2}{\circle{4}}
\thicklines
\put(6,18){\line(1,0){20}}
\put(0,16){\oval(8,36)[br]}
\put(0,20){\oval(8,36)[tr]}
\put(32,16){\oval(8,36)[bl]}
\put(32,20){\oval(8,36)[tl]}
\end{picture}}

\def\figh{\begin{picture}(0,0)
\multiput(0,0)(24,0){2}
{\multiput(0,0)(0,32){2}{\circle*{4}}}
\multiput(0,16)(24,0){2}{\circle{4}}
\multiput(0,1)(24,0){2}
{\multiput(0,0)(0,17){2}
{\multiput(-.4,0)(.8,0){2}{\thicklines\line(0,1){13}}}}
\put(2,16){\line(1,0){20}}
\multiput(-7,1)(36,0){2}{$^i$} 
\multiput(-6,20)(36,0){2}{$^j$} 
\end{picture}}

\def\figX{\begin{picture}(0,0)
\setlength{\unitlength}{1.3pt}
\multiput(4,0)(24,0){2}{\circle{4}}
\thicklines
\put(6,0){\line(1,0){20}}
\put(4.4,2.2){\line(2,3){22.7}}
\put(27.6,2){\line(-2,3){10}}
\put(4.4,36.2){\line(2,-3){10}}
\multiput(0,36)(27.3,0){2}{\line(1,0){4.1}}
\multiput(0,0)(30,0){2}{\line(1,0){2}}
\end{picture}}


In this section we describe an epimorphism from a finitely generated
group of {\it Feynman diagrams} (trivalent graphs/relations) to the
graded group $\calg_{\ell}(M)$. We then use this to evaluate a few examples
for small values of $\ell$. We show that for many $M$, the kernel of this
epimorphism is larger than one might naively predict based on the theory
for homology spheres \cite{GO3}, that is, there are relations in the group
of graphs which are not captured by the ``standard'' IHX and AS relations.

For each $m\ge0$, we describe a set $G^m$ of admissible abstract
graphs.  Feynman diagrams will be defined below as certain equivalence
classes of linear combinations of elements of $G^m$.

\begin{definition} An {\it $m$-admissible\/} graph $\Gamma$ is a finite
$1$-dimensional cell complex whose edge set is partitioned into the {\it
colored\/} edges $\calj = \calj_1\cup\cdots\cup\calj_m$ (where each
$\calj_i$ is nonempty with edges colored by the number $i$) and the {\it
white\/} edges $\calw$, and whose trivalent vertices are equipped with
a vertex {\it orientation\/} (an ordering of its incident edges up to
cyclic permutation), subject to the following conditions:
\begin{enumerate}
\item[a)] Each vertex is of valence $1$ or $3$.
\item[b)] Each edge has distinct vertices.
\item[c)] Each trivalent vertex is incident to at least one white
edge, and to at most one colored edge of any given color.
\item[d)] Each colored edge has at least one univalent vertex, and if it
has two such vertices (i.e.\ if it is {\it isolated\/}), then it is the
only edge of that color.  
\end{enumerate}
The edges with at least one univalent vertex will be called {\it
external\/}, while those with none will be called {\it internal}.  The
graph is said to be {\it closed} if all of its white edges are internal.  
\end{definition}

\begin{definition} Let $G^m$ be the set of all $m$-admissible graphs, and
$\cald^m$ be the free abelian group on $G^m$. The {\it degree\/} of
$\Gamma\in G^m$ is the number of white edges in $\Gamma$, that is, the
cardinality of $\calw$. Let $\cald^m_{\ell}$ be free abelian group on the
degree $\ell$ elements $G^m_{\ell}$ of $G^m$.  Note that $G^m_{\ell}$ is a
finite set.  Finally let $\calc^m_{\ell}$ denote the subgroup of
$\cald^m_{\ell}$ spanned by all closed graphs of degree $\ell$.
\end{definition}

Choose a base manifold $M$ in each $H_1$-bordism class and choose a
framed link description $M=S^3_J$ where $m$ (for $m$anifold)
denotes the number of components of $J$. Rational surgery framings are
allowed. We note in passing that $J$ may be chosen to be fairly simple.
For example, if $H_1(M)$ is torsion-free then $J$ can be chosen to be
$0$-framed and ``special'' (in the sense of 2.10) in that it can be
obtained from a trivial link by ``Borromean replacements'' \cite{CGO}. 
We define a map $\psi_J$ below and observe that the proof of 2.1
shows it is a surjection.

\begin{theorem} For any $($rationally\/$)$ framed $m$-component link $J$
for which $M=S^3_J$, as above, there is an associated epimorphism $\psi_J :
\cald^m_{\ell} \longrightarrow \calg_{\ell}(M)$.
\end{theorem}

\begin{proof}
For each $\Gamma\in G^m_{\ell}$, choose an immersion
$\Gamma \looparrowright D^2$ whose double points avoid vertices (for a
slight technical advantage we choose an over-crossing edge at each double
point) and such that each colored edge has one of its vertices on $\p
D^2$. Associate to this an unoriented tangle $T(\Gamma)$ in a $3$-ball
$B_1$ by the rules shown in Figure 5.4 (as in \cite{O}) in such a way
that each edge of $\Gamma$ corresponds to a single component of the
tangle with corresponding color when appropriate. This must be done in
such a way that the local orientations at the trivalent vertices can be
extended to a global orientation of the tangle. This explains the choice
5.4a or b.


\bigskip
\centerline{\kern10pt \figure{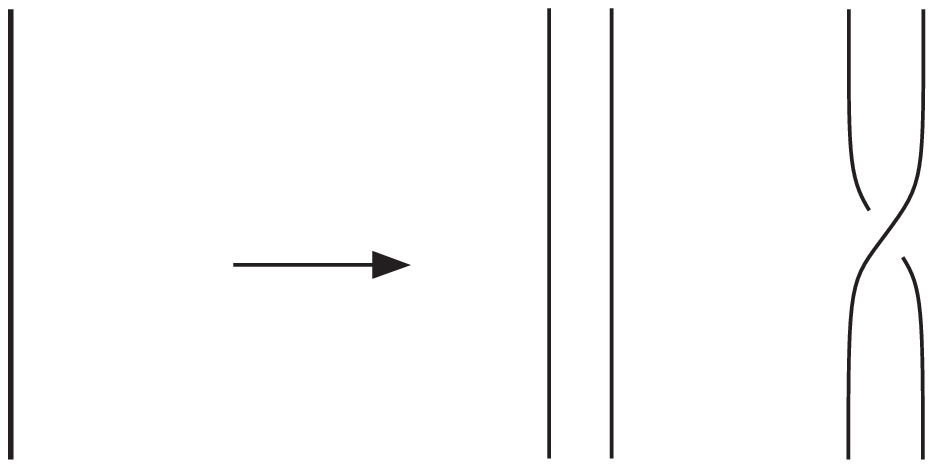}{.4} \kern60pt \figure{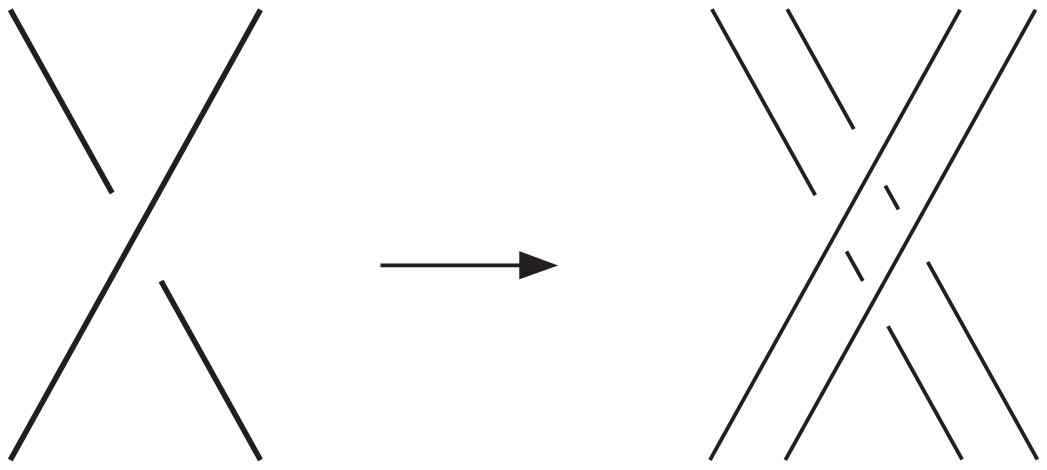}{.4}
\place{230}{-15}{a)} \place{198}{-15}{b)} \place{220}{20}{or}} 
\bigskip
\centerline{\figure{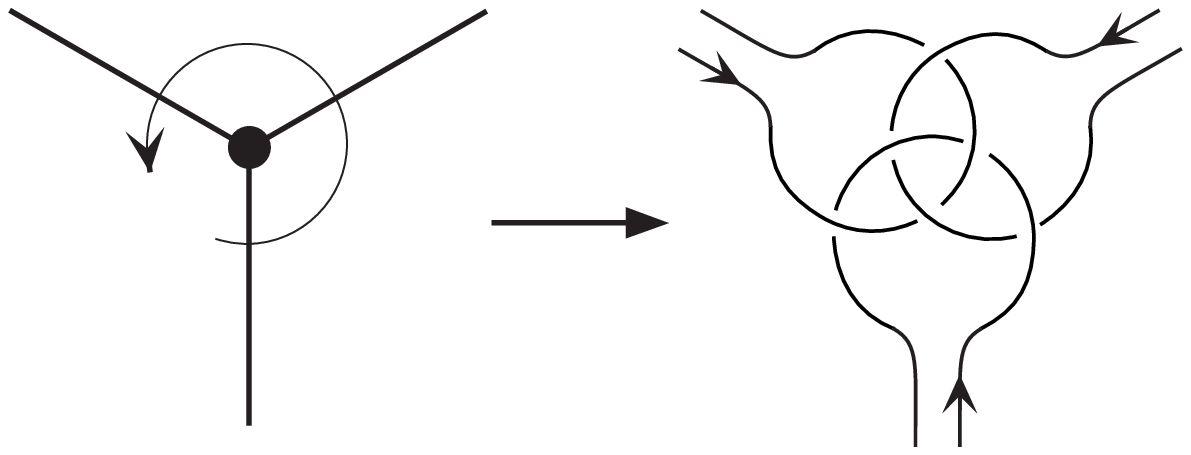}{.4} \kern40pt \figure{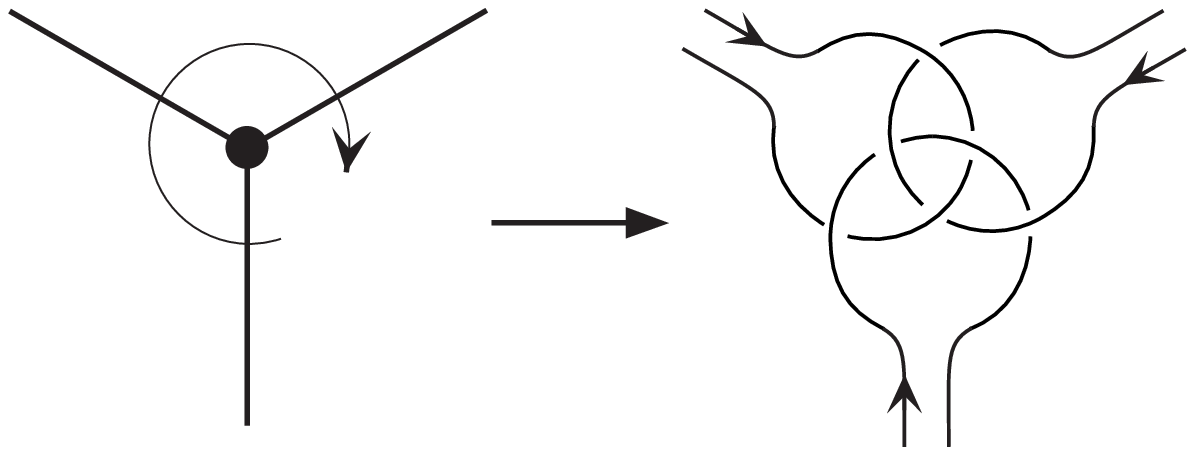}{.4}}
\centerline{Figure 5.4: $\Gamma \longrightarrow L(\Gamma)$}
\setcounter{theorem}{4}
\bigskip

Give each white component of $L(\Gamma)$ a $+1$ framing. Let $b_i$ be
the cardinality of $\calj_i$. Choose a $3$-ball $B_2$ in $S^3$ for which
the complementary tangle $(S^3-\intt B_2,(S^3-\intt B_2)\cap
J)$ is trivial and contains $b_i$ subarcs from the single
link component $J_i$. Then $(B_1,T(\Gamma))$ may be glued to
$(B_2,B_2\cap J)$ to form an unordered, unoriented framed link
$J\cup L(\Gamma)$ in $S^3$ which contains the link $J$ as sublink.
This gluing is not unique.

Now define
$
\psi_J:\cald^m_{\ell}\longrightarrow\calg_{\ell}(M)
$
to be the composition of the homomorphism $\cald^m_{\ell} \to
\calm_{\ell}(M)$, which sends $\Gamma$ to $M_{J\cup\de L(\Gamma)}$, with
the natural projection $\calm_\ell(M) \to \calg_\ell(M)$.  (Recall from
\S1 that $\de$ assigns to a framed link in $M$ the formal alternating sum
of its sublinks.)  It follows from the proof of Theorem~2.1 that
$\psi_J$ is surjective.

Observe that the map $\psi_J$ does not depend on the immersion of
$\Gamma$ since a ``band pass'' leads to equal elements in $\calg_{\ell}$
(cf.\ \cite{O}). For a similar reason it does not depend on the glueing
homeomorphism between $\p B_1$ and $\p B_2$ except for the information on
which components of $\calj_i$ are glued to which spots on $J_i$. If $J$
has zero linking numbers then even the latter does not matter (again by
the band-pass move or by the homotopy classification of links with zero
linking numbers by their $\ovmu(ijk)$). These statements will be
discussed more fully in \cite{CM2}. In any case, it may indeed be more
natural to average over all permutations of such glueings, but this
will not be needed in the present paper.
\end{proof}

Next we define a map
$$
d:\cald^m_{\ell} \longrightarrow \cald^m_{\ell}
$$ 
which is an extension of the ``deframing map'' of \cite{GO3}. For an
admissible graph $\Gamma$ and any subset $S$ of the set $T$ of all
trivalent vertices in $\Gamma$, let $\Gamma_S$ denote the admissible graph
obtained by ``splitting open''$\Gamma$ at each vertex in $S$ (creating $3s$
new univalent vertices) and deleting any resulting isolated colored edge
(unless it is the only edge with that color). Then set $d(\Gamma) =
\sum_{S<T}(-1)^s\Gamma_S$. Note that $d$ is the identity if $T$ is empty.

\begin{proposition} The deframing map $d$ is an isomorphism.
\end{proposition}

\begin{proof} The reader can verify that $d$ is its own inverse.
\end{proof}

In the remainder of this section we use the convention of \cite{GO3}
that a trivalent vertex of a graph $\Gamma$ lying the {\sl domain\/}
of the deframing map be denoted as in Figure~5.6a by a ``white
vertex,'' whereas for $\Gamma$ lying in the {\sl range} it will denoted
by a ``black vertex'' as in 5.6b.


\begin{center}
 \begin{picture}(175,50)
  \setlength{\unitlength}{1.5pt}
  \put(0,0){\usebox{\Y}}
  \put(85,0)
   {\begin{picture}(0,0)
    \put(0,0){\usebox{\Y}}
    \put(22,14){\circle*{4}}
    \end{picture}}
 \end{picture}
\end{center}

\medskip
\centerline{a) white vertex \kern 50pt b) black vertex}
\medskip
\centerline{Figure 5.6}
\setcounter{theorem}{6}
\medskip

We now identify five classes of relations on $\cald^m_{\ell}$ which lie in
the kernel of the composition of $\phi_J$ with the deframing map: AS
({\it antisymmetry\/}), S ({\it symmetry\/}), IHX, Y (an integrality
relation between Y-shaped graphs and closed graphs), and I
({\it isolated edge\/}).

\begin{theorem} The composition $\psi_J\circ d$ factors through an
epimorphism \penalty-1000$\phi_J:\cald^m_{\ell}/\{\mathrm{AS, S, IHX, I, Y}\}
\longrightarrow \calg_{\ell}(M)$ 
\end{theorem}

The relations AS, S, IHX, I and Y are defined in the proof.

\begin{definition} Let $\ov\cald^m_{\ell}\equiv\cald^m_{\ell} /
\{\mathrm{AS, S, IHX, I, Y}\}$.  The elements of $\ov\cald^m_{\ell}$ are
called {\it $m$-Feynman diagrams of degree $\ell$}.
\end{definition} 

\begin{proof}[Proof of $5.7$]
An element of $\eye$ is a graph $\Gamma$, one of whose white
edges is isolated. For such a graph we have $M_{J\cup\de L(\Gamma)}=0$
since $L(\Gamma)$ contains an isolated unknotted component. Since
$d(\eye)\sbq \eye$, it follows that $\psi_J\circ d(\eye)=0$.

The {\it antisymmetry relation\/} AS is shown in Figure 5.9 and
says that the effect of changing the vertex orientation at a single
trivalent vertex is the same as negation in $\cald$, as long as at least one
edge incident to that vertex is internal (i.e.\ ends in another trivalent
vertex). 


\def\pAS{\begin{picture}(0,0) \put(0,0){\usebox{\Y}}
\put(0,0){\usebox{\pcurl}} \put(0,0){\usebox{\tri}} \end{picture}}

\def\nAS{\begin{picture}(0,0) \put(0,0){\usebox{\Y}}
\put(0,0){\usebox{\ncurl}} \put(0,0){\usebox{\tri}} \end{picture}}

\begin{center}
 \begin{picture}(200,50)
  \setlength{\unitlength}{1.5pt}
  \put(0,0){\nAS}
  \put(60,12){$= \qquad -$}
  \put(110,0){\pAS}
 \end{picture}
\end{center}

\bigskip
\centerline{Figure 5.9:  Antisymmetry}
\medskip

\noindent
This is the same as Proposition~2.7 of \cite{GO3}, and the proof
that $\psi\circ d\mathrm{(AS)}=0$ also goes through as in \cite{GO3}, the
only essential ingredient being the half-twist lemma (\ref{lem:twist}).
Note that the ``marking lemma" (Lemma~2.1 of \cite{GO3}) also holds in the
present context, but since ``markings'' are not part of the structure of an
admissible graph (or a Chinese Character in the case of \cite{GO3}) it does
not directly indicate relations in $\cald^m_{\ell}$.

There are two types of {\it symmetry relations\/} S. The first is
shown in Figure~5.10 where $e$ is a {\sl white\/} edge of
$\Gamma$ with exactly one univalent vertex, and says that changing the
vertex orientation of the trivalent vertex of $e$ does not change the
image $\psi_J\circ d(\Gamma)$. The proof may be summarized as follows. A
change in vertex orientation leads to an insertion of an oppositely
oriented Borromean rings, changing a local $\ovmu(123)$ from $1$ to
$-1$, say. But the same effect on $\ovmu(123)$ can be achieved by
changing the orientation of the component arising from $e$. Since
these two are (locally) link homotopic, their images in $\calg_{\ell}$
are identical (see 2.9). But clearly the orientation of a link
component does not affect the surgered manifold.


\def\pS{\begin{picture}(0,0) \put(0,0){\usebox{\Y}}
\put(0,0){\usebox{\pcurl}} \put(0,0){\usebox{\uni}} \end{picture}}

\def\nS{\begin{picture}(0,0) \put(0,0){\usebox{\Y}}
\put(0,0){\usebox{\ncurl}} \put(0,0){\usebox{\uni}} \end{picture}}

\begin{center}
 \begin{picture}(190,50)
  \setlength{\unitlength}{1.5pt}
  \put(0,0){\nS}
  \put(60,12){$=$}
  \put(100,0){\pS}
 \end{picture}
\end{center}

\bigskip
\centerline{Figure 5.10:  Symmetry}
\medskip

The second type of symmetry relation is very similar and has an
identical proof. It states that, for any color $j$, changing the
vertex orientations at {\sl every\/} trivalent vertex which is incident to
an edge labelled by $j$ has no effect on $\psi_J\circ d(\Gamma)$.  This is
achieved by changing the orientation on the $j$-colored component of $J$.

The relation in Figure~5.11 is called the IHX {\it  relation\/} ---
assume clockwise vertex orientation in the plane of the picture (see
Figure~22 of \cite{GO3}). Note that any of the $4$ edges which leave the
picture can be colored or not colored. However, the $4$ edges leaving
the picture must be distinct edges, and no two may be colored alike.
This condition ensures that each of the $3$ graphs shown in 5.11 is
admissible. The proof of this set of relations is quite delicate and
will be postponed to \cite{CM2}. The case when none of the edges is
colored is due to Garoufalidis and Ohtsuki
\cite{GO3}.


\begin{center}
 \begin{picture}(200,60)
  \put(0,0){\figI}
  \put(60,22){$=$}
  \put(90,0){\figH}
  \put(145,22){$-$}
  \put(160,0){\figX}
 \end{picture}
\end{center}

\bigskip
\centerline{Figure 5.11:  The IHX Relation}
\medskip

The Y {\it relations\/} are shown in Figure 5.12, with the colored edges
drawn in thicker pen for clarity.  They are meant to say that if $\Gamma$
possesses any connected component which is Y-shaped, then $2\Gamma =
\Gamma'$ where $\Gamma'$ is obtained by replacing the Y-shaped component
(as shown) by the corresponding ``theta-shaped" closed graph\foot{Note
that the left hand side of each equation can be viewed as a half-theta \
\begin{picture}(6,6) \put(3,4){\oval(6,8)[l]} \put(0,4){\line(1,0){3}}
\end{picture} \ and the right hand side as a full theta \
\begin{picture}(6,8) \put(3,4){\oval(6,8)[l]} \put(3,4){\oval(6,8)[r]}
\put(0,4){\line(1,0){6}} \end{picture} \ with the colored edges (if any)
split open at the middle to conform to the definition of admissible
graphs.} with oppositely oriented trivalent vertices.  


\def\ya{\begin{picture}(0,0) \put(-15,11){$2$} \put(0,0){\usebox{\y}} 
\put(0,14){\usebox{\st}} \put(9,14){$^i$} \end{picture}}

\def\yb{\begin{picture}(0,0) \put(0,0){\usebox{\y}} \put(-15,11){$2$}
\end{picture}}

\def\yc{\begin{picture}(0,0) \put(-15,11){$2$} \put(0,0){\usebox{\y}} 
\put(0,14){\usebox{\st}} \put(9,14){$^i$} \put(22,14){\usebox{\stp}}
\put(24,21){$^j$} \end{picture}}

\begin{center}
 \begin{picture}(300,50)
  \put(0,0){\ya}
  \put(50,12){$=$}
  \put(70,-2){\usebox{\Wh}}
  \multiput(80,14)(47,0){2}{$^i$}
  \put(200,0){\yb}
  \put(250,12){$=$}
  \put(270,-2){\usebox{\Th}}
 \end{picture}
\end{center}

\smallskip
\centerline{a) \kern 150pt b)}\medskip

\begin{center}
 \begin{picture}(125,40)
  \put(0,0){\yc}
  \put(50,12){$=$}
  \put(80,0){\figh}
 \end{picture}
\end{center}

\smallskip
\centerline{c)}
\bigskip
\centerline{Figure 5.12:  Y Relations}
\setcounter{theorem}{12}
\medskip

\noindent
A sketch of the proof that $\psi_J\circ d=0$ for the
case 5.12c is as follows. Consider AS for {\sl one\/} of the
white vertices of the H-shaped graph on the right hand side of
the equation. Applying $\psi_J\circ d$ to this AS relation yields a
relation in $\calg_{\ell}$ wherein one sees two Borromean interactions of
opposite sign between the $i$, $j$ and white component. By link
homotopy considerations, as in \S2, these can be cancelled and
eliminated. The resulting relation in $\calg_{\ell}$ can then be seen to
be exactly $\psi_J\circ d$ applied to 5.12c. The other cases are
proved in exactly the same way. A more detailed proof will be included in
\cite{CM2}.

This completes the proof of Theorem 5.7 (modulo the IHX relations).
\end{proof}

Recall that $\calc^m_{\ell}$ is the subgroup of $\cald^m_{\ell}$ spanned by
closed graphs (all white edges are internal). One can speak of relations
AS, IHX and S among elements of $\calc^m_{\ell}$ since these relations
respect the defining condition for $\calc$. The following is then
immediate.

\begin{proposition} Let $\ov\calc^m_{\ell} = \calc^m_{\ell}/\{\mathrm{AS, S,
IHX}\}$. There is a commutative diagram of groups, as below, where the
horizontal maps are injective.
$$
\matrix
\ \ \calc^m_{\ell}\ \ &\hookrightarrow &\ \ \cald^m_{\ell}\\ \\
\downarrow & & \downarrow\\ \\
\ \ \ov\calc^m_{\ell}\ \ &\hookrightarrow &\ \ \ov\cald^m_{\ell}
\endmatrix
$$
\end{proposition}

One also has,

\begin{proposition} Let $\Gamma\in\cald^m_{\ell}$.  Then $2^\ell\ov\Gamma
\in \ov\calc^m_{\ell}$, where $\ov\Gamma$ denotes the equivalence class of
$\Gamma$ in $\ov\cald^m_{\ell}$, and so $\ov\calc \otimes \bz[\frac12] \cong
\ov\cald \otimes \bz[\frac12]$.  It follows that $\ov\calc^m_{\ell}$ is of
finite index in $\ov\cald^m_{\ell}$.
\end{proposition}

\begin{proof} Suppose $\Gamma$ has some external white edges. If any one of
these is {\sl not\/} part of a Y-shaped component, then, by AS and S (of
the first type), $2\ov\Gamma = 0$. On the other hand, if all of these
edges lie in Y-shaped components of $\Gamma$, then applying the $Y$
relations $k$ times (where
$k$ is the number of such components) shows that $2^k\ov\Gamma \in
\ov\calc^m_{\ell}$.  Clearly $k\le \ell$, and so the first statement
follows.  Since $\ov\cald^m_{\ell}$ is finitely generated, this implies
that
$\ov\calc^m_{\ell}$ is of finite index.
\end{proof}

\begin{corollary}\label{cor:cmn} The map
$\phi_J:\ov\calc^m_{\ell}\lra\calg_{\ell}(M)$ is an epimorphism after
tensoring with $\bz[\frac12]$ or $\bq$, and every element of the
cokernel of $\phi_J$ has order dividing $2^\ell$.
\end{corollary}

We shall see that, unlike the case of homology spheres, $\phi_J$ is
{\sl not\/} in general a rational isomorphism. In fact $\ov\calc^1_3$ has
rank one while $\calg_3(S^1\x S^2)$ has rank zero!

We compute some examples for the reader. Here $m=1$, $M=S^1\x S^2$, and
$J$ is the $0$-framed unknot in $S^3$. Recall $\calg_{\ell} =
\calm_{\ell}/\calm_{\ell+1}$. In the chart, $\bz_{5q}$ represents a non-zero
cyclic group of order a multiple of $5$ or $\infty$.

\bigskip


\centerline{\begin{tabular}{|c||c|c|c|c|c|c|} 
\hline \hline 
$\ell$ & $0$ & $1$ & $2$ & $3$ & $4$ & $5$  \\ 
\hline \hline 
$\ov\calc^1_{\ell}\ /\ 2$-torsion & $\bz$ & $0$ & $\bz$ & $\bz$ & $\bz^2$ 
&$\bz$ \\ 
\hline
generators & S & $-$ & W & $\Theta$ & C, W$*$W & W$*\Theta$ \\ 
\hline \hline
$\calg_{\ell}(S^1\x S^2)\ /\ 2$-torsion & $\bz$  & $0$ & $\bz$ & $0$ &
$\bz^2$ & $\bz_{5q}$ \\ 
\hline\hline
\end{tabular}}

\bigskip
\centerline{Figure 5.16: $\calg(S^1\x S^2)$ in low degrees}

\bigskip
Figure 5.17 shows pictures of the generators of $\ov\calc^1_{\ell}$ (mod
$2$-torsion). Since $m=1$, we do not need to label the colored components,
which are again shown in thicker pen. We shall briefly outline how the
table was derived. Let $\Gamma$ be an element of $\calc^1_{\ell}$ with $t$
trivalent vertices and $c$ non-isolated colored edges. Then it is easily
seen that $3t-c=2\ell$ by noting that two white edges emanate
from each of $c$ trivalent vertices while three emanate from each of the
other $(t-c)$ trivalent vertices, and that in this calculation each white
edge is counted twice.\foot{Note that the equation $3t-c=2\ell$ recovers
the result that, for homology spheres, $\calg_{\ell}\otimes\bq$ is zero
unless $\ell$ is a multiple of $3$ \cite{GL}\cite{GO3}.}  Hence $2\ell/3\le
t\le \ell$. This simplifies calculations, as does the following
observation.


\def\figS{\begin{picture}(0,0) \put(0,0){\usebox{\st}}
\multiput(0,0)(20,0){2}{\circle*{4}} \end{picture}}

\def\figW{\begin{picture}(0,0) \put(0,0){\usebox{\Wh}}
\end{picture}}

\def\figTh{\begin{picture}(0,0) \put(4,0){\figS}
\put(0,5){\usebox{\Th}} \end{picture}}

\def\figC{\begin{picture}(0,0) \put(0,0){\usebox{\Sq}}
\end{picture}}

\def\figWW{\begin{picture}(0,0) \put(0,0){\usebox{\Wh}}
\put(0,30){\usebox{\Wh}} \end{picture}}

\def\figWTh{\begin{picture}(0,0) \put(0,0){\usebox{\Wh}}
\put(20,30){\usebox{\Th}} \end{picture}}

\begin{center}
\begin{picture}(280,30)
\put(10,15){\figS}
\put(100,0){\figW}
\put(240,0){\figTh}
\end{picture}
\end{center}

\centerline{a) S \kern 85pt b) W \kern 95pt c) $\Theta$}

\begin{center}
\begin{picture}(330,60)
\put(10,0){\figC}
\put(120,-5){\figWW}
\put(250,-5){\figWTh}
\end{picture}
\end{center}

\centerline{d) C \kern 80pt e) W$*$W \kern 80pt f) W$*\Theta$}

\smallskip
\centerline{Figure 5.17}
\setcounter{theorem}{17}
\medskip

\begin{proposition} If $\Gamma\in\calc^m_{\ell}$ has an
odd number of trivalent vertices then $2\Gamma=0$ in $\ov\calc^m_{\ell}$.
More generally, if the number of non-isolated edges of some fixed color
$j$ is odd then $2\Gamma=0$.
\end{proposition}

\begin{proof} Let $c_i$ be the number of non-isolated $i$-colored
edges. The equation $3t-\Sigma c_i=2\ell$ derived above shows that if $t$ is
odd then some $c_j$ is odd. So it suffices to prove the second claim. Now
changing the vertex orientation at each of vertex incident to a
$j$-colored edge (denoted $\Gamma^*$) is a symmetry. On the other hand,
$\Gamma^*=(-1)^t\Gamma$ by anti-symmetry, since no component of $\Gamma$ is
Y-shaped. Hence, $2\Gamma=0$ in $\ov\calc^m_{\ell}$.
\end{proof}

Using the above considerations, one is led quite quickly by simple
combinatorics to see that $\ov\calc^1_{\ell}$ for $\ell\le4$ is {\sl
generated\/} by the graphs shown in the chart above. The case $\ell=5$
requires more work which we do not include here. It remains to show
that  W, $\Theta$, C and W$*$W are of infinite order (and linearly
independent) in $\ov\calc^1_{\ell}$. 

First consider the case $\ell=2$.  It was shown in \S3 that $\calg_2(S^1\x
S^2)$ has a map onto $\bz$ given by $C_2$, the coefficient of $z^2$ in the
Conway polynomial of the manifold. From Figure 5.12a we see that
W$=2\cdot$Y and then one calculates that $\phi_J(\text{Y})$ is $0$-surgery
on a trefoil knot minus $S^1\x S^2$. Hence $C_2(\phi_J(\text{Y}))=1$, and
the case $\ell=2$ is settled.

The case $\ell=3$ is the most interesting because here it will be seen that
$\phi_J$ has a non-trivial kernel. First we show that $\phi_J(\Theta)$
is zero by showing that $\phi_J$ of the graph Y$_1$ shown in
Figure~5.19a is $2$-torsion. We then apply the Y relation in Figure 5.12b
to see that $2$Y$_1=\Theta$ in $\ov\cald^1_3$.


\def\yst{\begin{picture}(0,0) \put(0,14){\figS}
\put(36,0){\usebox{\y}} \end{picture}}

\def\whst{\begin{picture}(0,0) \put(0,16){\usebox{\st}}
\put(0,0){\wh} \end{picture}}

\begin{center}
 \begin{picture}(210,50)
  \put(0,0){\yst}
  \put(140,0){\whst}
 \end{picture}
\end{center}

\centerline{a) Y$_1$ \kern 100pt b) Y$_2$}
\medskip
\centerline{Figure 5.19}
\medskip

Consider the framed links $L_1$ and $L_2$ in 5.20. These describe
homeomorphic $3$-manifolds as can be seen by ``sliding'' the smallest
$1$-framed circle over the $0$-framed circle.


\begin{center}
\begin{picture}(304,90)  
\setlength{\unitlength}{.8pt}
\put(0,0){\fcirc}
\put(-10,40){$1$}
\put(55,45){$0$}
\put(155,40){$1$}
\put(153,-10){$1$}
\put(60,0){\setlength{\unitlength}{.8pt}\Bor}
\put(240,0){\bolink}
\put(230,40){$0$}
\put(300,58){$1$}
\put(345,40){$1$}
\put(331,-10){$1$}
\end{picture}
\end{center}

\medskip
\centerline{a) $L_1$ \kern 150pt b) $L_2$}
\medskip
\centerline{Figure 5.20}
\setcounter{theorem}{20}
\medskip

\medskip\noindent
The reader can then work out that this implies that $\phi_J($Y$_1) =
-\phi_J($Y$_2)$, where Y$_2$ is the graph shown in 5.19b. But Y$_2$ is of
order~$2$ by an application of S and AS (see the proof of 5.14). Hence
we have shown that $\phi_J(\Theta)=0$.

To show that $\Theta$ is of infinite order, we use a little trick.
Observe that if $M=L(q,1)$ and $J'$ is the $q$-framed unknot then
$\phi_{J'}:\ov\calc^1_3\lra\calg_3(L(q,1))$ is a rational
epimorphism by Corollary \ref{cor:cmn}. So if $\calg_3(L(q,1))$ has
rank~$1$ then we are done. But this follows from 
\ref{cor:nontriviality}c. This is summarized as follows.

\begin{proposition} The map $\phi_J : \ov\calc^1_3 \lra \calg^1_3(S^1\x
S^2)$ is not a rational isomorphism. The graph denoted $\Theta$ in
Figure $5.17$ lies in the kernel. $($Here $J$ is the $0$-framed unknot$)$.
\end{proposition}

So the reader sees that more relations must be added to account for
handle slides. We shall not attempt a systematic treatment of this
in the present paper.

For the case $\ell = 4$, consider the image of W$*$W in
$\calg_4(S^1\x S^2)$. This is of infinite order as detected by $C_4$, the
coefficient of $z^4$ in the Conway polynomial; indeed it is represented by
the element $\la_4$ of Proposition \ref{prop:realization}.  Similarly 
$\phi_J($C) is the represented by the element $\hat\la_4$ introduced in
Remark \ref{rem:grading}, and is shown there to be of infinite order
(detected by $C_2^2$) and not a multiple of $\la_4$.  Therefore
$\calg_4(S^1\x S^2) = \bz^2$. 

Note that the the linear independence of C and W$*$W in $\ov\calc^1_4$
also follows from general principles, according to the following result.

\begin{theorem} Consider the set $\cala$ of all closed $m$-admissible degree
$\ell$ graphs with \,{\rm no} vertex orientations (for fixed $m$ and
$\ell$).  Let $\cale$ be the subset of $\cala$ consisting of graphs which
have an \,{\rm even} number of non-isolated edges of each color, and $\calo
= \cala - \cale$. Let $\calc(\cale)$ be the free abelian group on $\cale$
and $\calc(\calo)$ be group generated by $\calo$ with relations
$2\calo=0$. Then
$$
\ov\calc^m_{\ell}\cong\calc(\cale)/{\rm IHX}\oplus\calc(\calo)/{\rm IHX}
$$
where the IHX relations are as before, but restricted to the appropriate set
and with suitable sign changes $($see the proof\/$)$.
\end{theorem}

\begin{proof} We sketch a proof. Merely observe that the anti-symmetry
relations serve to eliminate generators and eliminate the vertex
orientations by choosing one for each abstract graph; one must of course
modify the signs in the IHX relations accordingly. The second symmetry
relation leads to a tautology if $\Gamma\in\cale$, or to $2\Gamma=0$ if
$\Gamma\in\calo$ (see Proposition 5.18).
\end{proof}

\begin{corollary} Consider the set $T$ of all $\Gamma\in\cale$, each of
which is a disjoint union of the closed ``theta-shaped" graphs that are the
right hand sides of the Y-relations $($Figure $5.12)$. Then $T$ is linearly
independent in $\ov\calc^m_{\ell}$. In particular, each such $\Gamma$ is of
infinite order.
\end{corollary}

\begin{proof} Note that $\langle{\rm IHX}\rangle\sbq\calc(\cale)$ is clearly
contained in the span of those $\Gamma$ which have some connected component
which either has $4$ different colors appearing, or has at least $3$
trivalent vertices. But the set $T$ is disjoint from this spanning set.
\end{proof}

This result can be refined to show C and W$*$W are linearly independent
in $\ov\calc^1_4$ by observing that C does not lie in the span of the IHX
relation since each embedding of an ``I-shaped graph'' in C has two
inputs colored alike. This was disallowed in IHX.

Observe that it follows from Corollary 5.23 that W$*\Theta$ is of infinite
order in $\ov\calc^1_5$. In fact $\phi_J($W$*\Theta)$ can be shown to be
non-trivial of either infinite order or order a multiple of $5$ in
$\calg_5(S^1\x S^2)$ by considering $\tau^2_5$ of section~4.


\section{Finite type invariants for spin manifolds}


The theory of invariants of \fty\ for closed spin
\m-manifolds was defined in 1.1--1.3 except for explaining how
the surgered $M_S$ inherits a spin structure from a spin structure on
$M$. The reader can compare the theory of N.~Shirokova \cite{Sh}. An
invariant of \fty\ for closed oriented \m-manifolds will be seen,
{\sl a fortiori\/}, to be an invariant of \fty\ for spin manifolds.
In addition the Rochlin invariant is a degree~$3\bmod16$ invariant of
\fty. The theory outlined by Shirokova in \cite{Sh} has neither of
these properties. As in \S2, we find that the group of invariants is
finitely generated within any fixed $H_1$-bordism class. In a
later paper we hope to investigate the mysterious invariants of spin
manifolds arising from quantum invariants as we have done in \S4 for
the non-spin invariants.

Here $\cals^{\spin}$ is the set of spin-structure-preserving
homeomorphism classes of spin \m-manifolds $(M,\sigma)$,
$\calm^{\spin}$ is the free abelian group on $\cals^{\spin}$, and
$\calm^{\spin}_{\ell}$ is the span of $[(M,\sigma),L]$ where $L$ is any
admissible link of $\ell$ components as in \S1. It is only necessary to
give a precise meaning to $[(M,\sigma),L]$ by assigning a spin
structure to the manifolds $M_S$ where $S<L$.

Given a spin manifold $M$ and an admissible link $S$, there is a
convenient way to specify the spin structure induced on $M_S$ using
the language of ``characteristic sublinks'' (see \cite{KM}; p.~541).
Namely, suppose $M=S^3_J$ and $J'\subseteq J$ is a characteristic
sublink corresponding to the given spin structure on $M$. Then the
appropriate spin structure on $M_S$ is the one corresponding to the
characteristic sublink $J'\cup S$. Note that since each component of
$S$ is $\pm1$-framed and has zero linking numbers with all other
components, $S$ {\sl must\/} be part of any characteristic sublink.
This ``framed surgery'' language is very convenient for checking
whether or not certain diffeomorphisms are actually spin
diffeomorphisms since most of the diffeomorphisms we employ are
described in terms of the ``Kirby calculus.'' 

If $A$ is a ring then $\calo^\spin$ is a filtered commutative
$A$-algebra (as shown in Proposition \ref{prop:algebra}). Since the
``forgetful map'' $\cals^{\spin}\to\cals$ respects the filtrations, the
following is clear.

\penalty-1000
\begin{proposition} If $\phi:\cals\to A$ is a \fti\ of degree
$\ell$ then $\phi':\cals^{\spin}\to\cals\to A$ (using the forgetful map) is
\fty\ of degree at most $\ell$, that is, there is a natural
monomorphism $\calo\hookrightarrow\calo^\spin$ which is an algebra
map.
\end{proposition}

Hence $\calo^\spin$ is large. There are also invariants not in the
subalgebra $\calo$.

\begin{proposition} The Rochlin invariant
$\mu:\cals^{\spin}\to\bz_{16}$ is a \fty\ degree $3$ invariant.
\end{proposition}

\begin{proof} Suppose $(M,\sigma)$ is a spin
\m-manifold. We claim that we may assume that $M$ is obtainable as
integral surgery on a link $J$ in $S^3$ which has all zero linking
numbers. For Murakami has shown that for any $M$ there exists a
connected sum of lens spaces $X$ such that $M\#X$ has such a surgery
description (\cite{M2}, Cor.~2.3). Moreover, if $L$ is not empty
$\mu([M,L])=\mu([M\#X,L])$ since the Rochlin invariant is additive
under connected sum and $[M\#X,L]$ is an alternating sum $[M,L]\#X$.
Thus we can assume $M=S^3_J$ as above.

Suppose $J'$ is the characteristic sublink of $J$ corresponding to
the spin structure $\sigma$ (see \cite{KM}; p.~541--544). Suppose
$L$ is an admissible link of $4$ components in $M$. By an
isotopy in $M$, we may assume $L$ lies in $S^3-J$ and has zero
linking numbers with each component of $J$. This uses the properties
of $J$ and the fact that each component of $L$ is null-homologous in
$M$. If $S<L$ then the characteristic sublink for the spin structure
on $M_S=S^3_{J\cup S}$ is $C_S=J'\cup S$, by definition. Recall that
the Rochlin invariant of $(S^3_{J\cup S},C_S)$ is given by
$\sigma(J'\cup S)-C_S\cdot C_S+8\text{Arf}(J'\cup S)\bmod16$
(\cite{KM}; p.~542). Here $\sigma$ is the signature of the linking
matrix and $\cdot$ is the total linking number. For brevity denote
this $\mu(M_S)$ by $\mu(S)$. We must show that
$\sum_{S<L}(-1)^s\mu(S)=0$, in other words that $\mu(\de L)=0$. Note
that $\sigma(J'\cup S)-C_S\cdot C_S=\sigma(J')+\sigma(S)-J'\cdot
J'-\tau(S)$ where $\tau$ is the trace of the linking matrix of $S$.
Since the latter matrix is diagonal with $\pm1$ entries on the
diagonal, $\sigma(S)=\tau(S)$. Thus $\sigma(J'\cup S)-C_S\cdot C_S$ is
{\sl independent\/} of $S$ and hence will not contribute to the
alternating sum. It remains to show that Arf~$(J'\cup\de
L)\equiv0\bmod2$ if $L$ has $4$ or more components. It has been shown
by Hoste, Murakami and Sturm that, for any ``totally proper'' link
$T$ in $S^3$, Arf~$(\de T)\equiv a_2(T)$, the coefficient of
$z^{t+1}$ in the Conway polynomial of $T$ \cite{Ho1}. Letting $T=\de
J'\cup L$ and using the fact that $\de\circ\de=\id$, we have
Arf~$(J'\cup\de L)\equiv$ Arf~$(\de\cdot\de J'\cup\de L)\equiv
a_2(\de J'\cup L)$. Now for any sublink $J''$ of $J'$, $J''\cup L$ is
an algebraically split link of more than $3$ components and Hoste has
shown that $a_2(J''\cup L)=0$ \cite{Ho2}. Hence Arf~$(J'\cup\de
L)\equiv0$ as desired. We remark in passing that J.~Levine's
generalization of Hoste's result has a proof which shows quite
clearly that $a_2\equiv0\bmod2$ if $J\cup L$ is algebraically split
$\bmod2$! (\cite{L2}, Proposition 4.1). Hence it is sufficient to
assume that $J$ is a ``totally proper'' link. Every \m-manifold is
surgery on a totally proper link in $S^3$ since any symmetric matrix
of integers can be diagonalized modulo~$2$ after stabilizing by
adding a $+1$.

Since $S^3-P$, where $P$ is the Poincar\'e homology sphere, lies in
$\calm^\spin_3$ and $\mu(S^3-P)\equiv8$, $\mu$ is of degree
precisely~$3$.
\end{proof}

\begin{theorem} For any closed spin \m-manifold $M$ and any integer
$\ell$, the group $\calg^\spin_\ell(M) = (\calm^\spin_\ell /
\calm^\spin_{\ell+1})(M)$ is finitely generated. Thus
$\calo^\spin_\ell(M)$ is finitely generated, and
$\calo^\spin_\ell=\Pi_{\calh^\spin}\calo^\spin_\ell(M_i)$ where
$\calh^\spin$ is the set of $H_1$-bordism classes of spin
\m-manifolds and $M_i$ is a representative from the class
$i\in\calh^\spin$.
\end{theorem}

\begin{proof} Lemma 2.2 remains true in the Spin category since it
is merely a combinatorial identity. Lemma~2.3 also holds using the
same proof. Lemma~2.4 remains true but the proof requires comment.
It is necessary to check that the diffeomorphism of the solid torus
used in the proof actually preserves the given spin structures. But
$S^1\x D^2$ has only two spin structures and these are determined by
looking at the spin structure on $S^1\x\p D^2$. Since the
diffeomorphism is the identity on the boundary, it preserves the
spin structure.

The ``Ohtsuki Lemmas'' 2.5 and 2.7 remain true. The only ingredients
of the proofs of 2.5 and 2.7 which are not definitions are the
diffeomorphisms associated to ``blowing up'' or ``blowing down'' an
unknotted circle which has zero linking numbers with all other
components. It must be checked that these diffeomorphisms preserve
the designated spin structures. Such $\pm1$ framed circles are
necessarily part of the characteristic sublink since they have zero
linking numbers with all other components, and for the same reason
it is known that blowing down such a curve does not change which of
the other components are in the characteristic sublink \cite{KM}.
For an identical reason, Lemma~2.9 remains true in the Spin
category. The rest of the proof of 2.1 works word for word,
reducing $\calg^\spin_\ell(M)$ to a finite spanning set which, indeed,
is obtained from the spanning set for $\calg_\ell(M)$ by including,
for each element $[M,L]$ of the latter, $[(M,\sigma),L]$ where
$\sigma$ varies over the $|H_1(M;\bz_2)|$ spin structures of $M$.
\end{proof}


\section{Finite type invariants for bounded manifolds}


We shall briefly discuss several theories for \fts\ for
compact \m-manifolds with boundary. The first theory leaves the boundary
``unmarked'' and the second and third assume the additional structure of
an orientation preserving homeomorphism $\phi:\p M\to S_g$ where $S_g$
is a fixed oriented surface in the homeomorphism class of $\p M$.
The first theory was defined in \S1 as the reader will note that no
assumption was made that $\p M$ is empty. In the second theory,
$\cals^\p$ is the set of triples $(M,\p M,\phi)$ as above where
$(M',\p M',\phi')\sim(M,\p M,\phi)$ if there is an orientation
preserving homomorphism $h:M\to M'$ such that $\phi'\circ h=\phi$ on
$\p M$. Given a link $L$ in $M$, a marking is induced on $\p M_L$ by
using the given product structure on the boundary of the cobordism
from $M$ to $M_L$. In the third theory, $\phi:\p M\to\p(H_g)$ ($H_g$
is the handlebody of genus~$g$) is required to induce $\phi_*:H_1(\p
M)\to H_1(\p H_g)$ which restricts to an isomorphism from the unique
$\bz^g$ summand containing kernel $(H_1(\p M) \overset{i_*}
{\longrightarrow}H_1(M))$ to the kernel of $H_1(\p H_g) \to H_1(H_g)$.

We deferred until now the proof of our ``Finiteness Theorem'' 2.1
for manifolds with boundary (unmarked). Let us indicate the changes
necessary in the proof given in \S2. The braiding and half-twist lemmas
need to be expanded to allow, in  Figures~2.6 and~2.8, ``pieces of the
boundary'' to run algebraically zero times through $L_1$. This is
made precise as follows. For each boundary component $S_{g_i}$ of
$M$, attach a handlebody $H_i$ with the same boundary to form a
closed oriented manifold $\what M$. Choose a spine for $H_i$ which is
abstractly homeomorphic to a union of $g_i$ circles, one base point
and $g_i$ arcs connecting the circles to the basepoint. Let the
image of this in $\what M$ be denoted $\what J_i$ and their union $\what
J$. As before $\what M$ can be expressed as surgery on a link $J$ in
$S^3$ which may be assumed to be disjoint from $\what J$. Hence $M$ is
recovered from $S^3_J$ by merely deleting a regular neighborhood of
$\what J$. $\what J_i$ should be viewed as a based $g_i$ component link
in $S^3$. Moreover if $L$ is an admissible link in $M$ then each
$L_i$ bounds a surface in $M$. Therefore we may assume that $L$ lies
in $\what M-\what J-J$ and that $L_i$ has zero linking number with each
component circle of $\what J$ (it bounds a surface in $\what M-\what J$),
as well as with each component of $J$ (as before). Now it is clear
that we have effectively changed a problem about manifolds with
boundary into a problem about closed manifolds with marked based
links $\what J$. Then Lemmas~2.5 and~2.7 remain true with ``strands''
of $\what J$ going through the disk spanned by $L_1$. Since $\what J$
merely records ``the location'' of $\p M$, this means these lemmas
hold in the category of manifolds with boundary. For the remainder
of the proof of Theorem~2.1 the reader should think of replacing the
link $J$ of the surgery lemma (\ref{lem:surgery}) and later by the
partially based link $J\cup\what J$. It is important to note that we needed
to choose a basing for our links in Definition~2.10 anyway, in order to use
Levine's work. Merely extend the partial basing to a full basing.
The rest of the proof of Theorem~2.1 now proceeds word for word with
$J\cup\what J$ replacing $J$.\hskip30pt $\square$

\medskip
Once again, invariants of degree~$0$ are precisely those functions
which are constant on surgery equivalence classes. These include
betti numbers, torsion numbers, the number of components of the
boundary, the genera of the boundary components, linking form
invariants, triple cup product forms and any invariants one might
choose to detect the isomorphism class of the {\sl pair\/}
$(H_1(M),H_1(\p M))$ (see \cite{CGO} for a fuller discussion).

We do not know if the second or third theories satisfy finite
generation.

Note that $\cals^\p\hookrightarrow\cals$ by ``plugging up'' $M$ via
solid handlebodies (using the marking). Hence
$\calo\hookrightarrow\calo^\p$, showing that $\calo^\p$ is large.


\section{Finite type invariants for marked manifolds}


Consider pairs $(M,\psi)$, where $M$ is a compact oriented \m-manifold
and $\psi$ is an isomorphism from $H_1(M)$ to a fixed abstract
abelian group $B$ (a ``marking'' of $H_1(M)$). Let $\cals^*$
be the set of equivalence classes of such pairs of marked
$3$-manifolds, where $(M_0,\psi_0)\sim(M_1,\psi_1)$ if and only if there
is an orientation-preserving homeomorphism $f:M_0\to M_1$ such that
$\psi_1\circ f_*=\psi_0$. Note that $(\#S^1\x S^2,\psi_0)\sim(\#S^1\x
S^2,\psi_1)$ for any $\psi_0$, $\psi_1$ so that if one is attempting to
distinguish $M$ from $\#S^1\x S^2$, there is no loss in marking $H_1$.
Now, if $S$ is an admissible link in $M$, then a marking of $H_1(M)$
extends naturally to a marking of $H_1(M_S)$, where $M_S$ is the surgered
manifold. Indeed it is clear that a marking of $H_1(M)$ extends over any
$H_1$-bordism. Thus there is a theory of \fts\ for this category (as
explained in section~1), which will be denoted by $\calo^*$. Note
that a theory based on pairs $(M,\alpha)$ where $\alpha\in
H^1(M;\bz_n)$ works similarly.

If $(M,\psi)$ is a marked \m-manifold then we can define many
group-valued invariants which would not be possible without the
marking. These include coefficients of the Conway polynomial,
Reidemeister torsion and Massey products (restricted to special
classes of manifolds so they are uniquely defined integers). Below we
shall show that the Conway coefficients are \fty. We shall not
address the Massey products here, although, since Massey products on
link exteriors are known to be of \fty, one must expect that they are
in this situation also. The extent to which Reidemeister torsion is
determined by \fts\ in this category will be detailed in a later
paper.

Suppose $(M,\psi)$ is a closed, marked \m-manifold with
$\betti(M)=m\ge1$. There is a canonical epimorphism
$B\twoheadrightarrow\bz$ given by sending each generator $1$ in each
$\bz$ factor of $B$ to $1$. The ``Alexander polynomial'' of
$(M,\psi)$ is the order of $H_1$ of the induced $\bz$-cover, divided
by $|{\rm torsion}\ H_1(M)|$. Any such manifold $M$ is $0$-framed
surgery on a link $K=\{K_1,\dots,K_k\}$ of null-homologous
components, with $\ell k(K_i,K_j)=0$, in a rational homology sphere
$\Sigma$. The Conway polynomial of $K$, $\nabla_K(z) =
z^{k-1}(a_0+a_2z^2+a_4z^4+\dots)$, is then defined and is related to the
Alexander polynomial of $\Sigma-K$ and hence to the Alexander polynomial
of $M$ in a similar fashion as explained in section~3 (see \S2.3.13 of
\cite{Les}). {\sl The Conway polynomial of $(M,\psi)$} is $\nabla_K(z)$.

\begin{theorem} Let $\cals^*$ be the set of equivalence classes of
closed marked \m-~manifolds $(M,\psi)$ with $\betti(M)=k\ge1$. Let
$C_\ell$ be the coefficient of $z^{k-1+\ell}$ in the Conway polynomial of
$(M,\psi)$. Then $C_\ell:\cals^*\to\bq$ is \fty\ of degree at most
$k-1+\ell$.
\end{theorem}

\medskip
\noindent{\bf Remark.} In fact if $\ell$ is odd then $C_\ell\equiv0$ so it
is degree $0$. If $\ell$ is even we claim the degree is precisely
$k-1+\ell$, but do not provide the proof here.

\begin{proof}[Proof of $8.1$] This follows immediately from Theorem
\ref{thm:divisibility}. The remark follows from Conjecture \ref{conj:div}.
\end{proof}

\begin{corollary} The Lescop invariant $\la_{\rm L}$ for (unmarked)
manifolds with $\betti=2$ is \fty\ of degree~1. The invariant
$\la_{\rm L}$ for manifolds with $\betti=3$ is \fty\ of
degree~0.
\end{corollary}

\begin{proof} $\la_{\rm L}$ equals $|{\rm torsion}\
H_1(M)|\cdot C_2(M)$ (\S5.1.6 of \cite{Les}). The corollary then follows
from Theorem 8.1 and the subsequent remark. The proof for
$\betti=3$ is easy and does not require 8.1 since in this case
$C_2$ is known to be the square of $\ovmu(123)$ \cite{C1} and this
is known to be constant on $H_1$-bordism classes (see section~1 and
also \cite{CGO}. Note that $\la_{\rm L}$ is independent of the
marking of $H_1(M)$.
\end{proof}

\begin{remark}\label{rem:conjecture}
Since we have invoked Conjecture \ref{conj:div} for $k=2$, $\ell=2$ in the
proof of 8.2 ($\betti=2$), we sketch the proof.
Theorem~\ref{thm:divisibility} guarantees that
$z^4$ divides $\nabla(\calm_4)$, whereas \ref{conj:div} claims $z^4$
divides
$\nabla(\calm_2)$ (restricted to $\betti=2$). Hence it suffices to
show $z^4$ divides the Conway polynomial of a generating set for
$\calg_2(\#^2_{i=1}S^1\x S^2)$ and $\calg_3$. Hence it suffices to
check this for the images of a generating set for the torsion free
part of $\ov\calc^2_3$ and $\ov\calc^2_4$, which is not difficult.
\end{remark}

For manifolds with $\betti = 0$, i.e.\ rational homology spheres,
Lescop's invariant agrees with the Casson-Walker invariant $\la$, which
is of degree $3$ (see Corollary 10.3 below).  Thus we have

\begin{corollary} The Lescop invariant $\la_{\rm L}:\cals\to\bq$ of
(unmarked) closed oriented \m-manifolds is \fty\ of degree~3.
\end{corollary}


\section{Further generalizations}


The theory we have presented is centered around the concept of
$H_1$-bordism. In effect, the \m-manifolds which are deemed
``close'' to $M$ are precisely those which are $H_1$-bordant to $M$
via a $4$-manifold $W$ which consists of a single $2$-handle
addition. The ``tangent vectors'' at $M$ to the ``space of
\m-manifolds'' are then the formal differences $\p_+W-\p_-W$, or
could even be thought of as the cobordisms themselves. This leads to
a theory in which the degree zero ``polynomials'' (being locally
constant on the space of \m-manifolds) are functions which are
constant on the $H_1$-bordism classes, which means they are
group-valued functions on the set of isomorphism classes of the
structure ($H_1$, linking form, triple cup product forms with
abelian coefficients). Hence our theory of \fts\ focusses on
distinguishing manifolds with isomorphic oriented cohomology rings,
separating this from the ``classical'' problem of distinguishing
cohomology rings.

There are additional ``classical'' invariants of \m-manifolds, namely
{\sl higher\/} Massey products, which could be included with the
cohomology rings, and there is a corresponding theory of \fts. We
summarize this theory below. Theories which fix even more aspects of
the homotopy type are possible but will not be discussed.

Let $k\ge2$ be an integer. We describe a family of theories of
$k$-\fts\ which agrees with our primary theory for $k=2$.

\begin{definition} A framed link $L$ in $M$ is called {\it
$k$-admissible\/} if
\begin{enumerate}
\item[a)] each component of $L$ lies in $(\pi_1(M))_k$, the $k^{\rm
th}$ term of the lower central series of $\pi_1(M)$
\item[b)] the pairwise linking numbers of $L$ are zero
\item[c)] the framings are $\pm1$.
\end{enumerate}
Clearly a sublink of a $k$-admissible link is itself $k$-admissible.
\end{definition}

\begin{definition} Let $\calm^k_\ell$ denote the subgroup of $\calm$
spanned by all $[M,L]$ where $L$ is a $k$-admissible link of
$\ell$ components in a $3$-manifold $M$.  A function $\phi:\cals\to A$ is
{\it $k$-\fty\  degree $\ell$\/} if $\phi(\calm^k_{\ell+1})=0$ and
$\phi(\calm^k_\ell)\neq0$, and $\calo^k_\ell = \ho(\calm/\calm_{\ell+1},A)$
is the algebra of all $k$-finite type invariants of degree at most $\ell$.
\end{definition}

Since $\calm^k_\ell \sbq \calm^{k-1}_\ell \sbq \dots \sbq \calm^2_\ell
\equiv \calm_\ell$ we have
$$
\calo^k_\ell \spq \calo^{k-1}_\ell \spq \dots \spq \calo^2_\ell \equiv
\calo_\ell,
$$ 
that is to say, there are {\sl more\/} invariants as $k$ increases.

\begin{definition} (see \cite{CGO}) Two $3$-manifolds $M$ and $N$
will be called {\it $k$-surgery equivalent\/} if there is a sequence
$M=M_0$, $M_1,\dots,M_r=N$ such that $M_{i+1}$ is obtained by
$\pm1$-surgery on a circle in $M_i$ which lies in ${\pi_1(M_i)}_k$.  They
are {\it $\pi/\pi_k$-bordant\/} if there is an oriented cobordism $W$
between $M$ and $N$, which is a ``product'' on $\pi_1/(\pi_1)_k$ (so for
$k=2$ this is
$H_1$-bordism).
\end{definition}

\begin{theorem} {\rm\cite{CGO}} Two $3$-manifolds $M$ and $N$ are
$k$-surgery equivalent if and only if $M$ and $N$ are
$\pi/\pi_k$-bordant $(k\ge2)$.
\end{theorem}

If one stipulates that the ``closest'' \m-manifolds to $M$ are ones
that are $\pi/\pi_k$-bordant via a single $2$-handle addition, and
that the tangent vectors at $M$ are formal differences of such, and
applies a notion of combinatorial derivative, then one generates
$\calo^k_\ell$ as the class of polynomials of degree at most
$\ell$.

\begin{proposition} Let $\calh_k$ denote the set of all $\pi/\pi_k$-bordism
classes of $3$-manifolds.  Then $\calm^k_\ell \cong
\bigoplus_{\alpha\in\calh_k} \calm^k_\ell(\alpha)$ and
$\calo^k_\ell \cong \Pi_{\alpha\in\calh_k} \calo^k_\ell(\alpha)$ where
$\calo^{(k)}_n(\alpha)$ is the corresponding theory restricted to
manifolds in the $\pi/\pi_k$-bordism class of $\alpha$.
\end{proposition}

It is shown in \cite{CGO} that $k$-surgery equivalence is related to
Massey products. It is shown that a manifold with $H_1\cong\bz^m$ is
$k$-surgery equivalent to $\#^m_{i=1}S^1\x S^2$ if and only if its
Massey products of order less than $2k-1$ vanish.

The proof that $\calo^k_\ell(\alpha)$ is finitely generated for
each $\alpha\in\calh_k$ is not complete even though almost all of
the steps of the proof of 2.1 carry over without difficulty.
Lemmas~1.4 and 2.2 hold without change, although a non-trivial
result from \cite{CGO} is required. Lemmas~2.3 and 2.4 hold with
$k$-admissible replacing admissible. Lemmas~2.5 and 2.7 hold without
alteration. Lemma~2.9 can be rephrased and partially recovered.

\begin{lemma} If $L$ and $L'$ are surgery equivalent links in a
$3$-manifold $M$ then $[M,L]\sim[M,L']$ in $\calg^k_\ell(M)$.
\end{lemma}

This is true because a surgery equivalence between {\sl links\/} in
$M$ is, {\sl by definition\/}, accomplished by a $\pm1$ surgery on a
circle $K$ which bounds a disk in $M$. Clearly more general
alterations are possible since $K$ could be allowed to represent a
non-trivial loop in $(\pi_1(M))_k$. Here the proof stops due to the
lack of an analogue of Levine's theorem. However note that it is
already possible to reduce to the case where the link $L\sbq L\cup
J\sbq S^3$ has only ``Borromean interactions'' and hence is given
by, loosely speaking, uni-trivalent graphs in $M$. This is entirely
consistent with the fact that $\pi/\pi_k$-bordism of manifolds is
classified by $H_3(\pi_1(M)/\pi_1(M)_k)$ modulo automorphism (see
\cite{CGO}). Since the latter group is finitely generated, it is
fairly clear that one can reduce to a finite set of {\sl
parameters\/} (presumably Massey products --- or Milnor's invariants ---
of weight less than $2k$). However the details have not yet been
considered. Moreover, it is less clear what is the analogue of the
final step (Lemma~2.5), that is to reduce from
$\ovmu(1122)=10$ $\ovmu(123)=6$, for example, to a sum of cases where
$\ovmu(1122)\in\{0,\pm1\}$ and $\ovmu(123)\in\{0,\pm1\}$.
Nonetheless it would be surprising if this was a serious problem.
Note that it is not necessary to {\sl classify\/} links modulo the
appropriate equivalence relation, just as it was not necessary for
us (in 2.1) to use the full strength of Levine's surgery equivalence
theorem. The ill-definedness of higher Massey products would be a
serious annoyance.

It seems clear, in light of recent work of Habegger and Masbaum
relating to Milnor's invariants to the Kontsevich integral, that
the $p$-order (see \ref{def:normal}) would vary less and less in a
$\pi/\pi_k$-bordism class as $k$ increases. This should allow for
the well-definedness of more invariants of $k$-\fty\ derived from
$\tau^{\so(3)}_p$.

The reader should note that $k$-\fty\ equals $2$-\fty\ for those
manifolds where $\pi_k=\pi_2$. This includes all manifolds with
cyclic first homology!

A theory based on control of {\sl all\/} the higher Massey products
at once seems attractive, but the finite generation (2.1) seems
unlikely for \m-manifolds whose lower central series strictly
descends.


\section{Relationships with other theories and other results}


In this section, we mention some relationships with other theories:
that of Garoufalidis-Ohtsuki \cite{GO1} for rational homology spheres, and
of Garoufalidis-Levine \cite{GL3} relating to the mapping class group.

The theory of Garoufalidis-Ohtsuki for rational homology spheres is based
on surgery on algebraically split links in homology spheres and as such
is not strongly related to our approach. In an attempt to get
$\calg_n$ finitely-generated they impose their ``Property~1'' which is
overly strong in our opinion. Morally, our theory should have
strictly more invariants. Certainly the $\bz_p$-rank of $H_1(M;\bz_p)$
is of \fty\ degree zero for us but not of \fty\ for them. However,
due to a slight flaw in their theory, we cannot show in generality
that an invariant which is of GO-\fty\ is \fty\ in our sense. Indeed,
Garoufalidis-Ohtsuki intended that $\calg_n$ should be
finitely-generated (consequence of their Theorem~2). However their
$\calg_0$ is not finitely generated: Suppose $M$ is a rational
homology sphere whose linking form is not isomorphic to the direct
sum of forms on cyclic groups (see \cite{KK}). Let $\phi$ be the
characteristic function on $M$. Then $\phi$ is \fty\ in the sense of
\cite{GO1}, because the only restrictions placed on $\phi$ by
\cite{GO1} involve Dehn surgery on algebraically split links in an
{\sl integral\/} homology sphere. But any manifold so obtained has a
linking form which is a direct sum of linking forms on {\sl cyclic\/}
groups (since its linking matrix is diagonal). Hence $\phi$ is zero
on all these manifolds. Since there are an infinite number of such
manifolds $M$ as above, their $\calg_0$ is infinitely generated.
(Indeed there are an infinite number of non-isomorphic linking forms
which are not ``diagonalizable''.)  But certainly $\phi$ is not \fty\
in our sense (for any $\ell$ there is a Brunnian $\ell$-component
link $L$ in $S^3$ on which surgery does not yield $S^3$ --- consider
$M\#[S^3,L]$).

Now we will show that, on the subclass of rational homology spheres,
any invariant which is \fty\ $n$ in the sense of \cite{GO1} and
which is additive on connected sums, is \fty\ of degree at most $n$ in our
sense.

\penalty-1000
\begin{theorem} Let $\calr\subset\calm$ be the span of the set
of rational homology spheres. Suppose that $\phi:\calr\otimes\bq\to\bq$ is
of \fty\ $n$ in the sense of Garoufalidis-Ohtsuki {\rm \cite[\S1.2]{GO1}}
and is additive on connected sums. Then the induced map $\phi:\calr\to\bq$
$($i.e.\ the composition of $\phi$ with the natural inclusion
$\calr\hookrightarrow\calr\otimes\bq)$ is \fty\ of degree at most $n$
in our sense.
\end{theorem}

\begin{corollary} The invariant of Casson-Walker for rational
homology \m-spheres is a rational valued \fti\ of degree~$3$.
\end{corollary}

\begin{proof}[Proof of $10.1$] In fact we need only assume that $\phi$
satisfies their ``Property $0$.''  Property $0$ says that
$\phi([\Sigma,L])=0$ for every {\sl integral\/} homology sphere $\Sigma$
and every rationally framed (with the proviso that the
framings be {\sl non-zero}) algebraically split link $L$ in $\Sigma$ with
more than $n$ components.  (Here ``algebraically split" means pairwise
linking numbers zero.)  Suppose $M$ is a fixed rational homology sphere
and $L$ is a fixed admissible $n+1$ component link in $M$. It suffices to
show that $\phi([M,L])=0$. Throughout we will identify $\calr$ with its
image in $\calr\otimes\bq$. 

First suppose that $M$ can be expressed as $S^3_J$ where $J$ is a
integrally framed algebraically split link in $S^3$. Then we have the
following combinatorial Lemma.

\begin{lemma} With the above notation,
$[S^3_J,L]=\sum_{S<J}(-1)^{s}[S^3,L\cup S]$.
\end{lemma}

The theorem follows immediately from the Lemma since, by Property~$0$
of \cite{GO1}, $\phi$ vanishes on $[S^3,L\cup S]$ since
$L\cup S$ has more than $n$~components. The Lemma is proved easily
by induction on $j$, the number of components of $J$. It is trivial for
$j=0$, so assume it for all links of $j\ge0$ components and consider a
link of $(j+1)$~components of the form $J\cup K$ where $K$ is the last
component. Then by Lemma 1.4, $[S^3_{J\cup K},L] = -[S^3_J,L\cup K] +
[S^3_J,L]$.  By induction this equals $\sum_{S<J}(-1)^{s} \bigl(
-[S^3,L\cup K\cup S] +[S^3,L\cup S] \bigr)$.  But this is
$\sum_{S<J\cup K}(-1)^{s}[S^3,L\cup S].$

Now consider the general case $[M,L]$. By a result of Murakami and
Ohtsuki \cite{M2}, there exists a rational homology sphere $X$ such that
$M\#X$ is integral surgery on some algebraically split link in $S^3$. But
$\phi([M\#X,L])=\phi([M,L])$ since $\phi$ is additive and
$L$ is not empty. Thus the above special case suffices to show that
$\phi$ is finite type.
\end{proof}

There is an interesting relation with the mapping class group. Recall
the subgroup $\calk$ of the mapping class group generated by Dehn
twists along bounding simple closed curves (see \cite{GL3}).

\begin{theorem} \textup{(\cite{CGO})} $M$ is $H_1$-bordant to $M'$
if and only if there is a Heegard splitting $M=H_1\cup_fH_2$ and a
homeomorphism $g\in\calk$ such that $M'=H_1\cup_{g\circ f}H_2$.
\end{theorem}

This indicates that one could filter all \m-manifolds using the type
of filtration discussed by Garoufalidis and Levine in (\cite{GL3}, 1.3)
corresponding to $\calk$, and that at least at the ``zero level'' it would
agree with our theory. However since Ohtsuki's theory for homology
spheres is a direct summand of our $\calm$, and since it is still unknown
even in this case if these theories agree (Ohtsuki versus \cite{GL3}), we
shall not pursue this here.

\bigskip
\begin{tabular}{lcl} 
{Tim D. Cochran}&\qquad\qquad&{Paul Melvin}\\
{Department of Mathematics}&\qquad\qquad&{Department of Mathematics}\\
{Rice University}&\qquad\qquad&{Bryn Mawr College}\\
{Houston, TX 77005--1892}&\qquad\qquad&{Bryn Mawr, PA 19010--2899}\\
{cochran@math.rice.edu}&\qquad\qquad&{pmelvin@brynmawr.edu}\\
\end{tabular}

\vfill

\end{document}